\documentclass{article}

\usepackage{amsfonts,amsmath,amssymb}
\usepackage{hyperref}
\usepackage{color}
\usepackage{graphicx}
\usepackage{subcaption}
\usepackage[section]{placeins}
\usepackage[utf8]{inputenc}
\usepackage{tikz}

\definecolor{darkblue}{rgb}{0,0.08,0.45}
\hypersetup{colorlinks,linkcolor=darkblue,citecolor=darkblue,filecolor=darkblue,urlcolor=darkblue}
\hypersetup{bookmarksopen,bookmarksnumbered}

\newcommand{\dy}{\ensuremath{\,\mathrm{d}\mathbf{y}}}
\newcommand{\normal}{\ensuremath{\hat{\mathbf{n}}}}
\newcommand{\ptot}{\ensuremath{p_{\mathrm{tot}}}}
\newcommand{\psca}{\ensuremath{p_{\mathrm{sca}}}}
\newcommand{\pinc}{\ensuremath{p_{\mathrm{inc}}}}
\newcommand{\pint}{\ensuremath{p_{\mathrm{int}}}}
\newcommand{\SLP}{\ensuremath{\mathcal{V}}}
\newcommand{\DLP}{\ensuremath{\mathcal{K}}}

\newcommand{\SL}{\ensuremath{V}}
\newcommand{\DL}{\ensuremath{K}}
\newcommand{\AD}{\ensuremath{T}}
\newcommand{\HS}{\ensuremath{D}}
\newcommand{\ID}{\ensuremath{I}}

\newcommand{\NtDi}{\ensuremath{\Lambda_\mathrm{NtD}^-}}
\newcommand{\DtNi}{\ensuremath{\Lambda_\mathrm{DtN}^-}}

\newcommand{\Hhalf}{\ensuremath{\mathcal{H}^{1/2}}}
\newcommand{\Hminushalf}{\ensuremath{\mathcal{H}^{-1/2}}}
\newcommand{\kint}{\ensuremath{k_-}}
\newcommand{\kext}{\ensuremath{k_+}}
\newcommand{\rhoint}{\ensuremath{\rho_-}}
\newcommand{\rhoext}{\ensuremath{\rho_+}}

\newcommand{\traceDe}{\ensuremath{\gamma_D^+}}
\newcommand{\traceDi}{\ensuremath{\gamma_D^-}}
\newcommand{\traceDei}{\ensuremath{\gamma_D^\pm}}
\newcommand{\traceNe}{\ensuremath{\gamma_N^+}}
\newcommand{\traceNi}{\ensuremath{\gamma_N^-}}
\newcommand{\traceNei}{\ensuremath{\gamma_N^\pm}}

\newcommand{\traceDme}{\ensuremath{\gamma_{D,m}^+}}

\newcommand{\traceNme}{\ensuremath{\gamma_{N,m}^+}}

\title{Boundary integral formulations for acoustic modelling of high-contrast media\footnote{© 2021. This manuscript version is made available under the CC-BY-NC-ND 4.0 license. This manuscript is published in Computers and Mathematics with Applications in final form at \url{https://doi.org/10.1016/j.camwa.2021.11.021}.}}
\author{Elwin van 't Wout\thanks{Institute for Mathematical and Computational Engineering, School of Engineering and Faculty of Mathematics, Pontificia Universidad Católica de Chile, Santiago, Chile. Contact: e.wout@uc.cl}
\and 
Seyyed R.~Haqshenas\thanks{Department of Mechanical Engineering, University College London, London, United Kingdom.}
\and
Pierre Gélat\footnotemark[2]
\and
Timo Betcke\thanks{Department of Mathematics, University College London, London, United Kingdom.}
\and
Nader Saffari\footnotemark[2]}
\date{December 6, 2021}

\begin{document}

\maketitle
	
\begin{abstract}
	The boundary element method is an efficient algorithm for simulating acoustic propagation through homogeneous objects embedded in free space. The conditioning of the system matrix strongly depends on physical parameters such as density, wavespeed and frequency. In particular, high contrast in density and wavespeed across a material interface leads to an ill-conditioned discretisation matrix. Therefore, the convergence of Krylov methods to solve the linear system is slow. Here, specialised boundary integral formulations are designed for the case of acoustic scattering at high-contrast media. The eigenvalues of the resulting system matrix accumulate at two points in the complex plane that depend on the density ratio and stay away from zero. The spectral analysis of the Calderón preconditioned PMCHWT formulation yields a single accumulation point. Benchmark simulations demonstrate the computational efficiency of the high-contrast Neumann formulation for scattering at high-contrast media.
\end{abstract}

\section{Introduction}

The boundary element method (BEM) numerically solves the Helmholtz equation by discretising a boundary integral equation at material interfaces~\cite{nedelec2001acoustic, steinbach2008numerical, sauter2010boundary}. The reformulation of a volumetric scattering model into a surface potential problem gives the BEM several computational advantages over numerical methods such as the finite element method that directly discretise the Helmholtz equation. The BEM does not require artificial boundary conditions for exterior scattering. Efficiency is obtained through preconditioned linear solvers and the fast multipole method or hierarchical matrix compression for matrix arithmetic~\cite{chew2001fast, betcke2017computationally}. Modern BEM implementations can accurately simulate large-scale wave scattering phenomena in acoustics, electromagnetics and elastodynamics. On the downside, the computational efficiency of BEM can deteriorate significantly for specific material configurations. This article studies geometries that involve material interfaces with high contrast in mass density and acoustic wavespeed for which standard boundary integral formulations become ill-conditioned.

Among the many different boundary integral formulations for acoustic transmission problems~\cite{wout2021benchmarking}, the first-kind PMCHWT (Poggio-Miller-Chang-Harrington-Wu-Tsai)~\cite{poggio1973integral, chang1974surface, wu1977scattering-bor} and the second-kind Müller boundary integral equations~\cite{muller1957grundprobleme} are among the most widely used. A well-established functional analysis of these boundary integral formulations is available~\cite{kress1978transmission, costabel1985direct}. This study considers the specific case of high-contrast media and analyses the influence of the mass density on the conditioning of the linear system obtained by Galerkin discretisation. Specifically, large jumps in density and wavespeed at material interfaces lead to slow convergence of iterative linear solvers such as GMRES. This deterioration in computational efficiency has been observed in numerical simulations reported in acoustic BEM literature (e.g.~\cite{niino2012preconditioning, wout2021benchmarking}) but no spectral analysis or specialised boundary integral formulations have been presented so far. The present work analyses the eigenvalue accumulation points for the standard boundary integral formulations and designs novel boundary integral formulations that remain well-conditioned for high-contrast media.

The simulation of acoustic scattering at high-contrast materials is of great interest to a variety of engineering applications. In biomedical ultrasound modelling, the presence of bone requires computationally efficient algorithms~\cite{haqshenas2021fast}. In underwater acoustics, resonances occur at water-air interfaces~\cite{wout2021proximity}. In material sciences, acoustic metamaterials can have arbitrary effective density characteristics~\cite{cummer2016controlling}. In the field of computational electromagnetics, high dielectric or magnetic contrasts are common in metamaterials, and specialised boundary integral formulations for Maxwell's equations have recently been designed that are efficient for high-contrast media~\cite{gossye2018calderon, gossye2019electromagnetic}.

The present study follows the same approach as in~\cite{gossye2018calderon, gossye2019electromagnetic} for electromagnetics by designing novel boundary integral formulations for acoustics based on a mix of direct interior and indirect exterior representation formulas. The only study that has considered these acoustic formulations (also called mixed-potential formulations) so far is~\cite{wout2021benchmarking}, where we provided benchmarks of a wealth of formulations, but only at low-contrast media and without a detailed analysis of their efficiency. Here, we present the design of the novel high-contrast formulations, analyse their spectrum and perform computational simulations. The high-contrast formulations are of the second kind with two accumulation points that depend on the density contrast. Furthermore, being indirect formulations, they require only half of the boundary integral operators present in standard formulations, such as the PMCHWT and Müller formulations. Hence, the novel formulations are quicker to assemble, have faster matrix-vector multiplications, and require less memory. Finally, Calderón preconditioning will be applied to the PMCHWT formulation, which has been shown to improve the conditioning at low-contrast media for acoustics~\cite{niino2012preconditioning, antoine2008integral} and electromagnetics~\cite{yan2010comparative, cools2011calderon, niino2012calderon}. The Calderón preconditioned PMCHWT has a single accumulation point and remains well conditioned for a wide range of density ratios and frequencies but has a high computational footprint at large-scale simulations. 

The boundary integral formulations for acoustic scattering at high-contrast media will be designed in Section~\ref{sec:formulation}. It will be shown in Section~\ref{sec:analysis} that the Calderón preconditioned PMCHWT and the high-contrast formulation are linear systems with accumulation points of eigenvalues, which yields well-conditioned formulations at high-contrast media. Section~\ref{sec:results} provides extensive numerical simulations of the BEM, corroborating the conditioning of the formulations at high-contrast media. The computational benchmarks also include numerical simulations at multiple objects as well as large-scale geometries.

\section{Formulation}
\label{sec:formulation}

The boundary integral formulations for high-contrast media will be derived for a single penetrable object. The extension to the more general case of multiple scattering will be explained at the end of this section.

\subsection{Model equations}

Let us consider a three-dimensional bounded object denoted by $\Omega^-$ whose surface $\Gamma$ is piecewise smooth. The unbounded exterior region is denoted by $\Omega^+$. The acoustic pressure field in the exterior is denoted by $\ptot$ and can be decomposed into an unknown scattered field $\psca$ and a known incident field $\pinc$ as $\ptot = \psca + \pinc$. In the interior, the acoustic pressure field is denoted by $\pint$. Harmonic wave propagation through material with a linear response can then by modelled by the Helmholtz equation as
\begin{equation}
	\label{eq:helmholtz}
	\begin{cases}
		-\Delta \psca - k_+^2 \psca = 0 & \text{in } \Omega^+, \\
		-\Delta \pint - k_-^2 \pint = 0 & \text{in } \Omega^-, \\
		\traceDe \ptot = \traceDi \pint & \text{on } \Gamma, \\
		\frac1{\rho_+} \traceNe \ptot = \frac1{\rho_-} \traceNi \pint & \text{on } \Gamma, \\
		\lim_{\mathbf{r} \to \infty} |\mathbf{r}| \left(\partial_{|\mathbf{r}|} \psca - \imath k_+ \psca\right) = 0.
	\end{cases}
\end{equation}
Here, $k_\pm$ and $\rho_\pm$ denote the wavenumber and mass density in region~$\Omega^\pm$, respectively. The Dirichlet and Neumann traces are defined by
\begin{subequations}
	\begin{align}
		\traceDei p(\mathbf{x}) &= \lim_{\mathbf{y} \to \mathbf{x}} p(\mathbf{y}) && \text{for } \mathbf{x} \in \Gamma \text{ and } \mathbf{y} \in \Omega^\pm, \\
		\traceNei p(\mathbf{x}) &= \lim_{\mathbf{y} \to \mathbf{x}} \nabla p(\mathbf{y}) \cdot \normal && \text{for } \mathbf{x} \in \Gamma \text{ and } \mathbf{y} \in \Omega^\pm,
	\end{align}
\end{subequations}
where the unit normal vector~$\normal$ points towards the exterior domain. The last equation is the Sommerfeld radiation condition for outgoing wave fields where $\imath$ denotes the imaginary unit. Furthermore, it is assumed that the incident wave field satisfies the Helmholtz equation with wavenumber $k_+$.

\subsection{High-contrast boundary integral formulations}

Any field that satisfies the Helmholtz system can be represented by surface potentials at the material interface~\cite{nedelec2001acoustic, steinbach2008numerical, sauter2010boundary}. Let us use a direct representation formula for the interior field, that is,
\begin{equation}
	\label{eq:representation:interior}
	\pint = \SLP^- \psi^- - \DLP^- \phi^-.
\end{equation}
Differently, the exterior field is represented by one of the following indirect representation formulas:
\begin{subequations}
	\label{eq:representation:exterior}
	\begin{align}
		\label{eq:representation:exterior:single}
		\psca &= \SLP^+ \psi^+, \\
		\label{eq:representation:exterior:double}
		\psca &= -\DLP^+ \phi^+.
	\end{align}
\end{subequations}
Here, $\SLP$ and $\DLP$ denote the single-layer and double-layer potentials operators, which are given by
\begin{subequations}
	\label{eq:operators:potential}
	\begin{align}
		\label{eq:slp}
		\left[\SLP^\pm\psi\right]\!(\mathbf{x}) &= \iint_\Gamma G^\pm(\mathbf{x},\mathbf{y}) \psi(\mathbf{y}) \dy && \text{for } \mathbf{x} \in \Omega^\pm, \\
		\label{eq:dlp}
		\left[\DLP^\pm\phi\right]\!(\mathbf{x}) &= \iint_\Gamma \frac{\partial}{\partial \normal(\mathbf{y})} G^\pm(\mathbf{x},\mathbf{y}) \phi(\mathbf{y}) \dy && \text{for } \mathbf{x} \in \Omega^\pm,
	\end{align}
\end{subequations}
for the Green's function
\begin{equation}
	\label{eq:green}
	G^\pm(\mathbf{x},\mathbf{y}) = \frac{e^{\imath k_\pm |\mathbf{x} - \mathbf{y}|}}{4\pi |\mathbf{x} - \mathbf{y}|} \quad \text{for } \mathbf{x},\mathbf{y} \in \Omega^\pm \text{ and } \mathbf{x} \ne \mathbf{y}.
\end{equation}
 Since the interior field is represented by the direct formula~\eqref{eq:representation:interior}, the surface potentials are the traces of the pressure field~\cite{steinbach2008numerical}, that is,
\begin{subequations}
	\begin{align}
		\phi^- &= \traceDi \pint, \\
		\psi^- &= \traceNi \pint.
	\end{align}
\end{subequations}
On the contrary, the exterior surface potentials do not have a direct physical meaning since the exterior field is given by the indirect representation formula~\eqref{eq:representation:exterior}. Taking the traces of the scattered field towards the surface yield
\begin{subequations}
	\label{eq:bie:exterior:slp}
	\begin{align}
		\traceDe \psca &= \SL^+ \psi^+, \\
		\traceNe \psca &= \left(\AD^+ - \frac12\ID\right) \psi^+
	\end{align}
\end{subequations}
for the indirect single-layer representation~\eqref{eq:representation:exterior:single}, and
\begin{subequations}
	\label{eq:bie:exterior:dlp}
	\begin{align}
		\traceDe \psca &= -\left(\DL^+ + \frac12\ID\right) \phi^+, \\
		\traceNe \psca &= \HS^+ \phi^+
	\end{align}
\end{subequations}
for the indirect double-layer representation~\eqref{eq:representation:exterior:double}.
Here, $\ID$ denotes the identity operator and
\begin{subequations}
	\label{eq:operators:boundary}
	\begin{align}
		\label{eq:sl}
		\left[\SL^\pm\psi\right](\mathbf{x}) &= \iint_\Gamma G^\pm(\mathbf{x},\mathbf{y}) \psi(\mathbf{y}) \dy && \text{for } \mathbf{x} \in \Gamma, \\
		\label{eq:dl}
		\left[\DL^\pm\phi\right](\mathbf{x}) &= \iint_\Gamma \frac{\partial}{\partial \normal(\mathbf{y})} G^\pm(\mathbf{x},\mathbf{y}) \phi(\mathbf{y}) \dy && \text{for } \mathbf{x} \in \Gamma, \\
		\label{eq:ad}
		\left[\AD^\pm\phi\right](\mathbf{x}) &= \frac{\partial}{\partial \normal(\mathbf{x})} \iint_\Gamma G^\pm(\mathbf{x},\mathbf{y}) \phi(\mathbf{y}) \dy && \text{for } \mathbf{x} \in \Gamma, \\
		\label{eq:hs}
		\left[\HS^\pm\phi\right](\mathbf{x}) &= -\frac{\partial}{\partial \normal(\mathbf{x})} \iint_\Gamma \frac{\partial}{\partial \normal(\mathbf{y})} G^\pm(\mathbf{x},\mathbf{y}) \phi(\mathbf{y}) \dy && \text{for } \mathbf{x} \in \Gamma,
	\end{align}
\end{subequations}
are the single-layer, double-layer, adjoint double-layer, and hypersingular boundary integral operators, respectively.

Notice that the set of two boundary integral equations~\eqref{eq:bie:exterior:slp} has three unknown surface potentials: $\traceDe\ptot$, $\traceNe\ptot$, and $\psi^+$. Similarly, the set of two boundary integral equations~\eqref{eq:bie:exterior:dlp} has three unknown surface potentials: $\traceDe\ptot$, $\traceNe\ptot$, and $\phi^+$. At the same time, the interior representation formula~\eqref{eq:representation:interior} is written in terms of the two unknown surface potentials $\traceDi\ptot$ and $\traceNi\ptot$. In the design of most boundary integral equations (cf.~\cite{wout2021benchmarking}), the two interface conditions couple the pairs of two traces at each interface, resulting in a well-defined set of boundary integral equations. This case is different, with an additional surface potential $\psi^+$ or $\phi^+$ present in the formulation. Hence, the standard design procedures will fail.

The design of the high-contrast formulations follows a different approach than usual and includes the Neumann-to-Dirichlet (NtD) and Dirichlet-to-Neumann (DtN) maps. The interior NtD and DtN maps are implicitly defined as
\begin{subequations}
	\label{eq:operators:ntdanddtn}
	\begin{align}
		\label{eq:ntd}
		\NtDi\left(\traceNi\pint\right) &= \traceDi\pint, \\
		\label{eq:dtn}
		\DtNi\left(\traceDi\pint\right) &= \traceNi\pint
	\end{align}
\end{subequations}
and are also known as the Poincaré-Steklov and Steklov-Poincaré operators. There are two equivalent expressions for these operators~\cite{sauter2010boundary}, that is,
\begin{subequations}
	\label{eq:ntd:expressions}
	\begin{align}
		\label{eq:ntd:sl}
		\NtDi &= \left(\frac12\ID + \DL^-\right)^{-1} \SL^-, \\
		\label{eq:ntd:hs}
		\NtDi &= \left(\HS^-\right)^{-1} \left(\frac12\ID - \AD^-\right)
	\end{align}
\end{subequations}
and
\begin{subequations}
	\label{eq:dtn:expressions}
	\begin{align}
		\label{eq:dtn:hs}
		\DtNi &= \left(\frac12\ID - \AD^-\right)^{-1} \HS^-, \\
		\label{eq:dtn:sl}
		\DtNi &= \left(\SL^-\right)^{-1} \left(\frac12\ID + \DL^-\right).
	\end{align}
\end{subequations}
Notice that no exterior NtD or DtN maps can be used since the exterior field is represented by an indirect formula. Furthermore, no closed-form expressions of these operators are available due to the presence of the inverse operators.

Now, the idea is to use the NtD and DtN maps, as well as the transmission conditions~\eqref{eq:helmholtz}, to design relations between the exterior traces. Specifically, one can write
\begin{subequations}
	\label{eq:transmission}
	\begin{align}
		\label{eq:transmission:ntd}
		&\traceDe \ptot = \traceDi \pint = \NtDi \traceNi \pint = \frac{\rho_-}{\rho_+} \NtDi \traceNe \ptot, \\
		\label{eq:transmission:dtn}
		&\traceNe \ptot = \frac{\rho_+}{\rho_-} \traceNi \pint = \frac{\rho_+}{\rho_-} \DtNi \traceDi \pint = \frac{\rho_+}{\rho_-} \DtNi \traceDe \ptot.
	\end{align}
\end{subequations}
These relations can be used to eliminate one of the unknown field traces in the sets of boundary integral equations~\eqref{eq:bie:exterior:slp} and~\eqref{eq:bie:exterior:dlp}. By either eliminating the Dirichlet trace or the Neumann trace from the two sets of boundary integral equations, four different boundary integral formulations appear~\cite{wout2021benchmarking}. Here, we will only consider the two formulations that result in well-conditioned systems.

In the case of the set of boundary integral equations~\eqref{eq:bie:exterior:slp} that correspond to the indirect single-layer representation~\eqref{eq:representation:exterior:single}, the surface potential $\psi^+$ is part of $\Hminushalf(\Gamma)$. Since the Neumann trace is part of $\Hminushalf(\Gamma)$ as well, the Dirichlet trace will be eliminated. That is, substituting relation~\eqref{eq:transmission:ntd} into the set of boundary integral equations~\eqref{eq:bie:exterior:slp} yields
\begin{subequations}
	\label{eq:bie:exterior:slp:ntd}
	\begin{align}
		\label{eq:bie:ntd}
		\frac{\rho_-}{\rho_+} \NtDi \traceNe \ptot - \traceDe \pinc &= \SL^+ \psi^+, \\
		\traceNe \ptot - \traceNe \pinc &= \left(\AD^+ - \frac12\ID\right) \psi^+.
	\end{align}
\end{subequations}
Similarly, the surface potential $\phi^+$ is part of $\Hhalf(\Gamma)$ which suggests eliminating the Neumann trace from the double-layer formulation. That is, substituting relation~\eqref{eq:transmission:dtn} into the set of boundary integral equations~\eqref{eq:bie:exterior:dlp} yields
\begin{subequations}
	\label{eq:bie:exterior:dlp:dtn}
	\begin{align}
		\traceDe \ptot - \traceDe \pinc &= -\left(\DL^+ + \frac12\ID\right) \phi^+, \\
		\label{eq:bie:dtn}
		\frac{\rho_+}{\rho_-} \DtNi \traceDe \ptot - \traceNe \pinc &= \HS^+ \phi^+.
	\end{align}
\end{subequations}
Both sets of boundary integral equations~\eqref{eq:bie:exterior:slp:ntd} and~\eqref{eq:bie:exterior:dlp:dtn} are well-defined with two independent equations for two unknown surface potentials. However, these formulations cannot be discretised yet due to the presence of the NtD and DtN operators.

To solve the issue of having the NtD and DtN operators, which have no closed-form expressions, their expressions in terms of (inverse) boundary integral operators will be substituted. That is, two expressions~\eqref{eq:ntd:expressions} are available for the NtD map and two expressions~\eqref{eq:dtn:expressions} for the DtN map. Hence, two different versions of the boundary integral equation~\eqref{eq:bie:ntd} can be designed, as well as two different versions of the equation~\eqref{eq:bie:dtn}. Here, the expression~\eqref{eq:ntd:hs} will be substituted into equation~\eqref{eq:bie:ntd} to obtain
\begin{equation}
	\frac{\rho_-}{\rho_+} \left(\HS^-\right)^{-1} \left(\frac12\ID - \AD^-\right) \traceNe \ptot - \traceDe \pinc = \SL^+ \psi^+.
\end{equation}
When multiplying the equation from the left by the hypersingular operator, the set of boundary integral equations~\eqref{eq:bie:exterior:slp:ntd} becomes
\begin{subequations}
	\label{eq:bie:slp:ntd}
	\begin{align}
		\left(\frac12\ID - \AD^-\right) \traceNe \ptot - \frac{\rho_+}{\rho_-} \HS^- \SL^+ \psi^+ &= \frac{\rho_+}{\rho_-} \HS^- \traceDe \pinc, \\
		\traceNe \ptot + \left(\frac12\ID - \AD^+\right) \psi^+ &= \traceNe \pinc.
	\end{align}
\end{subequations}
Notice that these boundary integral equations are both part of $\Hminushalf(\Gamma)$, which would not have been the case if the expression~\eqref{eq:ntd:sl} was chosen.
Similarly, substituting expression~\eqref{eq:dtn:sl} for the DtN map into the boundary integral equation~\eqref{eq:bie:dtn} yields
\begin{equation}
	\frac{\rho_+}{\rho_-} \left(\SL^-\right)^{-1} \left(\frac12\ID + \DL^-\right) \traceDe \ptot - \traceNe \pinc = \HS^+ \phi^+.
\end{equation}
Multiplication from the left by the single-layer operator converts the set of boundary integral equations~\eqref{eq:bie:exterior:dlp:dtn} into
\begin{subequations}
	\label{eq:bie:dlp:dtn}
	\begin{align}
		\left(\frac12\ID + \DL^+\right) \phi^+ + \traceDe \ptot &= \traceDe \pinc, \\
		-\frac{\rho_-}{\rho_+} \SL^- \HS^+ \phi^+ + \left(\frac12\ID + \DL^-\right) \traceDe \ptot &= \frac{\rho_-}{\rho_+} \SL^- \traceNe \pinc.
	\end{align}
\end{subequations}
In summary, the two sets of boundary integral equations read
\begin{subequations}
	\label{eq:bie:highcontrast}
	\begin{align}
		\label{eq:bie:highcontrast:neu}
		\begin{bmatrix} \frac12\ID - \AD^- & -\frac{\rho_+}{\rho_-} \HS^- \SL^+ \\ \ID & \frac12\ID - \AD^+ \end{bmatrix}
		\begin{bmatrix} \traceNe \ptot \\ \psi^+ \end{bmatrix}
		&= \begin{bmatrix} \frac{\rho_+}{\rho_-} \HS^- \traceDe \pinc \\ \traceNe \pinc \end{bmatrix}, \\
		\label{eq:bie:highcontrast:dir}
		\begin{bmatrix} \frac12\ID + \DL^+ & \ID \\ -\frac{\rho_-}{\rho_+} \SL^- \HS^+ & \frac12\ID + \DL^- \end{bmatrix}
		\begin{bmatrix} \phi^+ \\ \traceDe \ptot \end{bmatrix}
		&= \begin{bmatrix} \traceDe \pinc \\ \frac{\rho_-}{\rho_+} \SL^- \traceNe \pinc \end{bmatrix}
	\end{align}
\end{subequations}
which will be called the high-contrast Neumann and Dirichlet boundary integral formulation, respectively. Notice that the Neumann version maps $\Hminushalf(\Gamma) \times \Hminushalf(\Gamma)$ into $\Hminushalf(\Gamma) \times \Hminushalf(\Gamma)$ and the Dirichlet version maps $\Hhalf(\Gamma) \times \Hhalf(\Gamma)$ into $\Hhalf(\Gamma) \times \Hhalf(\Gamma)$. Out of the eight different boundary integral formulations that can be designed by this framework of mixed potentials (cf.~\cite{wout2021benchmarking}), these are the only two that are completely defined in the same function space and are, therefore, the only second-kind boundary integral formulations.

\subsubsection{Extension to multiple scattering}

The high-contrast formulations can readily be extended to the case of multiple scattering from a set of disjoint objects. That is, let us consider objects $\Omega_m$ for $m=1,2,\dots,\ell$ with interior wavenumber $k_m$ and density $\rho_m$, that are all embedded in the exterior region $\Omega_0$ with wavenumber $k_0$ and density $\rho_0$. Then, the boundary integral formulation reads
\begin{equation}
	\label{eq:system:multiplescattering}
	\begin{bmatrix}
		B_{11} & B_{12} & \cdots & B_{1\ell} \\
		B_{21} & B_{22} & \cdots & B_{2\ell} \\
		\vdots & \vdots & \ddots & \vdots \\
		B_{\ell 1} & B_{\ell 2} & \cdots & B_{\ell\ell}
	\end{bmatrix}
	\begin{bmatrix} u_1 \\ u_2 \\ \vdots \\ u_\ell \end{bmatrix}
	= \begin{bmatrix} f_1 \\ f_2 \\ \vdots \\ f_\ell \end{bmatrix}
\end{equation}
with the blocks given by
\begin{subequations}
	\begin{align}
		B_{mm} &= \begin{bmatrix} \frac12\ID - \AD_m^- & -\frac{\rho_0}{\rho_m} \HS_m^- \SL_{mm}^+ \\ \ID & \frac12\ID - \AD_{mm}^+ \end{bmatrix}, \\
		B_{mn} &= \begin{bmatrix} 0 & -\frac{\rho_0}{\rho_m} \HS_m^- \SL_{mn}^+ \\ 0 & -\AD_{mn}^+ \end{bmatrix} \quad \text{if } m \ne n, \\
		u_m &= \begin{bmatrix} \traceNme \ptot \\ \psi_m^+ \end{bmatrix}, \\
		f_m &= \begin{bmatrix} \frac{\rho_0}{\rho_m} \HS_m^- \traceDme \pinc \\ \traceNme \pinc \end{bmatrix}
	\end{align}
\end{subequations}
for the Neumann high-contrast formulation, and
\begin{subequations}
	\begin{align}
		B_{mm} &= \begin{bmatrix} \frac12\ID + \DL_{mm}^+ & \ID \\ -\frac{\rho_m}{\rho_0} \SL_m^- \HS_{mm}^+ & \frac12\ID + \DL_m^- \end{bmatrix}, \\
		B_{mn} &= \begin{bmatrix} \DL_{mn}^+ & 0 \\ -\frac{\rho_m}{\rho_0} \SL_m^- \HS_{mn}^+ & 0 \end{bmatrix} \quad \text{if } m \ne n, \\
		u_m &= \begin{bmatrix} \phi_m^+ \\ \traceDme \ptot \end{bmatrix}, \\
		f_m &= \begin{bmatrix} \traceDme \pinc \\ \frac{\rho_m}{\rho_0} \SL_m^- \traceNme \pinc \end{bmatrix}
	\end{align}
\end{subequations}
for the Dirichlet high-contrast formulation.

\subsection{Standard formulations}

The most widely used boundary integral formulations for Helmholtz transmission problems are the PMCHWT and Müller formulations, which are given by
\begin{equation}
	\label{eq:pmchwt}
	\begin{bmatrix} -\DL^+ - \DL^- & \SL^+ + \frac{\rho_-}{\rho_+} \SL^- \\ \HS^+ + \frac{\rho_+}{\rho_-}\HS^- & \AD^+ + \AD^- \end{bmatrix}
	\begin{bmatrix} \traceDe \ptot \\ \traceNe \ptot \end{bmatrix}
	=
	\begin{bmatrix} \traceDe \pinc \\ \traceNe \pinc \end{bmatrix}
\end{equation}
and
\begin{equation}
	\label{eq:muller}
	\begin{bmatrix} \ID - \DL^+ + \DL^- & \SL^+ - \frac{\rho_-}{\rho_+} \SL^- \\ \HS^+ - \frac{\rho_+}{\rho_-}\HS^- & \ID + \AD^+ - \AD^- \end{bmatrix}
	\begin{bmatrix} \traceDe \ptot \\ \traceNe \ptot \end{bmatrix}
	=
	\begin{bmatrix} \traceDe \pinc \\ \traceNe \pinc \end{bmatrix},
\end{equation}
respectively. The PMCHWT formulation is a first-kind boundary integral equation and the Müller formulation is of second kind.

Since the PMCHWT formulation is of first kind, it is often preconditioned to improve the convergence of the linear solver. One of the most effective techniques is to use the Calderón identities, in specific the projection property of the Calderón operators. The Calderón preconditioned PMCHWT formulation reads
\begin{align}
	\label{eq:pmchwt:calderon}
	&\begin{bmatrix} -\DL^+ - \DL^- & \SL^+ + \frac{\rho_-}{\rho_+} \SL^- \\ \HS^+ + \frac{\rho_+}{\rho_-}\HS^- & \AD^+ + \AD^- \end{bmatrix}^2
	\begin{bmatrix} \traceDe \ptot \\ \traceNe \ptot \end{bmatrix} \nonumber \\
	&\quad =
	\begin{bmatrix} -\DL^+ - \DL^- & \SL^+ + \frac{\rho_-}{\rho_+} \SL^- \\ \HS^+ + \frac{\rho_+}{\rho_-}\HS^- & \AD^+ + \AD^- \end{bmatrix}
	\begin{bmatrix} \traceDe \pinc \\ \traceNe \pinc \end{bmatrix}
\end{align}
which is a well-conditioned formulation~\cite{antoine2008integral, niino2012preconditioning}.
All of these formulations can readily be extended to multiple scattering~\cite{wout2021benchmarking}.

Notice that the standard formulations, which are based on direct representation formulas, include both density contrasts $\rho_-/\rho_+$ and $\rho_+/\rho_-$ in the system matrix. Hence, an imbalance in matrix elements may occur at high-contrast media. Differently, the high-contrast Dirichlet and Neumann formulations include a single density contrast only: either $\rho_-/\rho_+$ or $\rho_+/\rho_-$, respectively.

\subsection{Numerical discretisation}

The numerical discretisation of the boundary integral formulations follows a Galerkin method with piecewise linear (P1) elements, associated to each node in a triangular surface mesh. Notice that P0 elements could be used for the $\Hminushalf(\Gamma)$ spaces as well. However, since operator products are present, these P0 elements would have to be defined on the dual mesh~\cite{betcke2020product}, which increases the computational overhead.

\section{Spectral analysis}
\label{sec:analysis}

The linear system resulting from the BEM is a dense matrix. Solving this system for large-scale simulation is computationally expensive and iterative linear solvers are preferred over direct factorisations~\cite{marburg2003performance, sakuma2008fast}. Since the complex-valued matrix is not Hermitian, the GMRES solver~\cite{saad1986gmres} will be used. The convergence of Krylov solvers depends on the spectrum of the matrix, where a low condition number and clustering of eigenvalues typically lead to small numbers of iterations. The second-kind boundary integral equations have good convergence properties since these formulations are in the form of an identity operator plus a compact operator, which yields a spectrum with eigenvalues accumulating near a fixed point. For this reason, first-kind boundary integral formulations are often preconditioned such that the preconditioned system is of second kind. First, we will analyse the spectral properties of operator products. Then, this analysis will be applied to the high-contrast formulations and the Calderón preconditioned PMCHWT formulation. Both have accumulation points that depend on the density ratio across the material interface.

The following spectral analysis relies on the compactness of the (adjoint) double-layer operator. Hence, the material interface~$\Gamma$ is assumed to be $\mathcal{C}^2$ smooth, for which the boundary integral operators are compact~\cite{costabel1985direct, colton2013integral, sauter2010boundary}. Notice that the boundary integral formulations can be applied to Lipschitz surfaces as well (cf.~\cite{costabel1988boundary}) but at the expense of corner singularities and lack of compactness properties.

\subsection{Operator products}
\label{sec:operatorproducts}

\begin{figure}[!ht]
	\centering
	\begin{subfigure}[b]{0.49\columnwidth}
		\includegraphics[width=\textwidth]{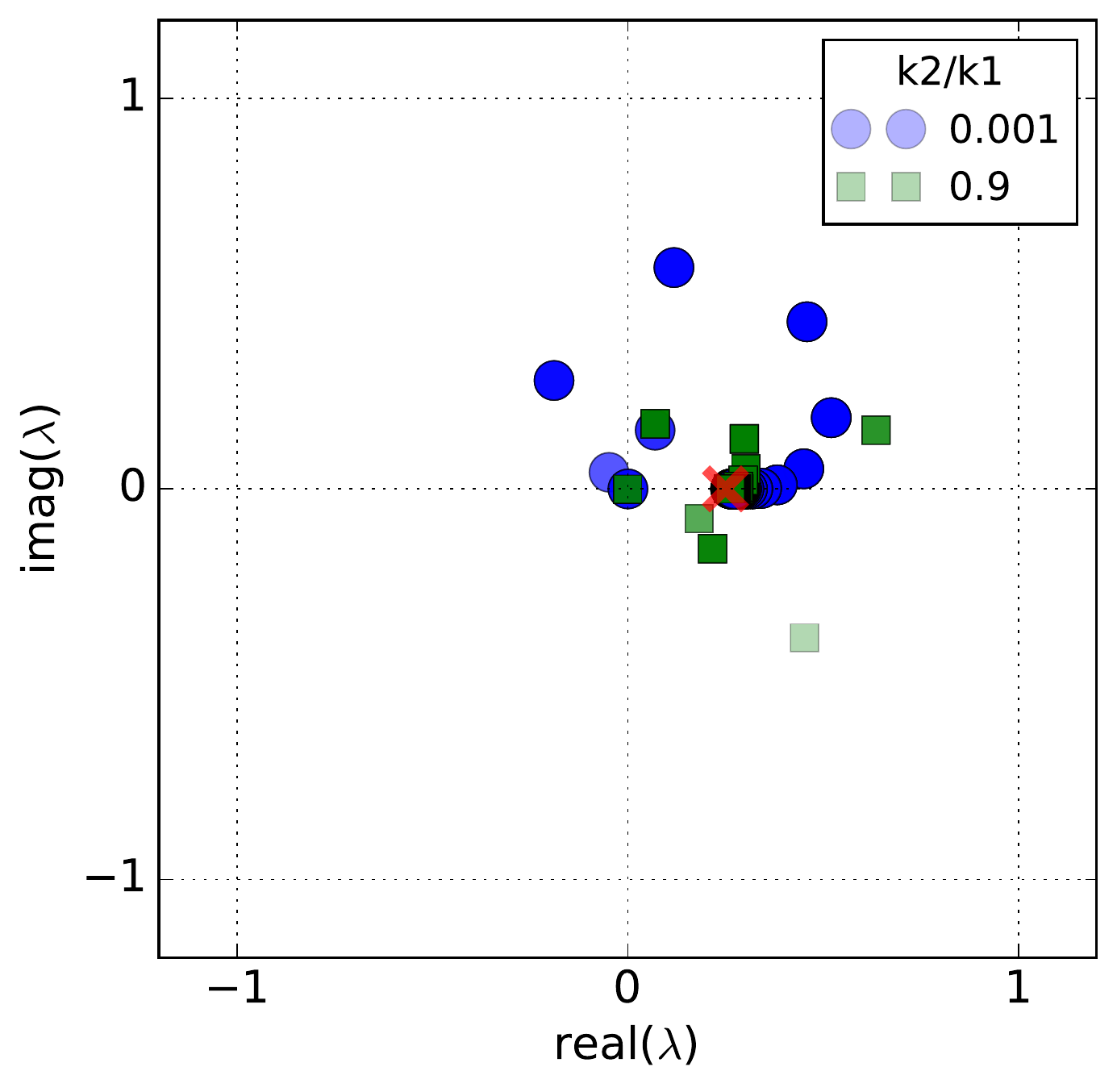}
		\caption{Operator $\SL_1 \HS_2$ or, equivalently, $\HS_2 \SL_1$.}
	\end{subfigure}
	\hfill
	\begin{subfigure}[b]{0.49\columnwidth}
		\includegraphics[width=\textwidth]{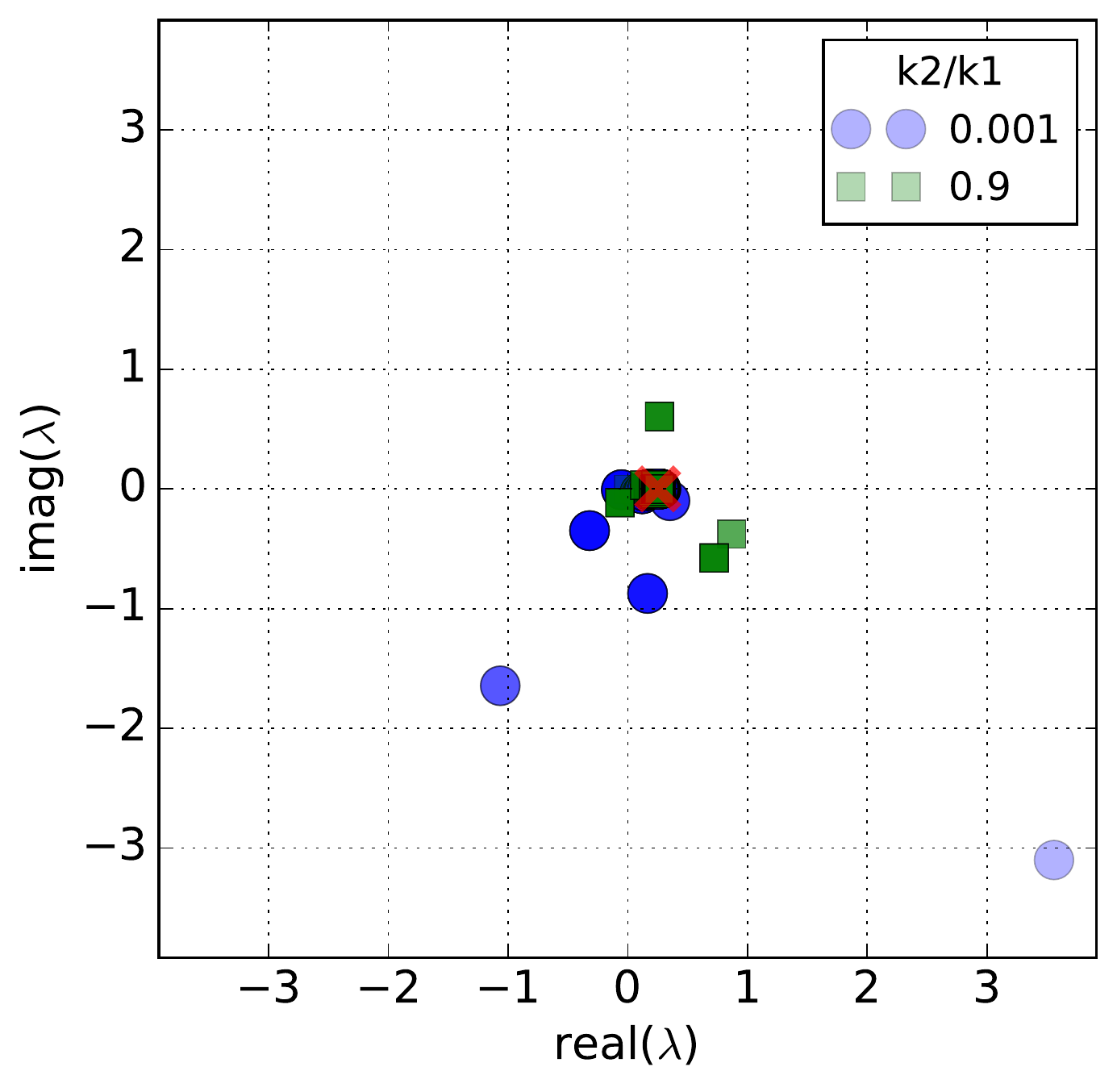}
		\caption{Operator $\SL_2 \HS_1$ or, equivalently, $\HS_1 \SL_2$.}
	\end{subfigure}
	\caption{The eigenvalues of operator products for $k_1 = 7$ and $k_2 \in \{0.007, 6.93\}$, calculated at a sphere with unit radius and 1262 degrees of freedom. The red cross visualises the expected accumulation point at $1/4$.}
	\label{fig:spectrum:products}
\end{figure}

Let us consider the operator products of single-layer and hypersingular operators that are present in the formulations. The Calderón identities~\cite{steinbach2008numerical} state that
\begin{subequations}
	\label{eq:calderonidentities}
	\begin{align}
		\label{eq:slhs}
		\SL \HS &= \frac14\ID - \DL^2, \\
		\label{eq:hssl}
		\HS \SL &= \frac14\ID - \AD^2.
	\end{align}
\end{subequations}
Since the double-layer and adjoint double-layer operators are compact, there is an accumulation point at $1/4$. In other words, the eigenvalues of these operator products are clustered around the point $1/4$ in the complex plane. Since the Calderón identities hold for operators with the same wavenumber only, these results cannot be used directly for operator products with different wavenumbers. Still, single-layer and hypersingular operators with wavenumber $k_1$ are compact perturbations of single-layer and hypersingular operators with wavenumber $k_2$, respectively~\cite{kleinman1988single, antoine2008integral, claeys2013multi, boubendir2015integral}. For this reason, the spectrum of operator products $\SL_1 \HS_2$ and $\HS_1 \SL_2$ also accumulate at $1/4$. This is computationally validated in Figure~\ref{fig:spectrum:products}, where one can also observe a larger spread in eigenvalues for higher contrasts in wavenumber.

\subsection{High-contrast formulations}

The matrix of the Neumann high-contrast formulation~\eqref{eq:bie:highcontrast:neu} can be written as
\begin{equation}
	\begin{bmatrix} \frac12\ID - \AD^- & -\frac{\rho_+}{\rho_-} \HS^- \SL^+ \\ \ID & \frac12\ID - \AD^+ \end{bmatrix}
	= \begin{bmatrix} \frac12\ID & -\frac{\rho_+}{\rho_-} \HS^- \SL^+ \\ \ID & \frac12\ID \end{bmatrix}
	- \begin{bmatrix} \AD^- & 0 \\ 0 & \AD^+ \end{bmatrix}.
\end{equation}
Since the adjoint-double layer is a compact operator, the accumulation points will be determined by the first matrix on the right-hand side. By definition, any eigenvalue pair $(\lambda,\mathbf{v})$ of this matrix satisfies
\begin{align}
	\begin{bmatrix} \frac12\ID & -\frac{\rho_+}{\rho_-} \HS^- \SL^+ \\ \ID & \frac12\ID \end{bmatrix}
	\begin{bmatrix} \mathbf{v}_1 \\ \mathbf{v}_2 \end{bmatrix}
	&= \lambda \begin{bmatrix} \mathbf{v}_1 \\ \mathbf{v}_2 \end{bmatrix}.
\end{align}
The second row states that $\mathbf{v}_1 = (\lambda - \frac12) \mathbf{v}_2$. Substitution of which into the first row yields
\begin{equation}
	-\frac{\rho_+}{\rho_-} \HS^- \SL^+ \mathbf{v}_2 = \left(\lambda - \frac12\right)^2 \mathbf{v}_2.
\end{equation}
Since the operator product $\HS^- \SL^+$ has an accumulation point at $1/4$, the accumulation point~$\tilde\lambda_N$ of the Neumann high-contrast formulation satisfies
\begin{equation}
	-\frac{\rho_+}{\rho_-} \frac14 = \left(\tilde\lambda_N - \frac12\right)^2.
\end{equation}
This quadratic equation is solved for
\begin{equation}
	\label{eq:accumulation:neumann}
	\tilde\lambda_N = \frac12 \pm \frac\imath2 \sqrt{\frac{\rho_+}{\rho_-}}.
\end{equation}
With an equivalent analysis, one can show that the accumulation point of the Dirichlet high-contrast formulation~\eqref{eq:bie:highcontrast:dir} is given by
\begin{equation}
	\label{eq:accumulation:dirichlet}
	\tilde\lambda_D = \frac12 \pm \frac\imath2 \sqrt{\frac{\rho_-}{\rho_+}}.
\end{equation}
This spectral analysis concludes that the high-contrast formulations have two accumulation points that are independent of the wavenumber and stay away from zero. Furthermore, the accumulation points of the high-contrast Neumann and Dirichlet formulations only grow towards infinity when the interior or exterior density tends to zero, respectively, and only with a speed proportional to the square root of the density ratio. Hence, these formulations are expected to be well-conditioned for high-contrast media. An exceptional case where the accumulation point could become zero is when acoustic metamaterials with a negative effective mass density are considered~\cite{lu2009phononic}. The accumulation points are computationally validated in Section~\ref{sec:results:spectrum}.

\subsection{Calderón preconditioned PMCHWT}

The system matrix of the Calderón preconditioned PMCHWT formulation~\eqref{eq:pmchwt:calderon} can be written as

\begin{equation}
	\begin{bmatrix} -\DL^+ - \DL^- & \SL^+ + \frac{\rho_-}{\rho_+} \SL^- \\ \HS^+ + \frac{\rho_+}{\rho_-}\HS^- & \AD^+ + \AD^- \end{bmatrix}^2
	=
	\begin{bmatrix} C_{11} & C_{12} \\ C_{21} & C_{22} \end{bmatrix}
\end{equation}
with
\begin{subequations}
	\begin{align*}
		C_{11} &= \DL^+ \DL^+ + \DL^+ \DL^- + \DL^- \DL^+ + \DL^- \DL^- \\ &\quad + \SL^+ \HS^+ + \frac{\rho_+}{\rho_-} \SL^+ \HS^- + \frac{\rho_-}{\rho_+} \SL^- \HS^+ + \SL^- \HS^-, \\
		C_{12} &= -\DL^+ \SL^+ - \frac{\rho_-}{\rho_+} \DL^+ \SL^- - \DL^- \SL^+ - \frac{\rho_-}{\rho_+} \DL^- \SL^- \\ &\quad + \SL^+ \AD^+ + \SL^+ \AD^- + \frac{\rho_-}{\rho_+} \SL^- \AD^+ + \frac{\rho_-}{\rho_+} \SL^- \AD^-, \\
		C_{21} &= -\HS^+ \DL^+ - \HS^+ \DL^- - \frac{\rho_+}{\rho_-} \HS^- \DL^+ - \frac{\rho_+}{\rho_-} \HS^- \DL^- \\ &\quad + \AD^+ \HS^+ + \frac{\rho_+}{\rho_-} \AD^+ \HS^- + \AD^- \HS^+ + \frac{\rho_+}{\rho_-} \AD^- \HS^-, \\
		C_{22} &= \HS^+ \SL^+ + \frac{\rho_-}{\rho_+} \HS^+ \SL^- + \frac{\rho_+}{\rho_-} \HS^- \SL^+ + \HS^- \SL^- \\ &\quad + \AD^+ \AD^+ + \AD^+ \AD^- + \AD^- \AD^+ + \AD^- \AD^-.
	\end{align*}
\end{subequations}
By substituting the Calderón identities~\eqref{eq:calderonidentities}, these expressions simplify to
\begin{subequations}
	\label{eq:cp:terms}
	\begin{align}
		C_{11} &= \frac12\ID + \DL^+ \DL^- + \DL^- \DL^+ + \frac{\rho_+}{\rho_-} \SL^+ \HS^- + \frac{\rho_-}{\rho_+} \SL^- \HS^+, \\
		C_{12} &= \frac{\rho_-}{\rho_+} \left(\SL^- \AD^+ - \DL^+ \SL^-\right) + \SL^+ \AD^- - \DL^- \SL^+, \\
		C_{21} &= \AD^- \HS^+ - \HS^+ \DL^- + \frac{\rho_+}{\rho_-} \left(\AD^+ \HS^- - \HS^- \DL^+\right), \\
		C_{22} &= \frac12\ID + \AD^+ \AD^- + \AD^- \AD^+ + \frac{\rho_-}{\rho_+} \HS^+ \SL^- + \frac{\rho_+}{\rho_-} \HS^- \SL^+.
	\end{align}
\end{subequations}
This leaves us with expressions of operator products with different wavenumbers only.

Let us first consider the off-diagonal blocks, for which we need the following Calderón identities~\cite{steinbach2008numerical}:
\begin{subequations}
	\label{eq:calderonidentities:2}
	\begin{align}
		\SL \AD - \DL \SL &= 0, \\
		\AD \HS - \HS \DL &= 0.
	\end{align}
\end{subequations}
As discussed in Section~\ref{sec:operatorproducts}, a change in wavenumber is a compact perturbation of the boundary integral operator and these Calderón identities also hold for operator products with different wavenumbers, except for a compact perturbation. Hence, operators $C_{12}$ and $C_{21}$ are compact and do not influence the accumulation points of the spectrum~\cite{antoine2008integral}.

Concerning the diagonal blocks, products of the double-layer and adjoint double-layer operators are compact and products of single-layer and hypersingular operators accumulate at $1/4$, as was shown above. Hence, the Calderón preconditioned PMCHWT formulation has a single accumulation point given by
\begin{equation}
	\label{eq:accumulation:calderon}
	\tilde\lambda_C = \frac12 + \frac{\rho_+}{4\rho_-} + \frac{\rho_-}{4\rho_+}.
\end{equation}
This spectral analysis concludes that the formulation is well-conditioned since the accumulation point stays away from zero for any frequency and density ratio. Furthermore, the accumulation point only tends to infinity when either of the density ratios tends to infinity. These accumulation points are computationally validated in Section~\ref{sec:results:spectrum}.

The same analysis also shows that when considering only the exterior or interior Calderón operator as preconditioner~\cite{wout2021benchmarking}, two accumulation points are present: $\frac14 + \frac{\rho_+}{4\rho_-}$ and $\frac14 + \frac{\rho_-}{4\rho_+}$. Since such a spectrum will lead to ill-conditioned systems at high contrast media, this Calderón preconditioning approach will not be considered here.

\section{Results}
\label{sec:results}

This section provides the results of computational benchmarks that validate the efficiency of the boundary integral formulations at high-contrast media.

\subsection{Settings}

Let us first summarise the computational settings of the benchmarks.

\subsubsection{Boundary integral formulations}

Five different boundary integral formulations will be considered for the computational benchmarks: the standard PMCHWT and Müller formulations, the Calderón preconditioned PMCHWT formulation, and the novel high-contrast formulations, as summarized in Table~\ref{table:formulations}. Since the high-contrast formulations are based on an indirect representation of the fields, only half the number of boundary integral operators need to be assembled, compared to the direct formulations. Notice that the computational complexity of the standard formulations can be improved by exploiting the symmetry between the double-layer and adjoint double-layer operators.

The high-contrast Neumann and Dirichlet formulations require less operations for each matrix-vector multiplication. That is, for a single matrix-vector multiplication of the entire system, the number of individual matrix-vector multiplications of boundary integral operators is $\ell + 3\ell^2$. Differently, the standard formulations require $4\ell + 4\ell^2$ dense operations for each multiplication of the system matrix. Furthermore, Calderón preconditioning doubles the cost of a matrix-vector multiplication. While operator sums can be calculated explicitly to reduce the time for each matrix-vector multiplication, this cannot be performed anymore for fast arithmetic with the fast multipole method or hierarchical matrix compression. Furthermore, operator products are never calculated explicitly. Instead, separate matrix-vector multiplications are performed. Finally, the presence of identity operations in the boundary integral formulations does not incur significant computation time since these are sparse operators.

\begin{table}[!ht]
	\centering
	\caption{The boundary integral formulations for scattering at multiple penetrable objects, with the number of dense boundary integral operators to be assembled and the number of dense matrix-vector multiplications, where $\ell$ denotes the number of scatterers.}
	\label{table:formulations}
	\begin{tabular}{lcccccc}
		\hline\hline
		formulation & \#operators & \#matvecs \\
		\hline
		High-contrast Neumann \eqref{eq:bie:highcontrast:neu} & $2\ell + 2\ell^2$ & $\ell + 3\ell^2$ \\
		High-contrast Dirichlet \eqref{eq:bie:highcontrast:dir} & $2\ell + 2\ell^2$ & $\ell + 3\ell^2$ \\
		PMCHWT \eqref{eq:pmchwt} & $4\ell + 4\ell^2$ & $4\ell + 4\ell^2$ \\
		Müller \eqref{eq:muller} & $4\ell + 4\ell^2$ & $4\ell + 4\ell^2$ \\
		Calderón PMCHWT \eqref{eq:pmchwt:calderon} & $4\ell + 4\ell^2$ & $8\ell + 8\ell^2$ \\
		\hline\hline
	\end{tabular}
\end{table}

\subsubsection{Numerical parameters}

All boundary integral formulations are implemented with version~3 of the open-source BEMPP library~\cite{smigaj2015solving, scroggs2017software, betcke2020product}. The meshes are generated with Gmsh~\cite{geuzaine2009gmsh} and MeshLab~\cite{meshlab}. The spectra and condition numbers are calculated with dense matrix assembly, while hierarchical matrix compression with a tolerance of $10^{-5}$ is used when solving the system. The linear solver is GMRES with a relative tolerance of $10^{-7}$ as convergence criterion, without restart, and implemented with the library SciPy~\cite{virtanen2020scipy}.

\subsubsection{Material parameters}

No attenuation will be considered for the acoustic propagation. Hence, the wavenumber is given by
\begin{equation}
	k = 2\pi f / c
\end{equation}
for a given frequency~$f$ in Hz and a wavespeed $c$ in m/s.
See Table~\ref{table:parameters:physical} for characteristic values of materials commonly found in acoustical engineering.

\begin{table}[!ht]
	\centering
	\caption{Physical parameters of the scattering media~\cite{duck1990physical, itis2018, gray1963american, wohletz1992volcanology}.}
	\label{table:parameters:physical}
	\begin{tabular}{lrr}
		\hline\hline
		material & $\rho$ [kg/m$^3$] & $c$ [m/s] \\
		\hline
		air & 1.225 & 340 \\
		fat & 917 & 1412 \\
		water & 1025 & 1500 \\
		bone & 1912 & 4080 \\
		basalt & 2740 & 3350 \\
		iron & 7725 & 4094 \\
		\hline\hline
	\end{tabular}
\end{table}

\begin{figure}[!ht]
	\centering
	\includegraphics[width=0.45\columnwidth]{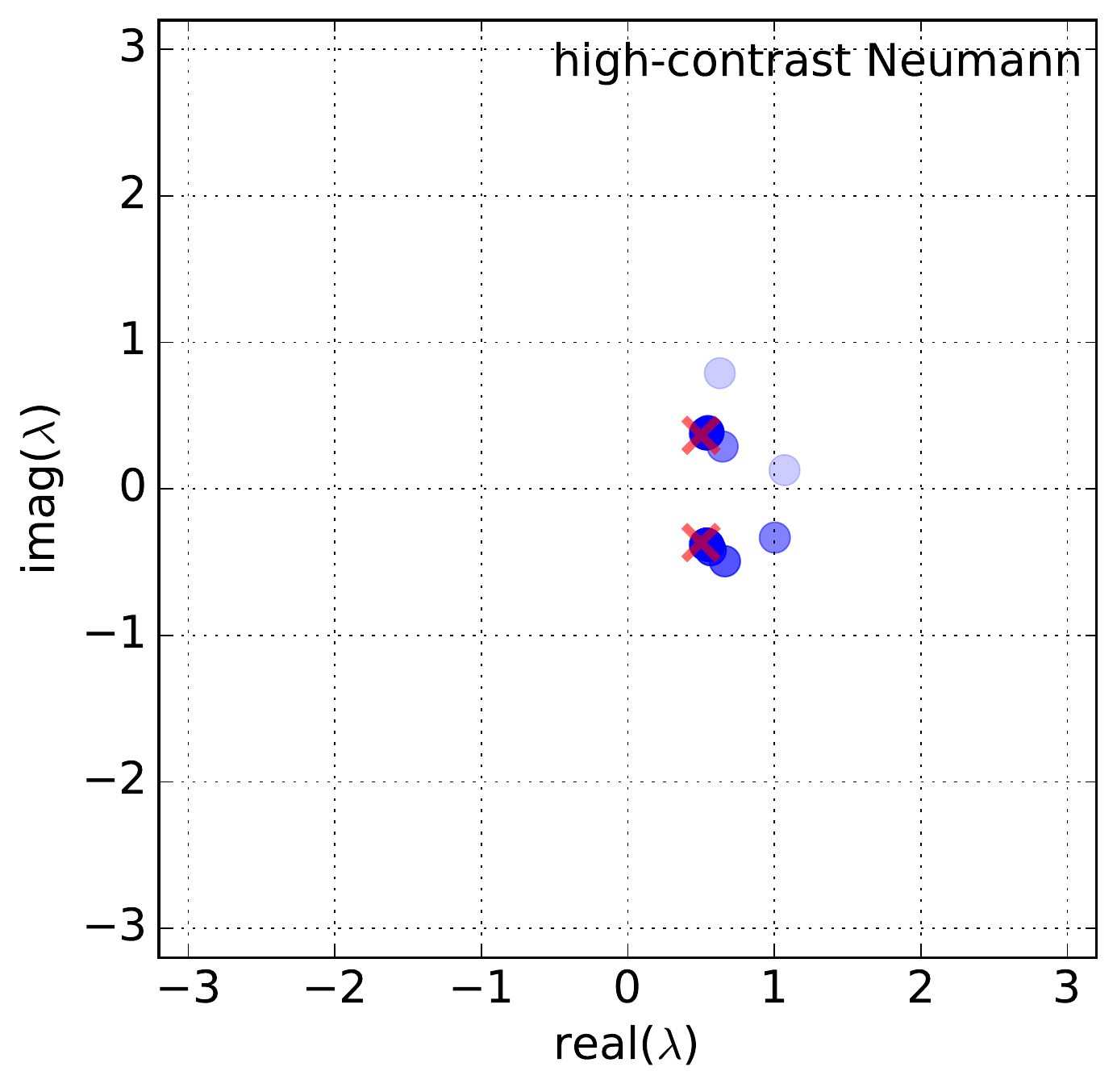}
	\includegraphics[width=0.45\columnwidth]{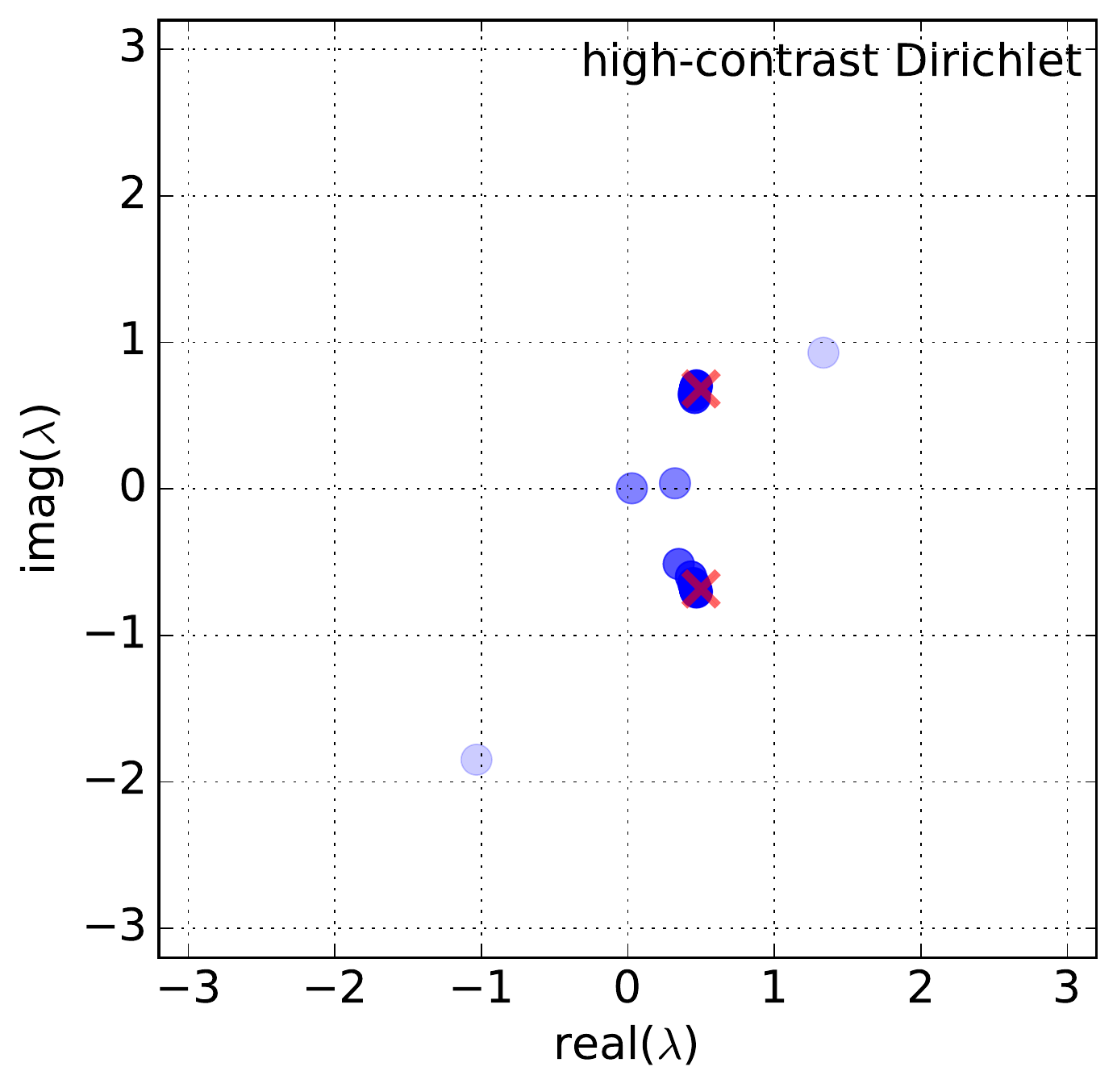}
	\includegraphics[width=0.32\columnwidth]{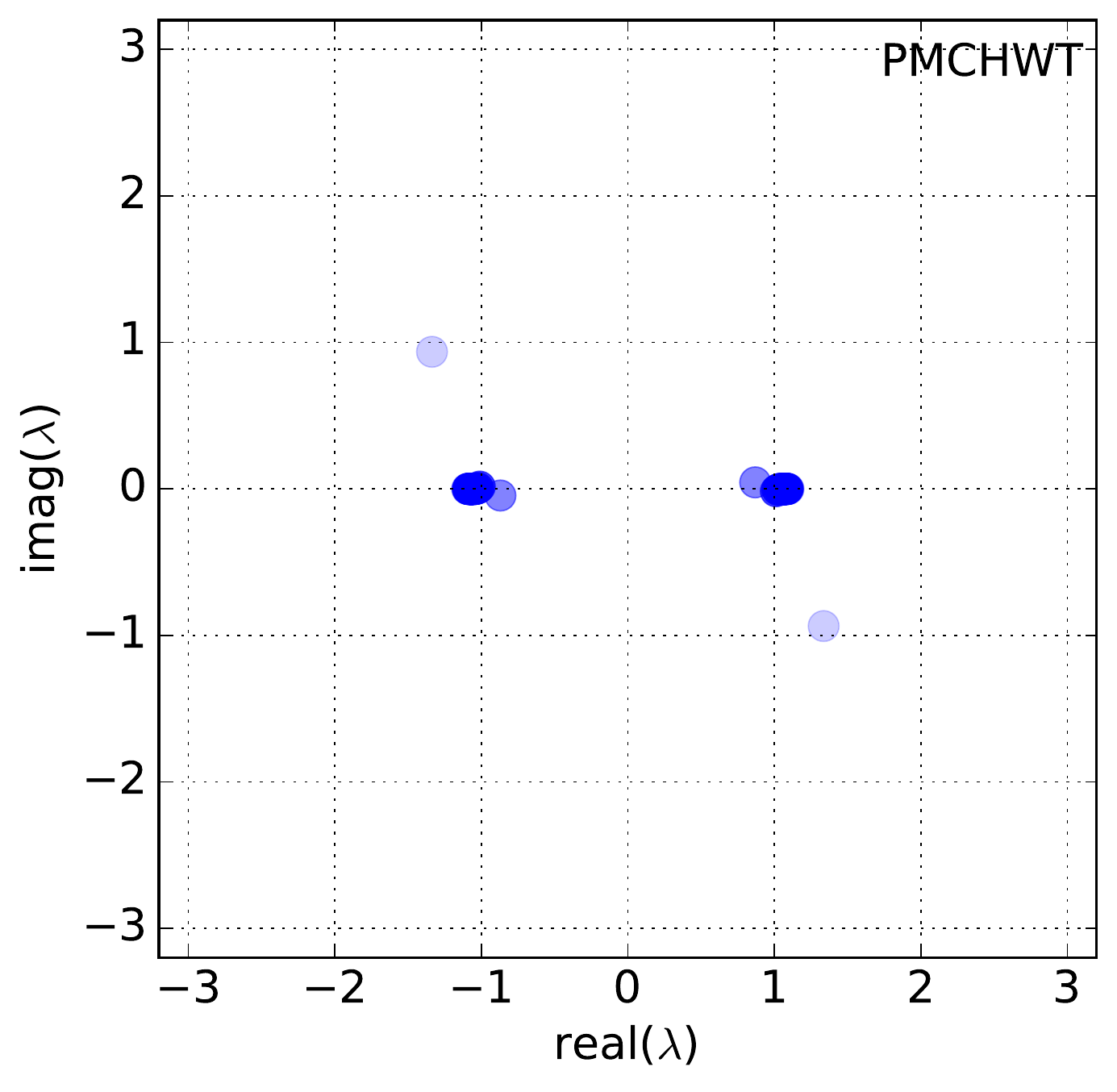}
	\includegraphics[width=0.32\columnwidth]{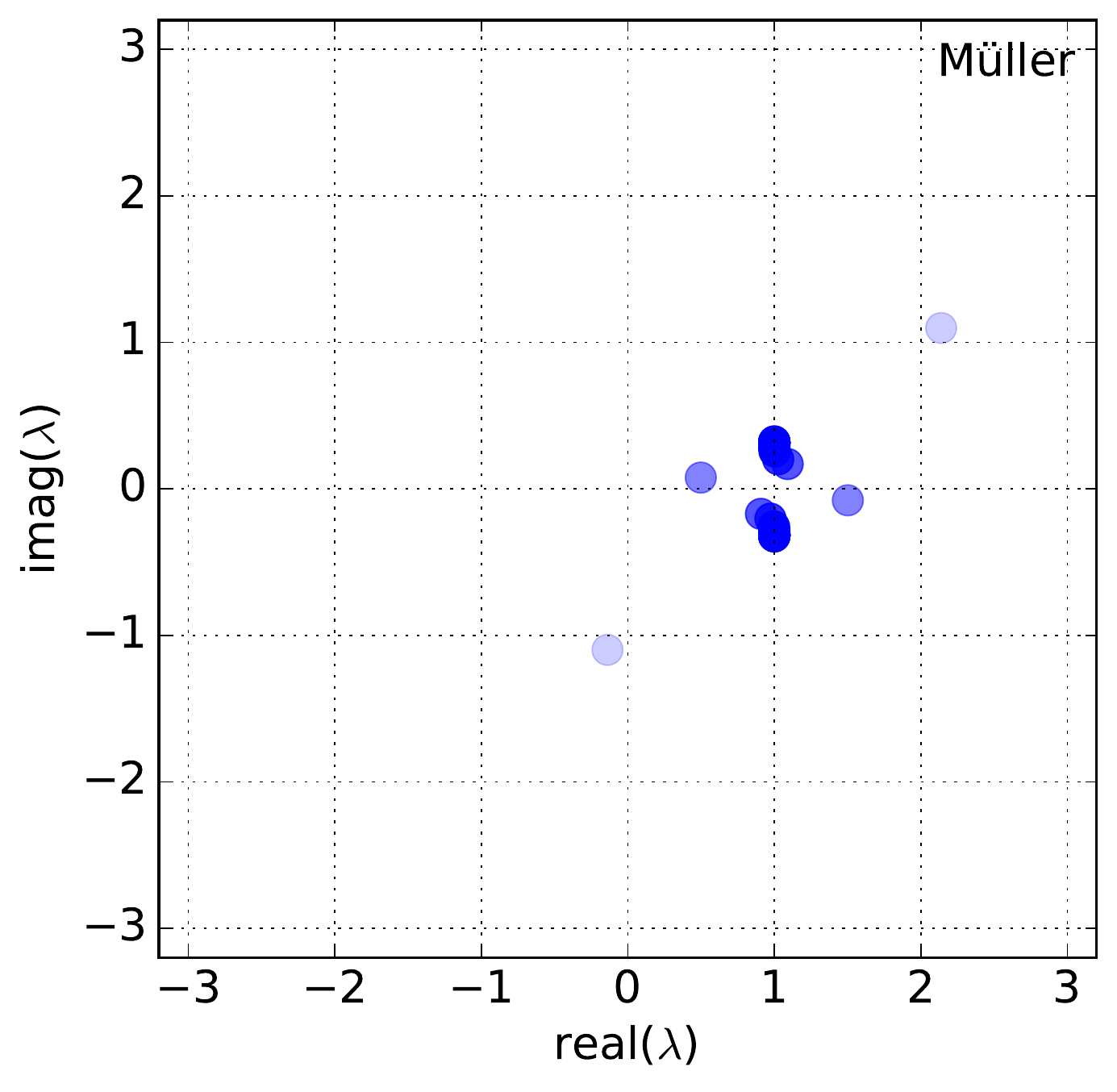}
	\includegraphics[width=0.32\columnwidth]{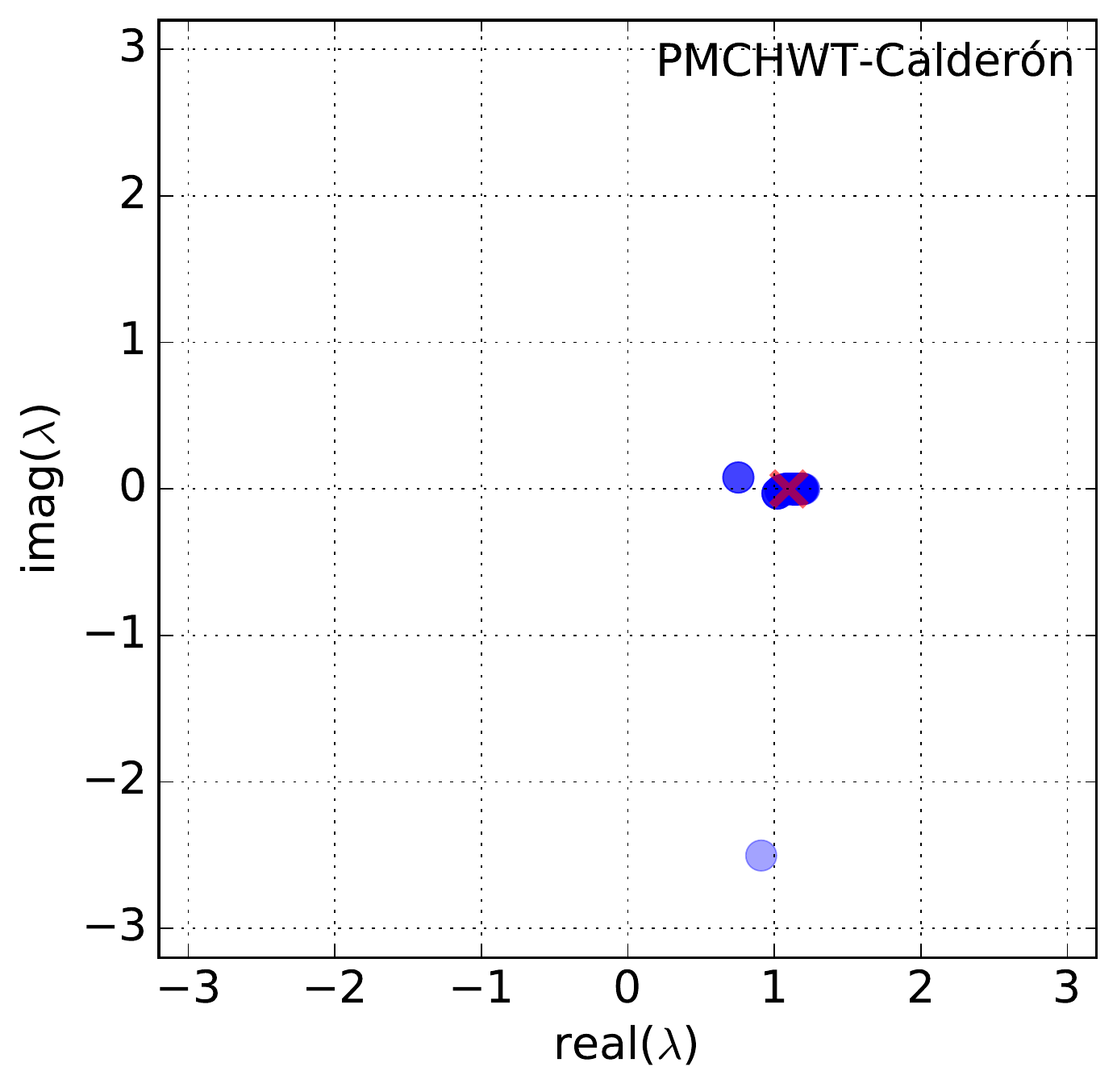}
	\caption{The eigenvalues of the boundary integral formulations for exterior water and interior bone. The frequency is 500~Hz, $\kext=2.09$, $\kint=0.77$, the geometry is a sphere with unit radius, and the surface mesh has 177~nodes. The red crosses depict the expected accumulation points.}
	\label{fig:spectrum:formulations}
\end{figure}

\subsection{Spectrum of formulations}
\label{sec:results:spectrum}

The high-contrast Neumann and Dirichlet formulations and the Calderón preconditioned PMCHWT formulation have eigenvalue accumulation points that depend on the density ratio, as given by Eqns.~\eqref{eq:accumulation:neumann}, \eqref{eq:accumulation:dirichlet}, and~\eqref{eq:accumulation:calderon}, respectively. For this purpose, let us consider a spherical domain and calculate the eigenvalues of the system matrix of the boundary integral formulations. The density ratio between water and bone at the interface is moderately high. The results in Figure~\ref{fig:spectrum:formulations} clearly show the accumulation of the eigenvalues at the expected points. Furthermore, this benchmark also suggests clustering of eigenvalues for the standard PMCHWT and Müller formulations.

\subsection{Conditioning with density contrast}

The accumulation points of the high-contrast materials directly depend on the density ratio across a material interface. For this reason, the condition number of the system matrix is expected to depend on the density ratio as well. Let us benchmark the influence of the density ratio on the condition number by considering a spherical object. Two benchmarks will be performed, one where the wavespeed remains constant with changing density and one where the wavespeed is related to the density by a constant material compressibility. In the following, the relative density will be defined as $\rho_-/\rho_+$.

\begin{figure}[!ht]
	\centering
	\includegraphics[width=0.99\columnwidth]{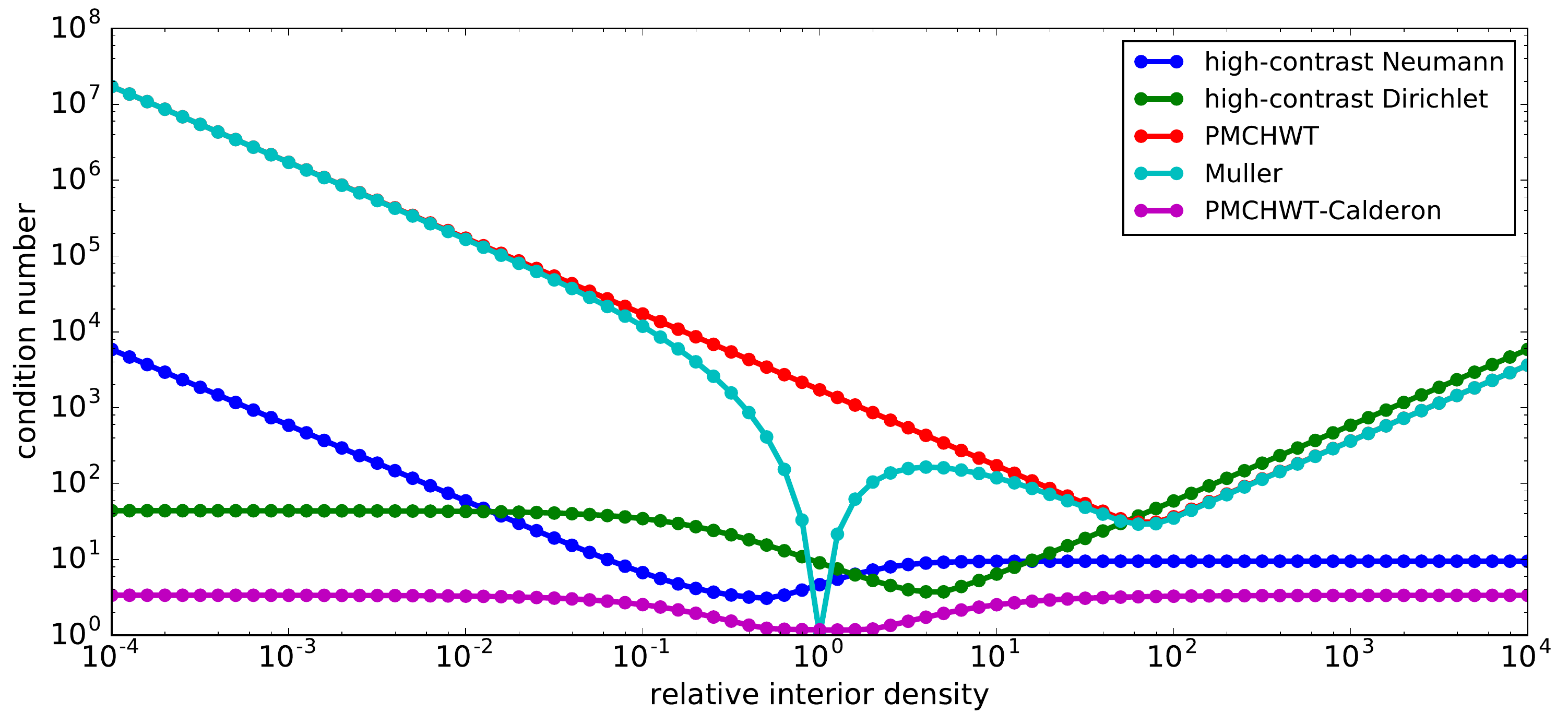}
	\caption{The condition number of the boundary integral formulations with respect to the relative density of the interior material: $\rho_-/\rho_+$. The surface mesh on the unit sphere has 784~nodes, which holds 8 elements per wavelength at a frequency of 1250~Hz. The wavespeed is constant and $\kext = \kint = 5.24$.}
	\label{fig:cond:rho:c}
\end{figure}

\begin{figure}[!ht]
	\centering
	\begin{subfigure}[b]{\columnwidth}
		\centering
		\includegraphics[width=0.34\columnwidth]{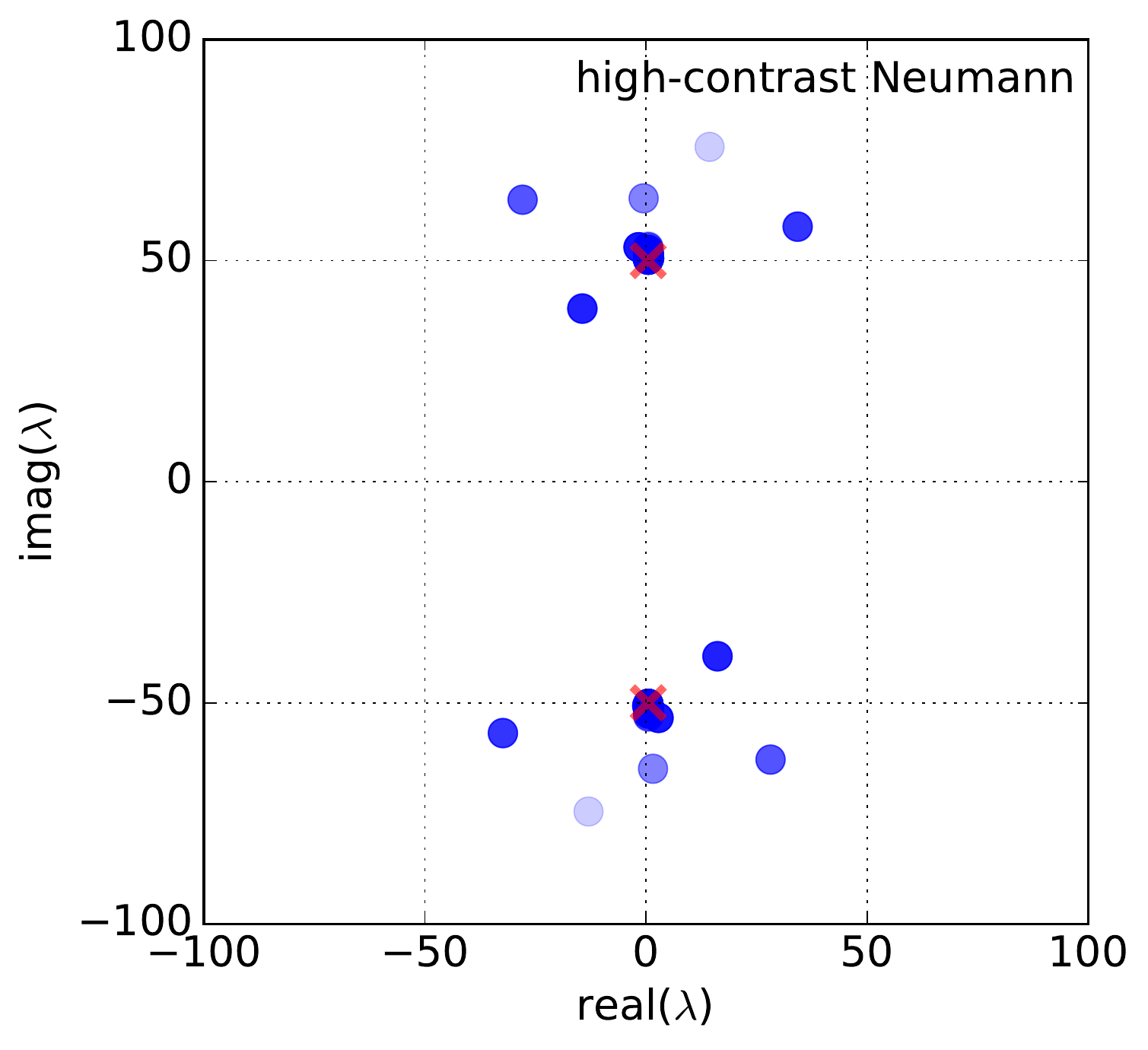}
		\includegraphics[width=0.34\columnwidth]{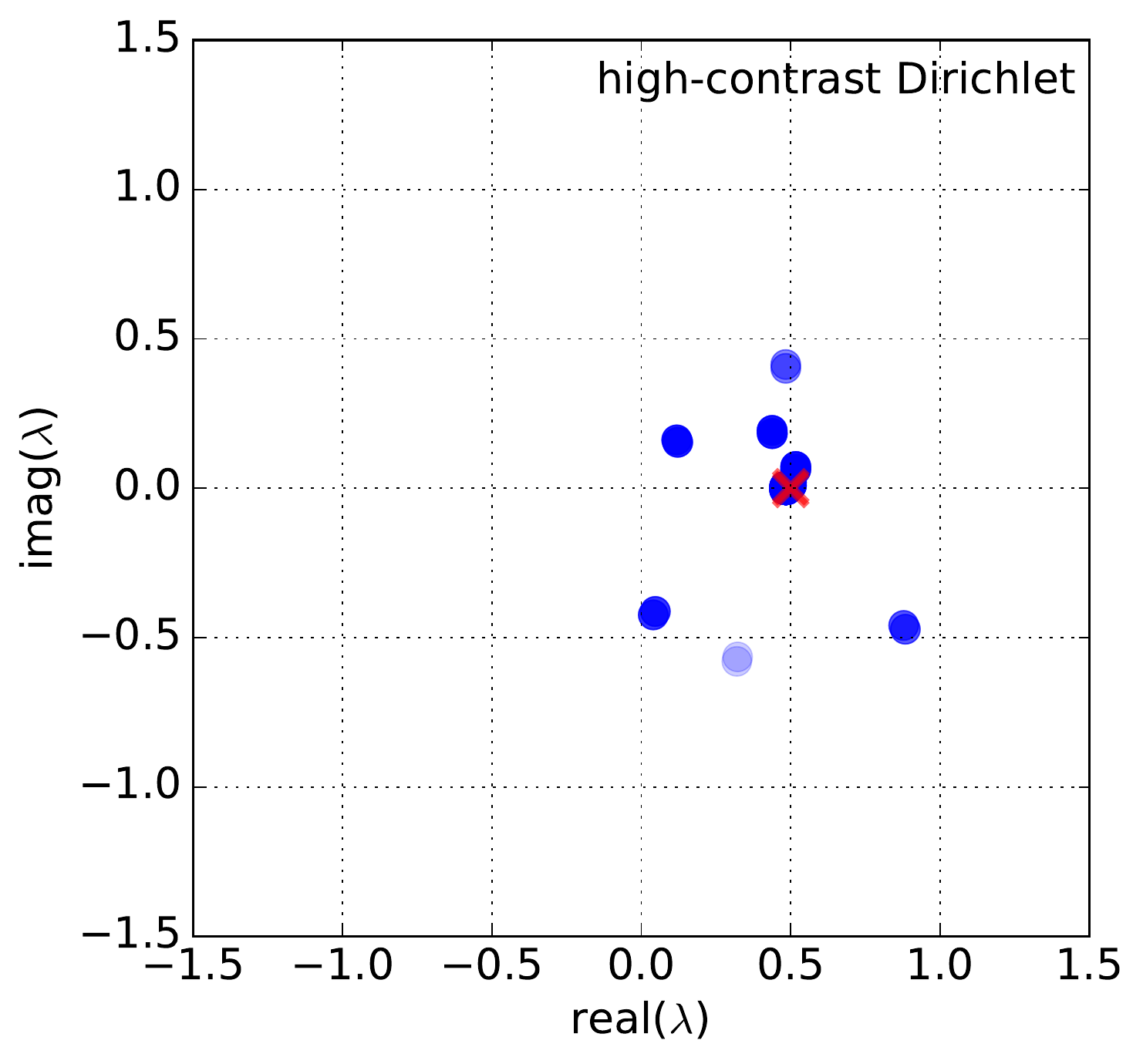}
		\includegraphics[width=0.32\columnwidth]{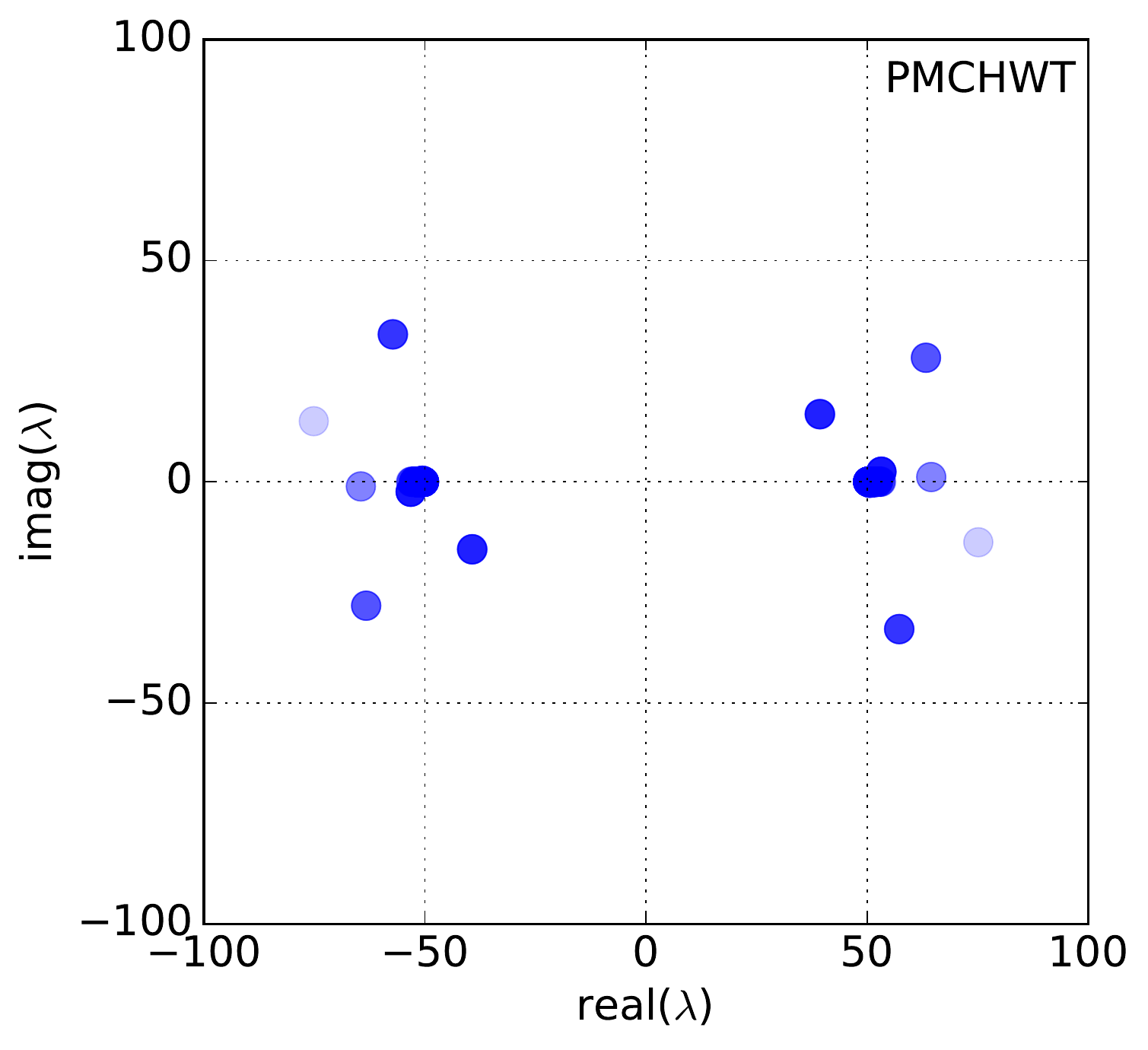}
		\includegraphics[width=0.32\columnwidth]{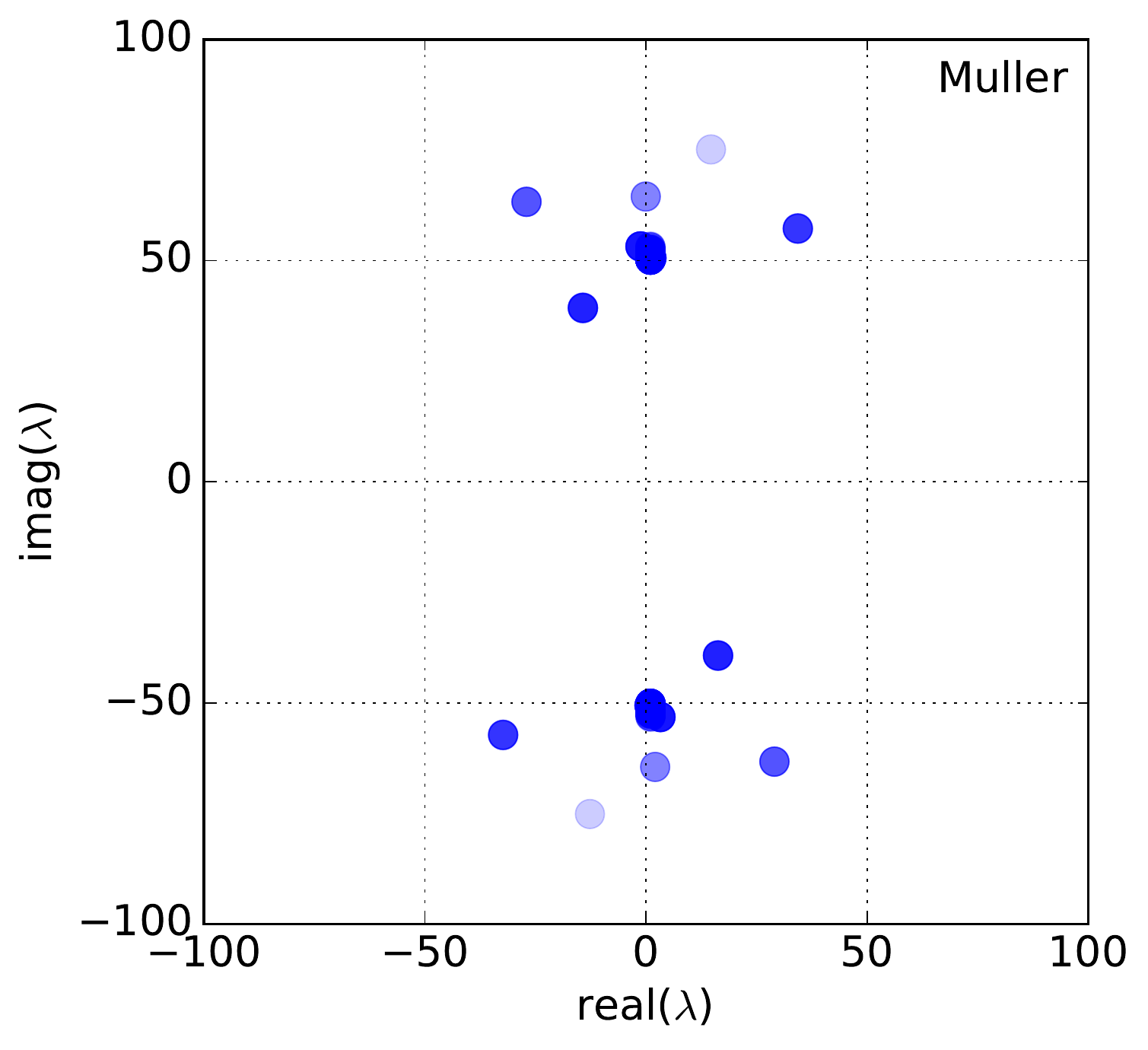}
		\includegraphics[width=0.32\columnwidth]{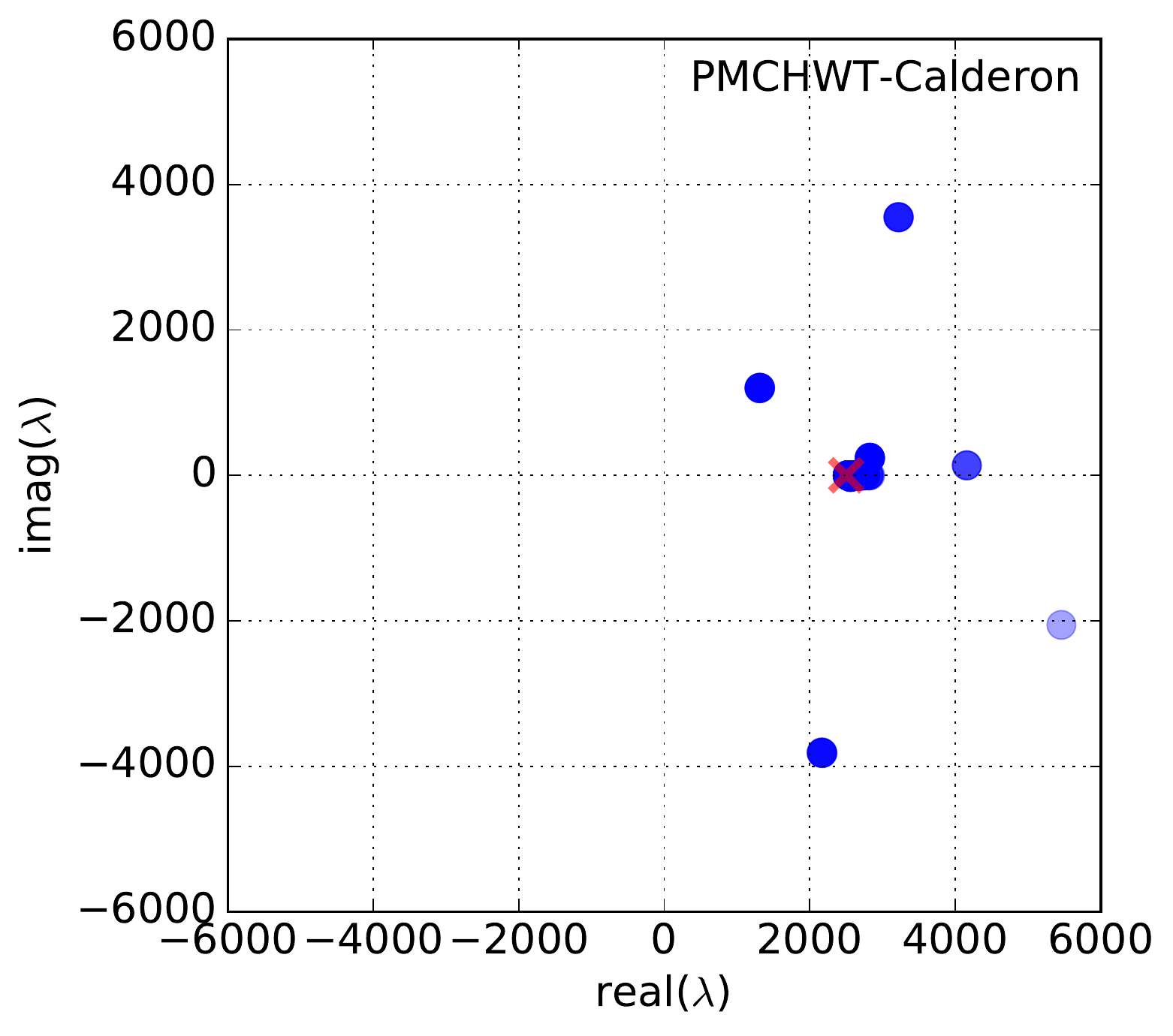}
		\caption{Density ratio $\rho_-/\rho_+ = 10^{-4}$.}
		\label{fig:spectrum:rho:c:0}
	\end{subfigure}
	\begin{subfigure}[b]{\columnwidth}
		\centering
		\includegraphics[width=0.34\columnwidth]{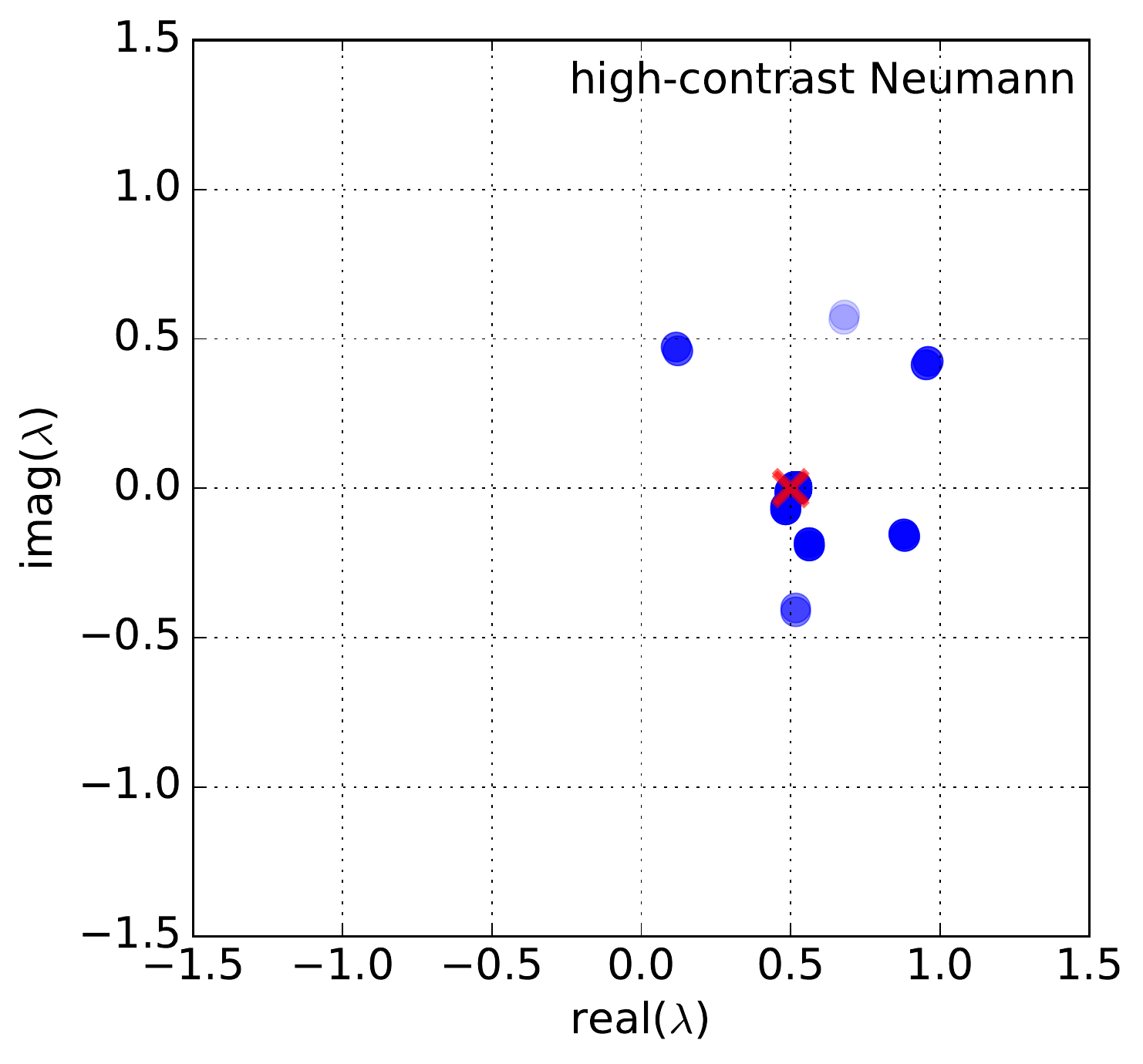}
		\includegraphics[width=0.34\columnwidth]{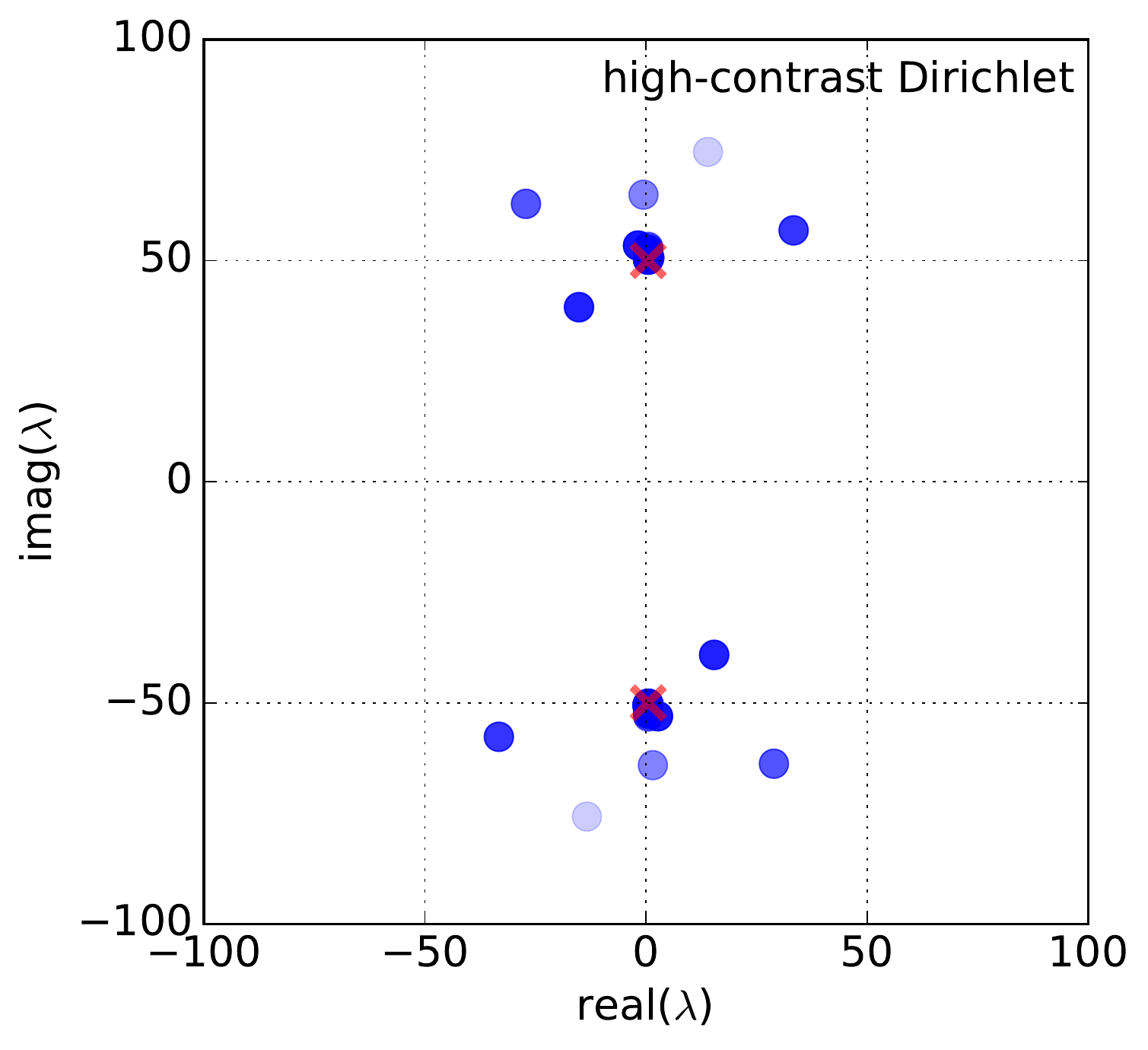}
		\includegraphics[width=0.32\columnwidth]{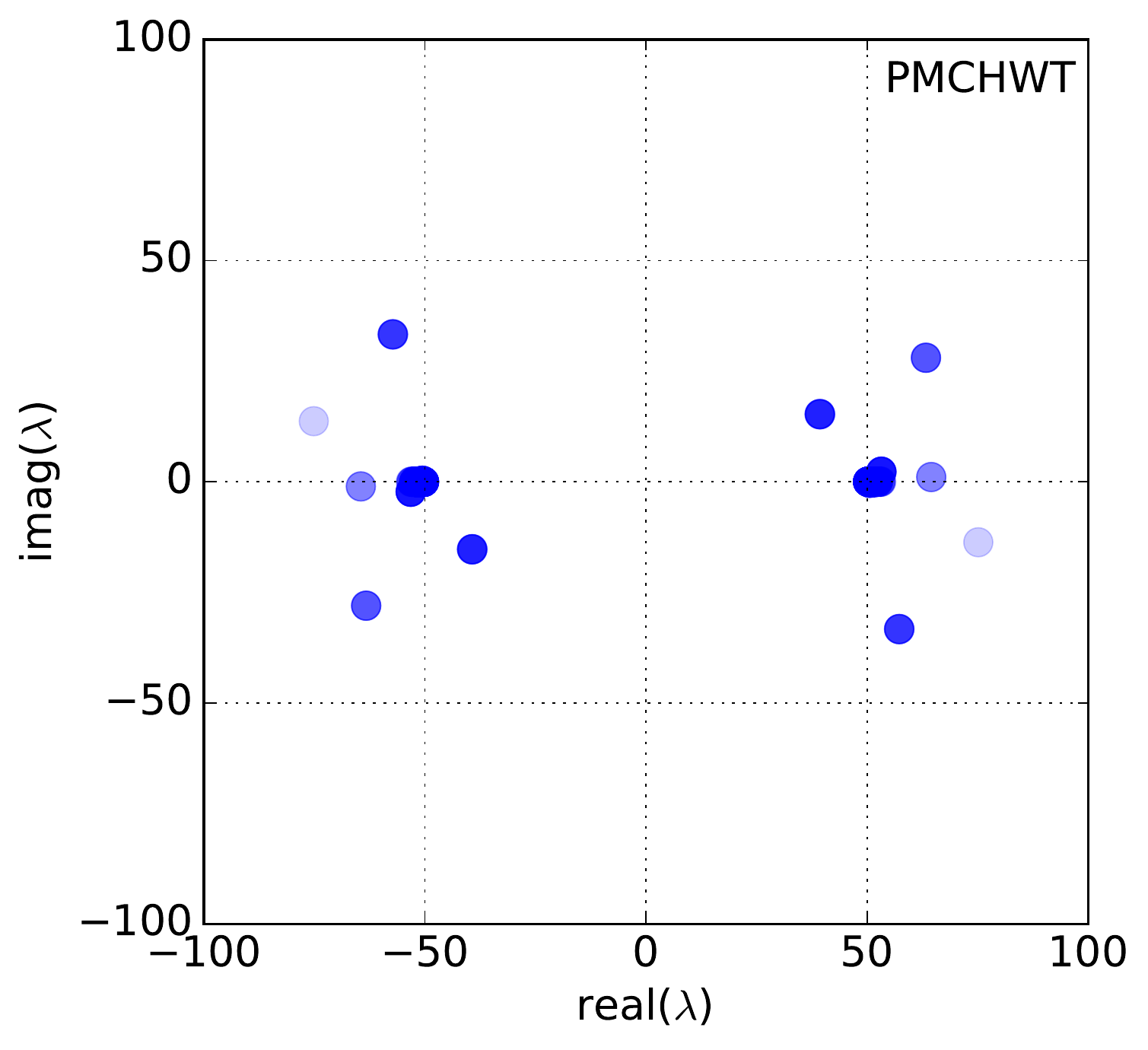}
		\includegraphics[width=0.32\columnwidth]{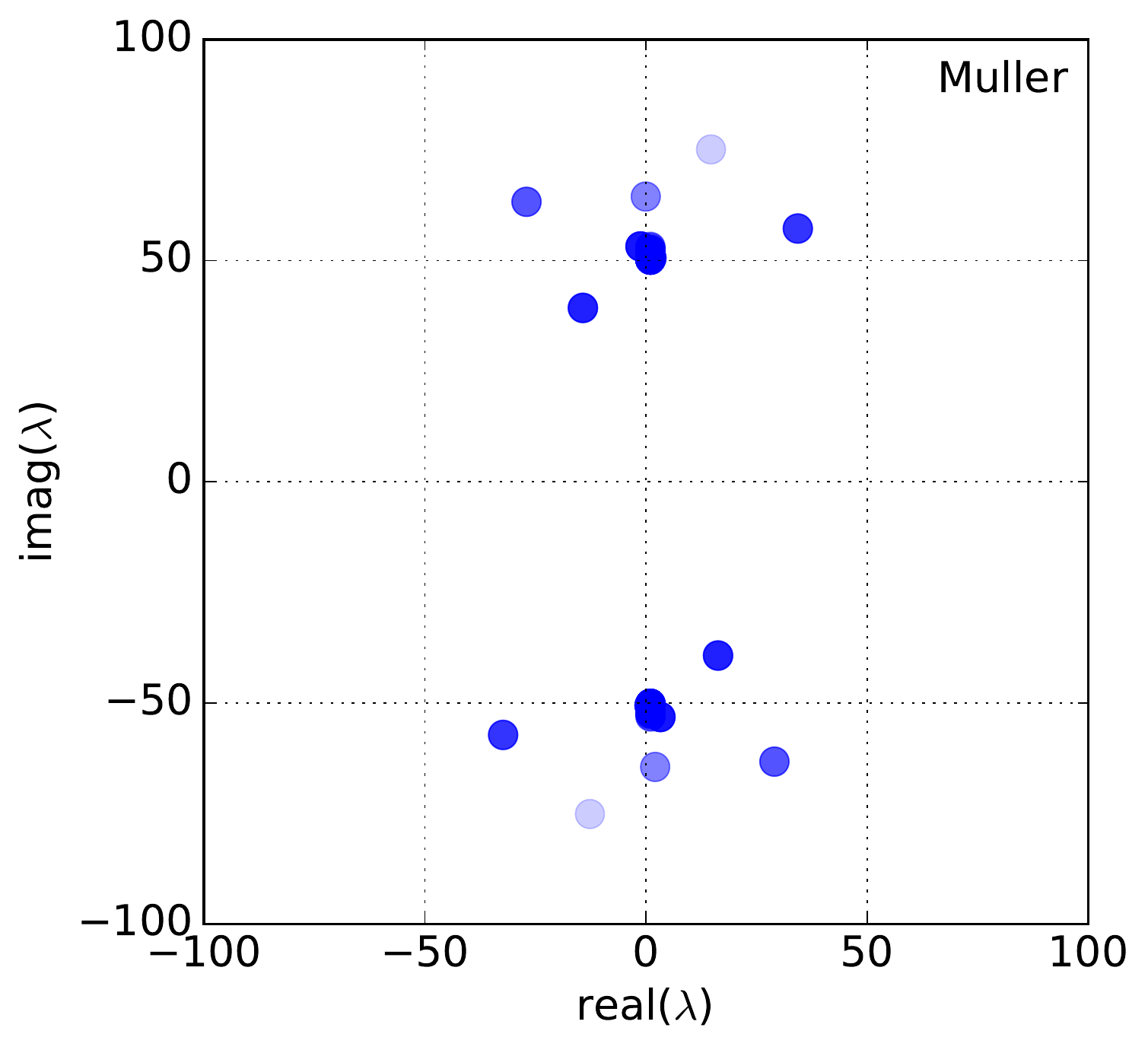}
		\includegraphics[width=0.32\columnwidth]{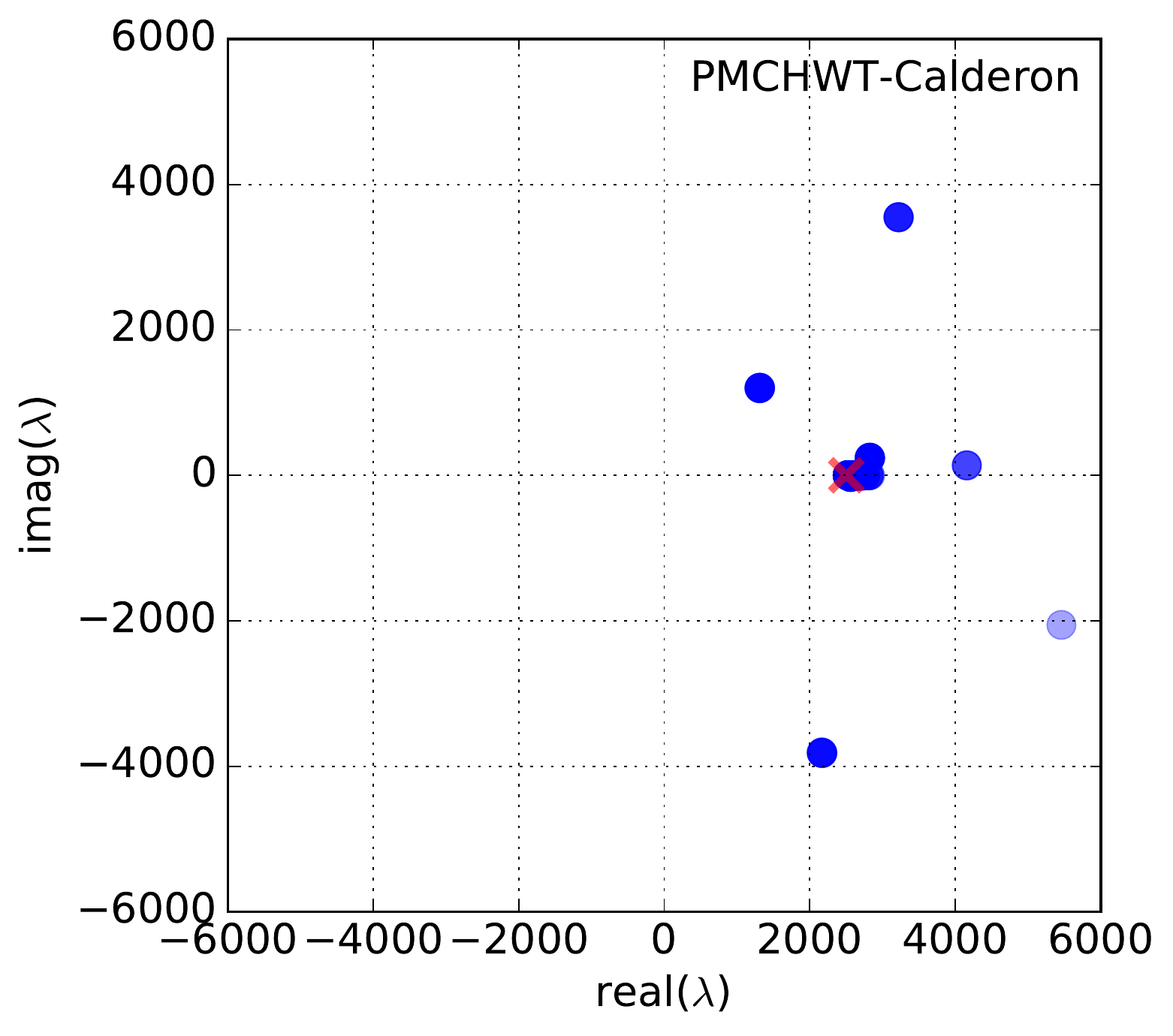}
		\caption{Density ratio $\rho_-/\rho_+ = 10^{4}$.}
		\label{fig:spectrum:rho:c:8}
	\end{subfigure}
	\caption{The eigenvalues of the boundary integral formulations for the benchmark described in Figure~\ref{fig:cond:rho:c}. The red crosses depict the expected accumulation points.}
	\label{fig:spectrum:rho:c}
\end{figure}

First, let us consider a benchmark where the interior density changes and the exterior density remains constant. The surface mesh and frequency are fixed. In both materials, the wavespeed will be fixed to $c=1500$ (resembling water), thus yielding a constant wavenumber across the entire benchmark.
Figure~\ref{fig:cond:rho:c} depicts the condition number of the boundary integral formulations, for a wide range of interior densities. Figure~\ref{fig:spectrum:rho:c} depicts the spectra for the high-contrast cases of $\rho_-/\rho_+ = 10^{-4}$ and  $\rho_-/\rho_+ = 10^{4}$.

The benchmarks show that the PMCHWT and Müller formulations become ill-conditioned for high density contrasts. Differently, the Calderón preconditioned PMCHWT formulation remains well-conditioned for the entire range of density contrasts. This behaviour is expected since this formulation has a single accumulation point. Moreover, since the interior and exterior wavenumbers are equal, the cross terms of the operator products~\eqref{eq:cp:terms} cancel out exactly.

Regarding the high-contrast formulations, the Neumann version remains well-conditioned for high interior densities while the Dirichlet version remains well-conditioned for low interior densities. This behaviour is consistent with the spectral analysis that shows a set of two accumulation points that depend on the density ratio, see Eq.~\eqref{eq:accumulation:neumann}. For the high-contrast Neumann formulation, both accumulation points converge towards $1/2$ when the interior density is relatively high while these two accumulation points diverge when the interior density is relatively low. The accumulation point of the Dirichlet version behaves in opposite direction, that is, the two accumulation points converge towards $1/2$ for relatively low interior densities and diverge for relatively high interior densities.

\begin{figure}[!ht]
	\centering
	\includegraphics[width=0.99\columnwidth]{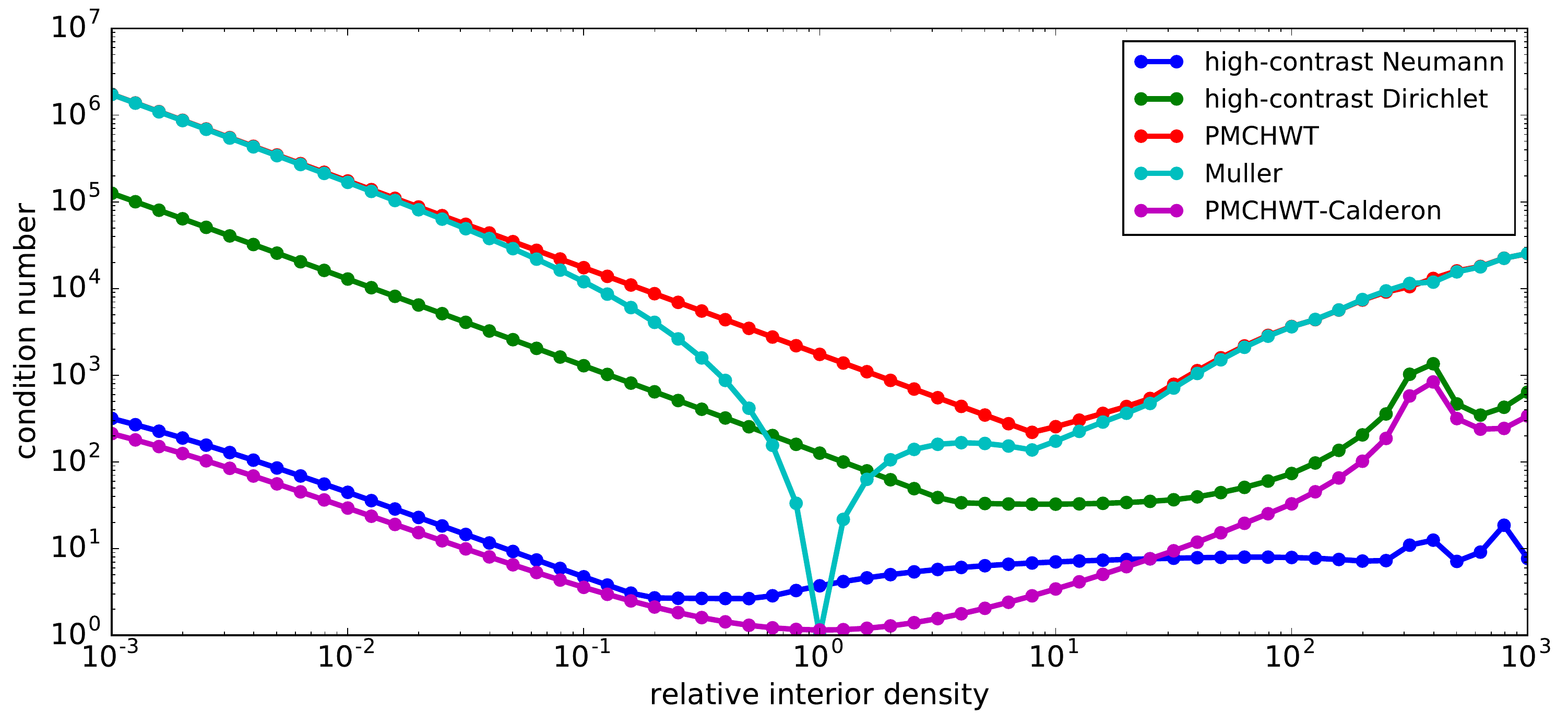}
	\caption{The condition number of the boundary integral formulations with respect to the relative density of the interior material: $\rho_-/\rho_+$. The surface mesh of the sphere has 784~nodes and the frequency of 39.5~Hz is chosen such that at least 8 elements per wavelength are present in all cases. The wavespeed depends on the density with a constant compressibility so that $\kext = 0.166$ and $0.00524 \leq \kint \leq 5.24$.}
	\label{fig:cond:rho:beta}
\end{figure}

\begin{figure}[!ht]
	\centering
	\begin{subfigure}[b]{\columnwidth}
		\centering
		\includegraphics[width=0.34\columnwidth]{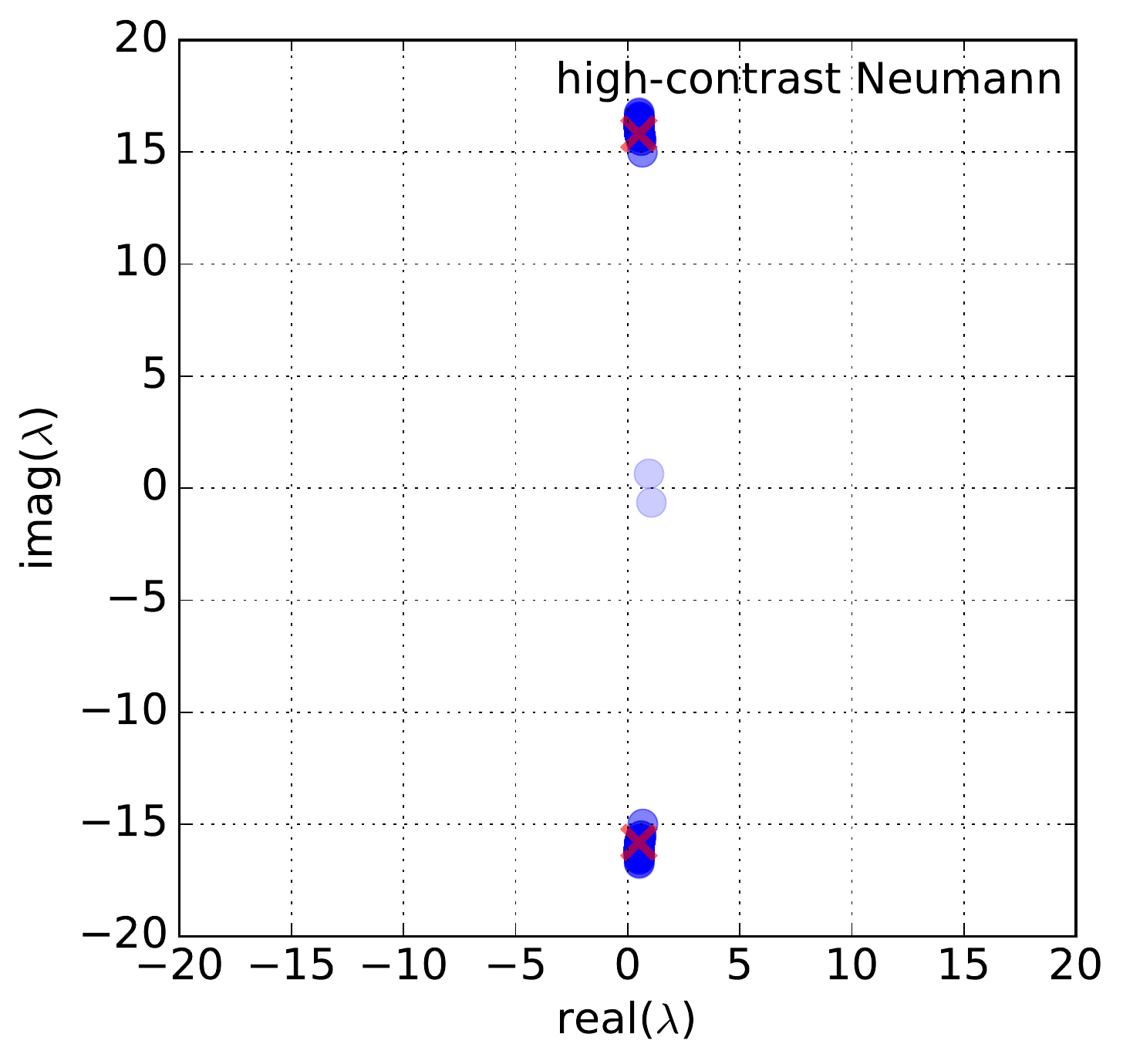}
		\includegraphics[width=0.34\columnwidth]{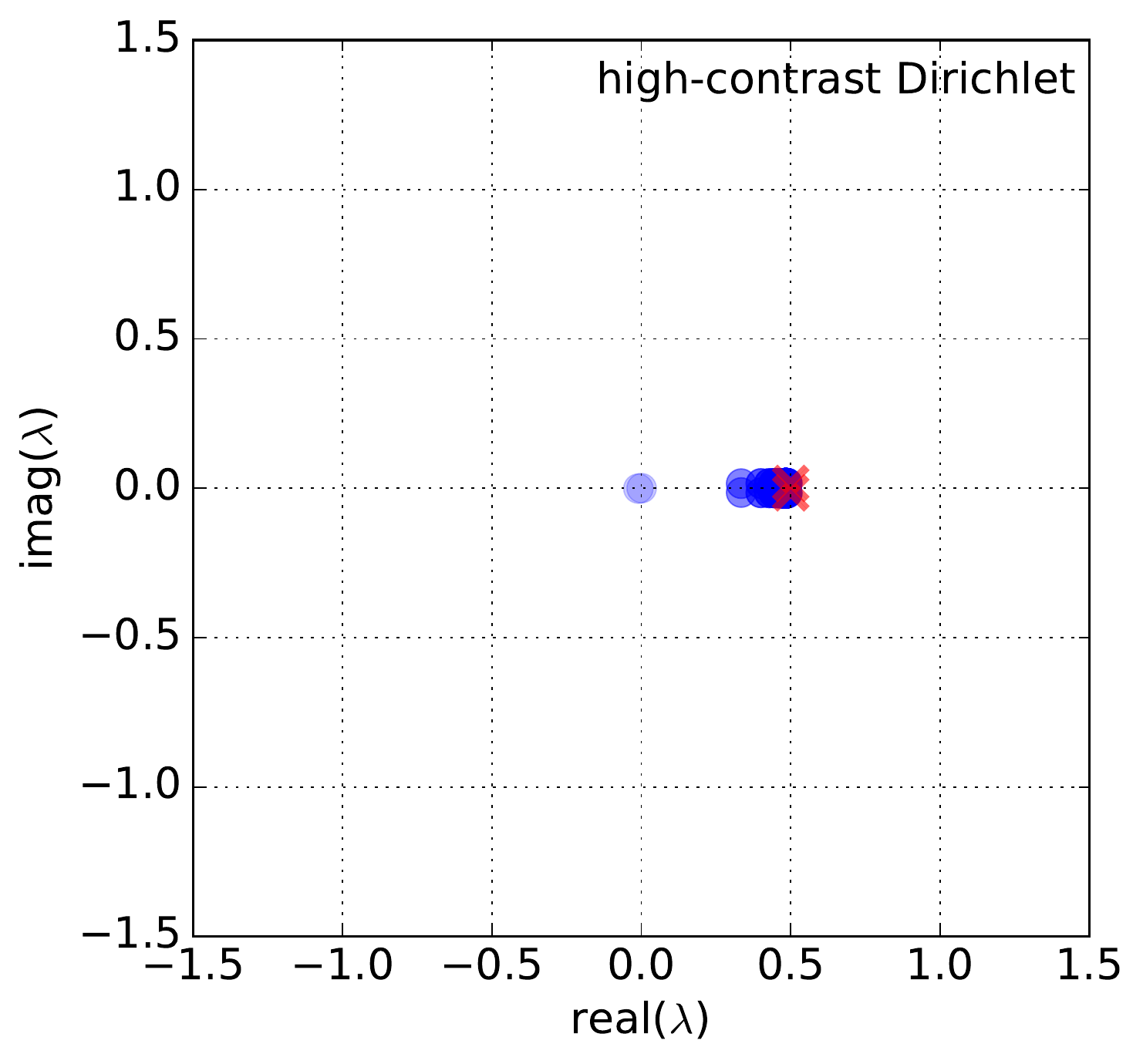}
		\includegraphics[width=0.32\columnwidth]{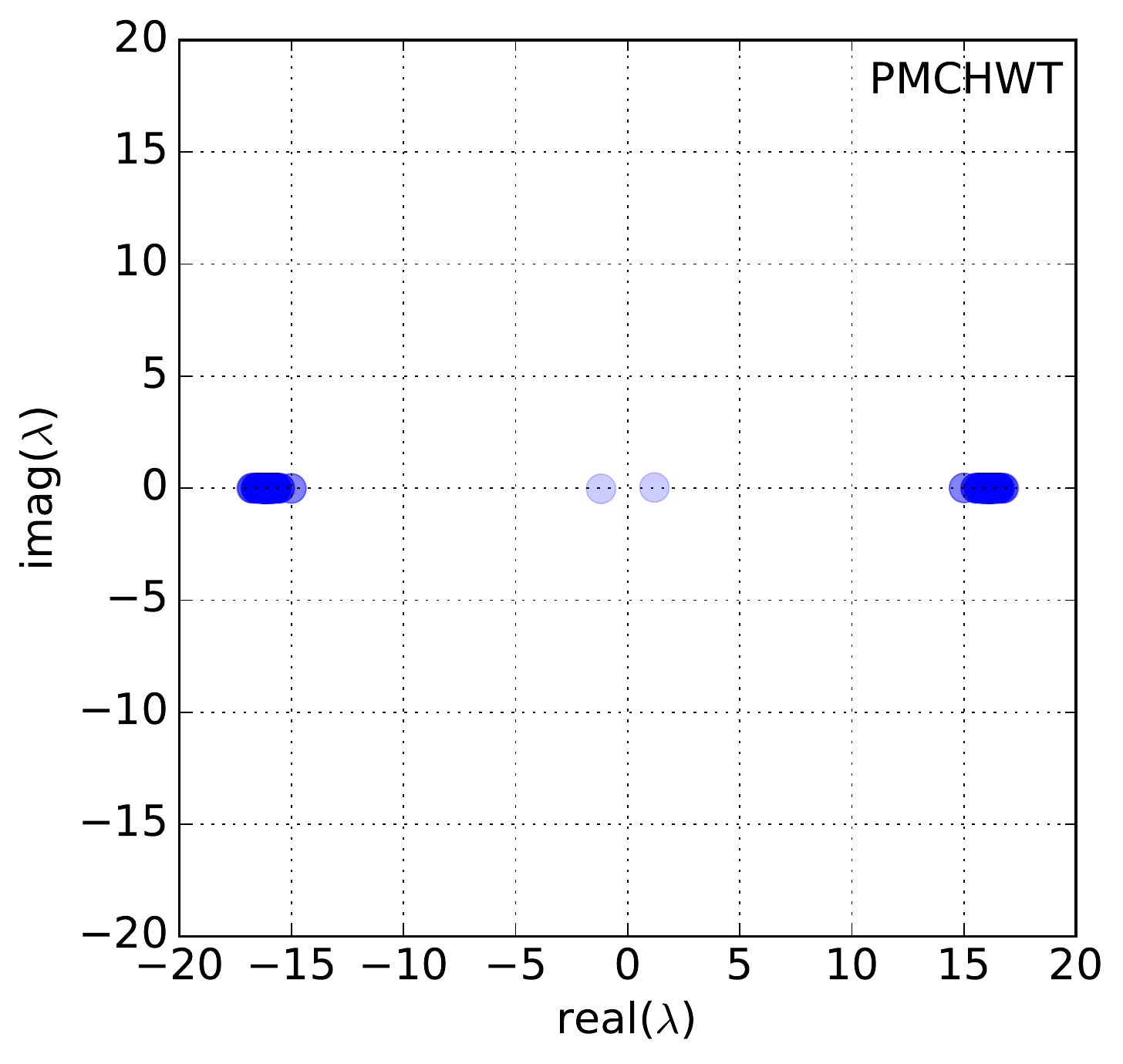}
		\includegraphics[width=0.32\columnwidth]{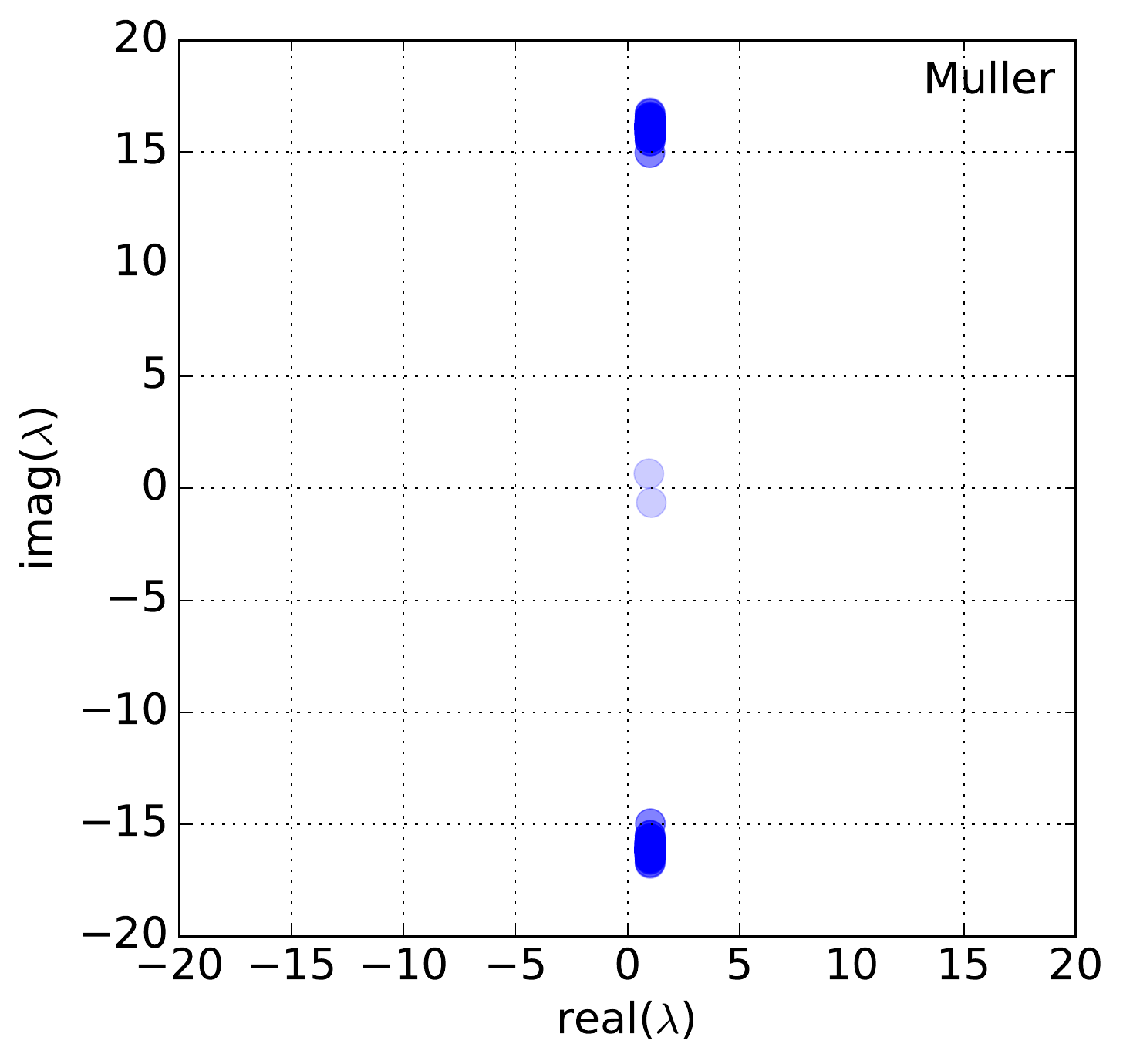}
		\includegraphics[width=0.32\columnwidth]{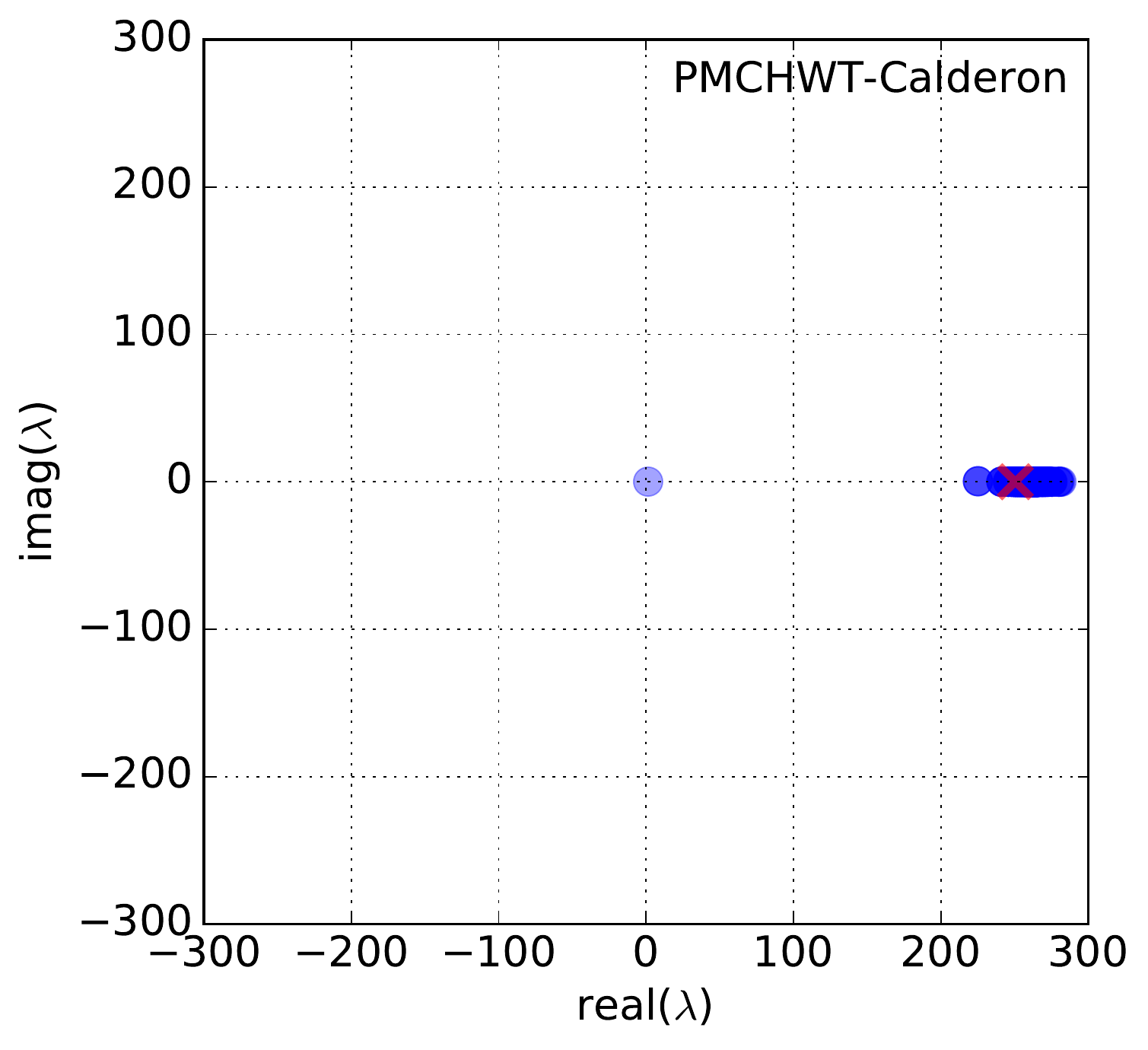}
		\caption{Density ratio $\rho_-/\rho_+ = 10^{-3}$.}
		\label{fig:spectrum:rho:beta:0}
	\end{subfigure}
	\begin{subfigure}[b]{\columnwidth}
		\centering
		\includegraphics[width=0.34\columnwidth]{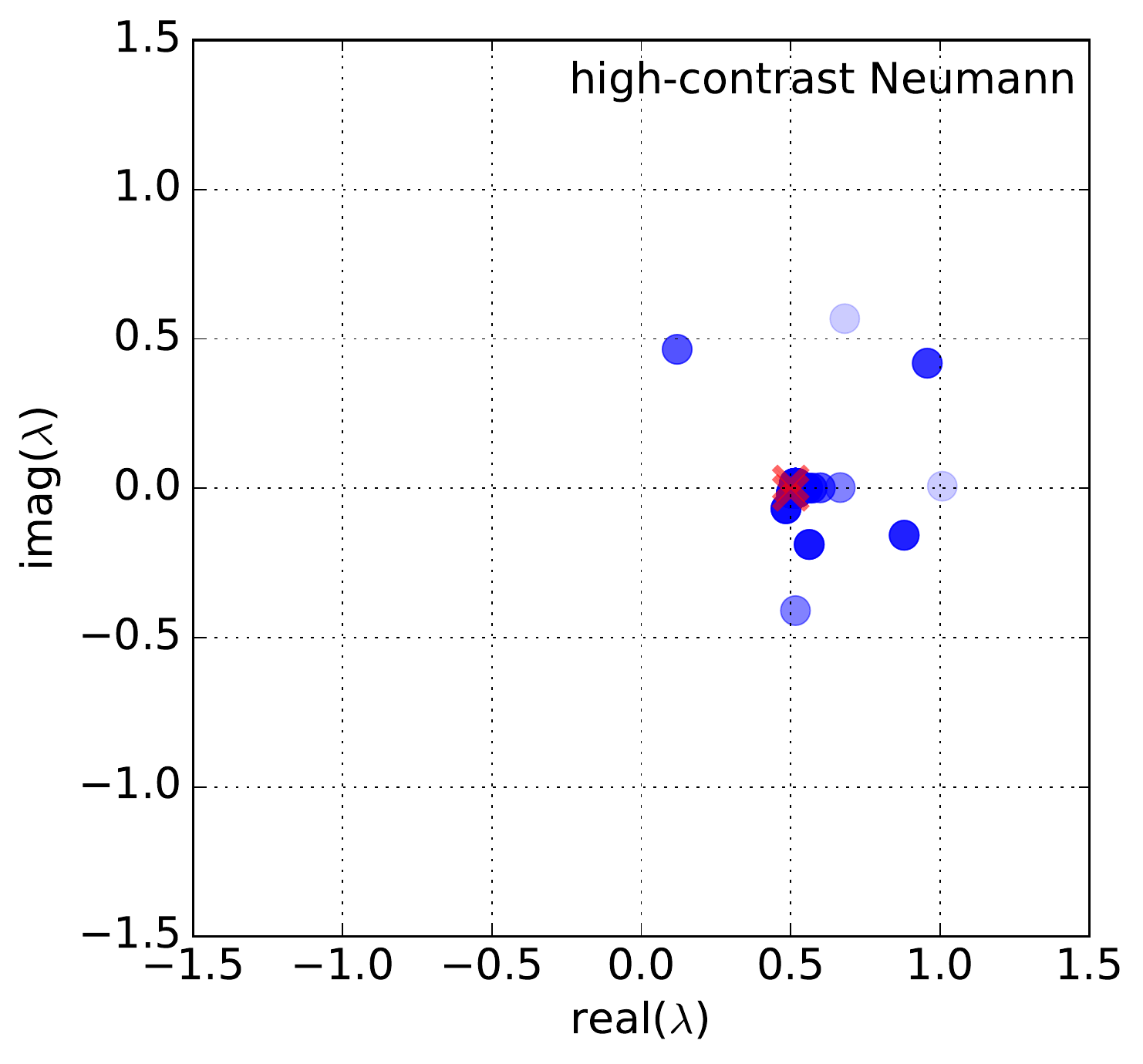}
		\includegraphics[width=0.34\columnwidth]{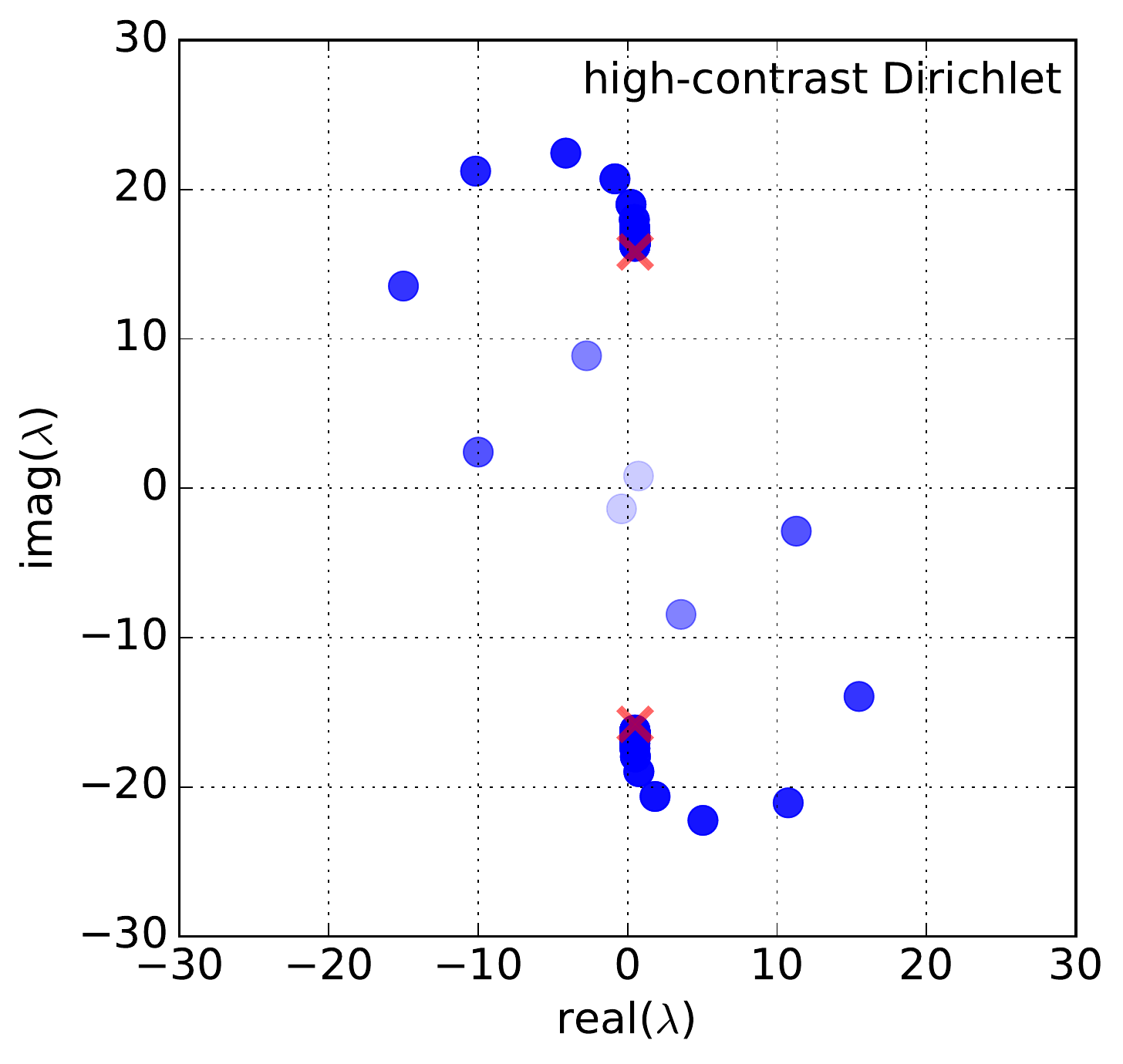}
		\includegraphics[width=0.32\columnwidth]{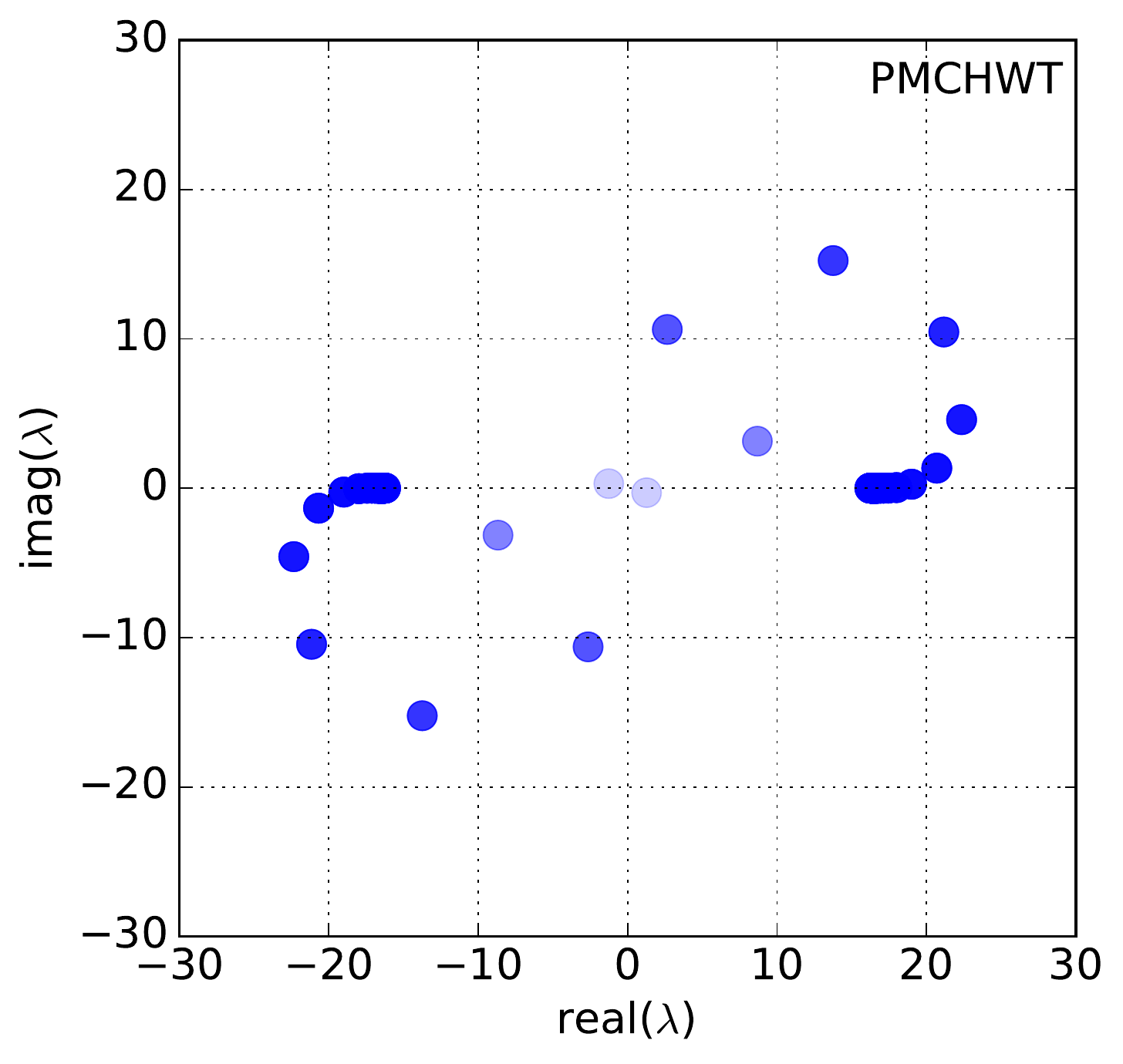}
		\includegraphics[width=0.32\columnwidth]{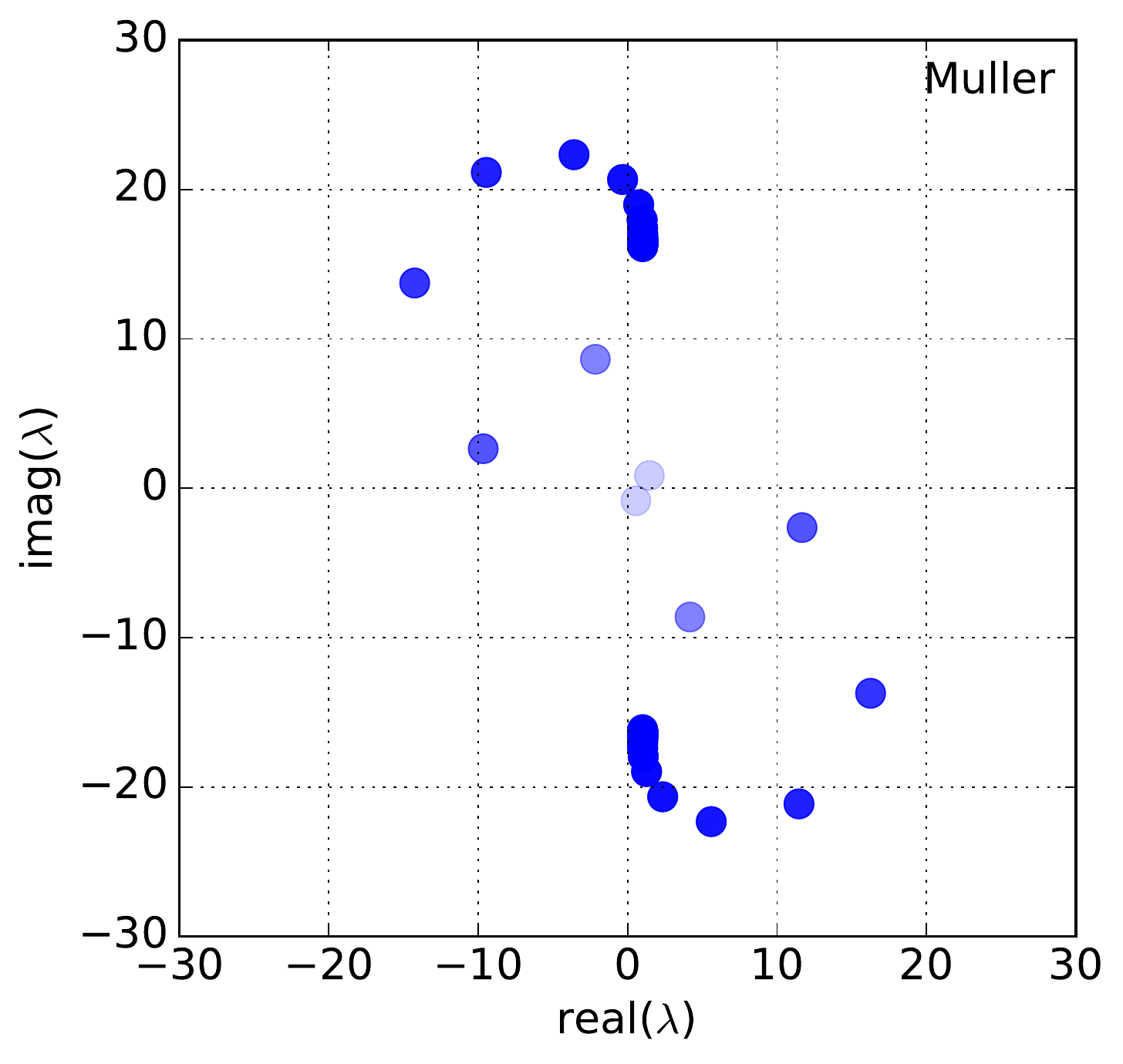}
		\includegraphics[width=0.32\columnwidth]{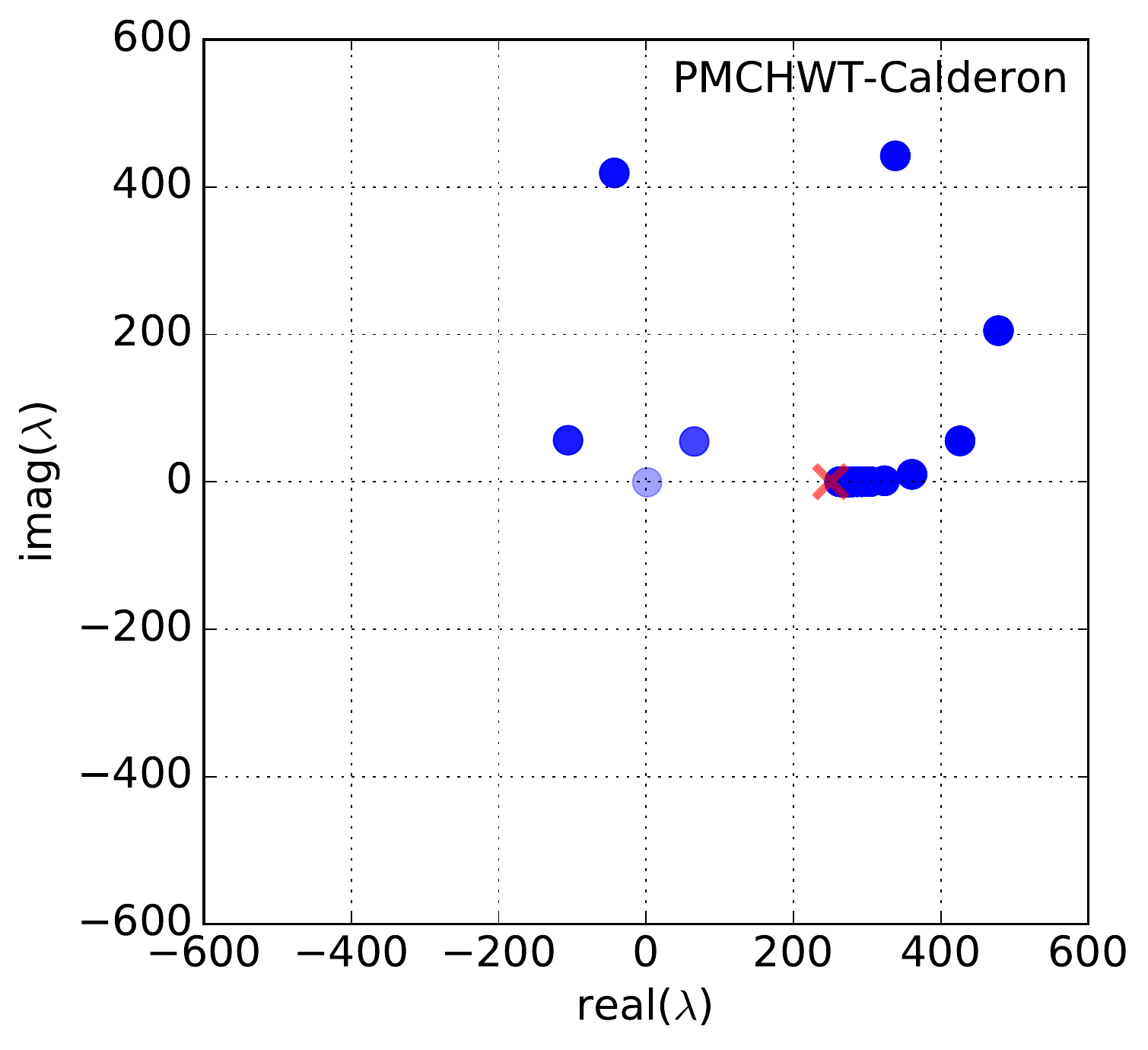}
		\caption{Density ratio $\rho_-/\rho_+ = 10^{3}$.}
		\label{fig:spectrum:rho:beta:6}
	\end{subfigure}
	\caption{The eigenvalues of the boundary integral formulations for the benchmark described in Figure~\ref{fig:cond:rho:beta}. The red crosses depict the expected accumulation points.}
	\label{fig:spectrum:rho:beta}
\end{figure}

The previous benchmark assumed a constant wavespeed for changing density. This has physical limitations since the acoustic wavespeed in a material depends on the mass density. A commonly used relation that models this dependency is by considering
\begin{equation}
	\label{eq:wavespeed}
	c = \frac1{\sqrt{\beta \rho}}
\end{equation}
where $\beta$ denotes the compressibility of the material~\cite{morse1986theoretical}. The second benchmark for the influence of the interior density uses a constant compressibility. Hence, the density influences the wavenumber as well. All physical parameters are taken to resemble water, and the mesh and frequency are fixed. The mesh was generated with at least eight triangles per wavelength for the highest frequencies in the benchmark.
Figure~\ref{fig:cond:rho:beta} presents the conditioning of the boundary integral formulations and Figure~\ref{fig:spectrum:rho:beta} the spectra for the highest density contrasts.

As before, the standard PMCHWT and Müller formulations become ill-conditioned at high density ratios, as expected. Differently than before, the other formulations also show deterioration of the conditioning at high-contrast media, even though much better conditioned than the standard formulations. Looking at the spectra, clustering of the eigenvalues is still visible but now with eigenvalues close to zero as well. Since the accumulation points of the eigenvalues depend on the densities only, the worse conditioning has to be attributed to a higher contrast in wavenumber. The operator products present in both the Calderón preconditioned PMCHWT formulation and the high-contrast formulations diverge from the Calderón identities that are only valid for zero material contrast. The only formulation that keeps a constant condition number with respect to density ratio is the high-contrast Neumann formulation when the interior material has a higher density than the exterior.

\subsection{Conditioning with frequency}

The results presented above clearly demonstrate a strong influence of the density ratio on the conditioning of the discretisation matrix, as well as the ratio in wavespeed. Another physical parameter that strongly influences the conditioning of the boundary integral formulation is the frequency of the wave field. The following benchmarks will demonstrate the influence of the frequency and material contrast in terms of the condition number of the system matrix and the number of iterations required for the GMRES iterative solver. Different materials will be used, with a low-contrast ratio between water and fat, an intermediate contrast ratio between water and bone, and a high-contrast situation with air and iron materials. See Table~\ref{table:parameters:physical} for the physical parameters.

\begin{figure}[!ht]
	\centering
	\begin{subfigure}[b]{0.99\columnwidth}
		\includegraphics[width=\textwidth]{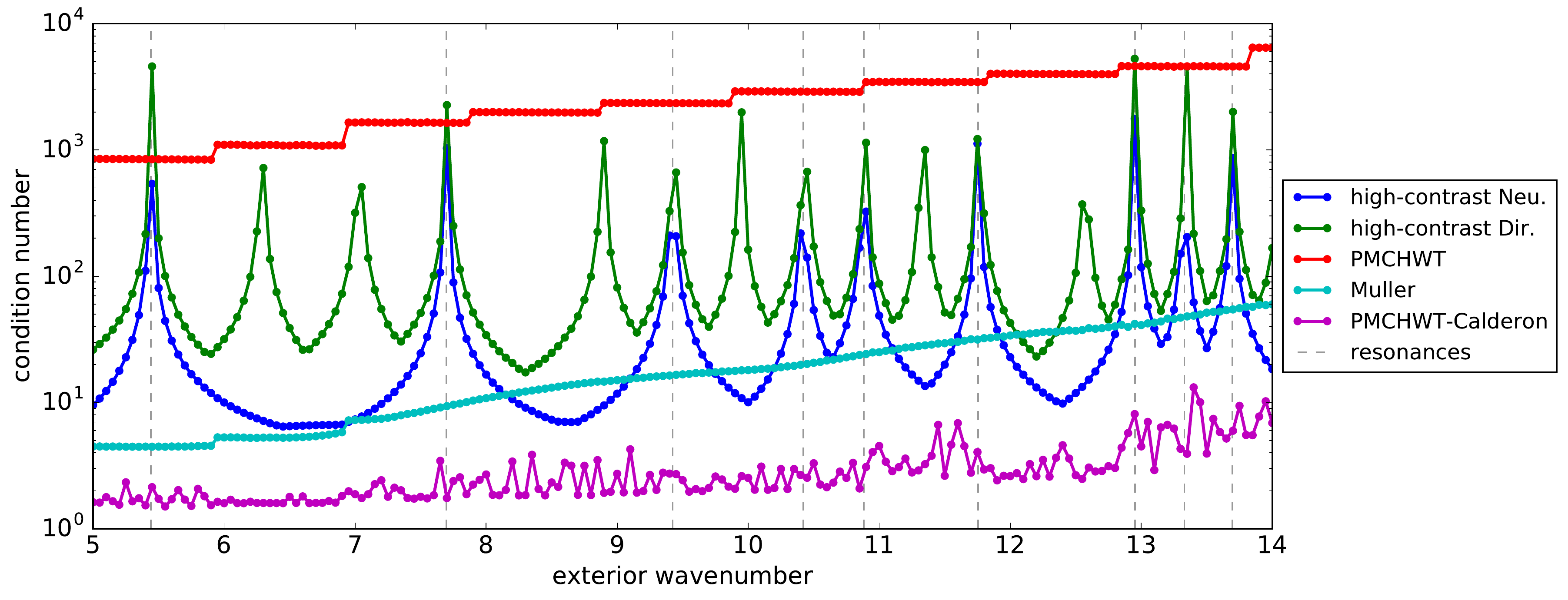}
		\caption{Exterior water and interior fat ($\kint/\kext=1.06$, $\rhoint/\rhoext=0.895$).}
	\end{subfigure}
	\begin{subfigure}[b]{0.99\columnwidth}
		\includegraphics[width=\textwidth]{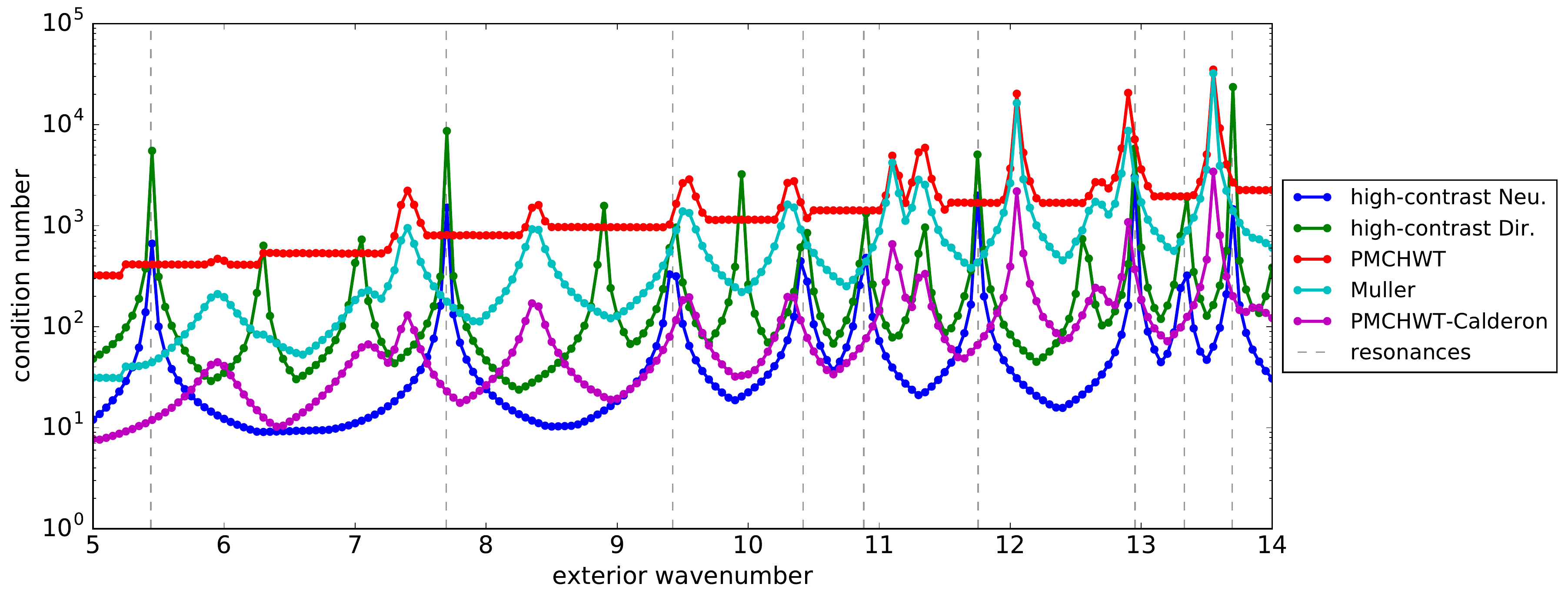}
		\caption{Exterior water and interior bone ($\kint/\kext=0.368$, $\rhoint/\rhoext=1.87$).}
	\end{subfigure}
	\begin{subfigure}[b]{0.99\columnwidth}
		\includegraphics[width=\textwidth]{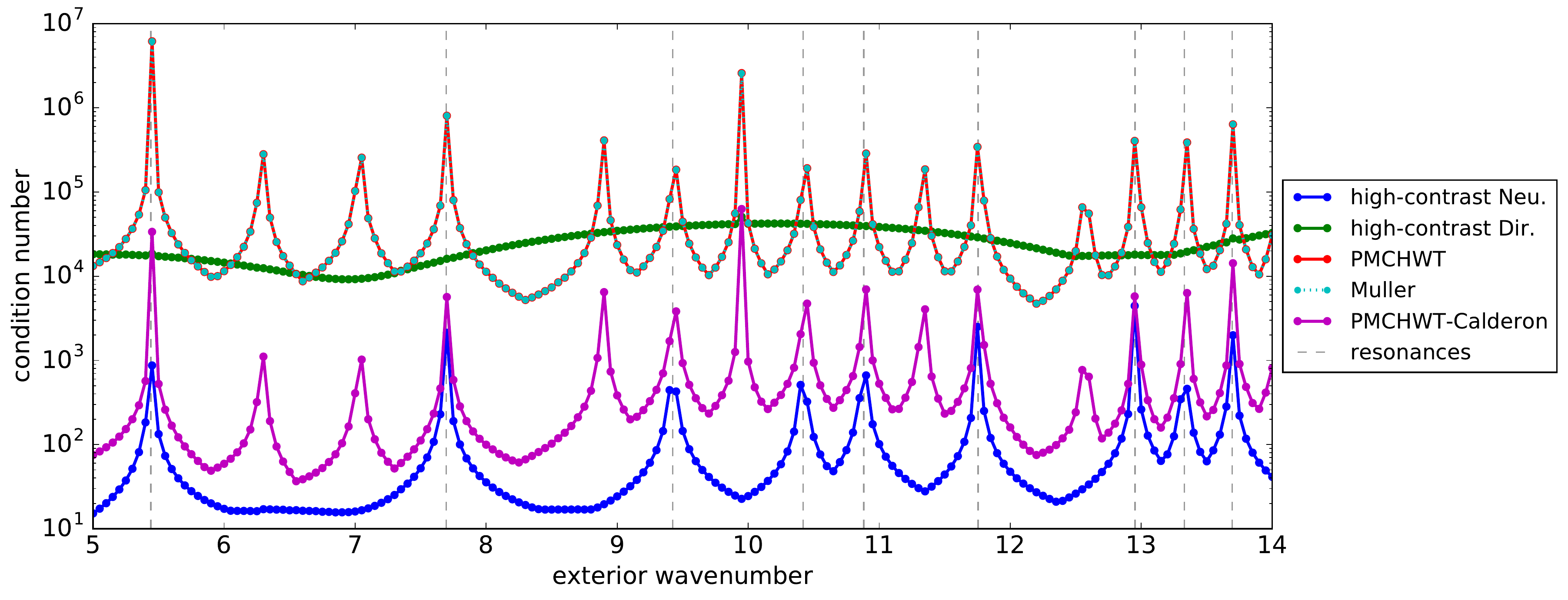}
		\caption{Exterior air and interior iron ($\kint/\kext=0.083$, $\rhoint/\rhoext=6306$).}
	\end{subfigure}
	\caption{The condition number with respect to the frequency of the wave field. The geometry is a unit cube with a mesh density of 6 triangles per wavelength.}
	\label{fig:freq:cond}
\end{figure}

\begin{figure}[!ht]
	\centering
	\begin{subfigure}[b]{0.99\columnwidth}
		\includegraphics[width=\textwidth]{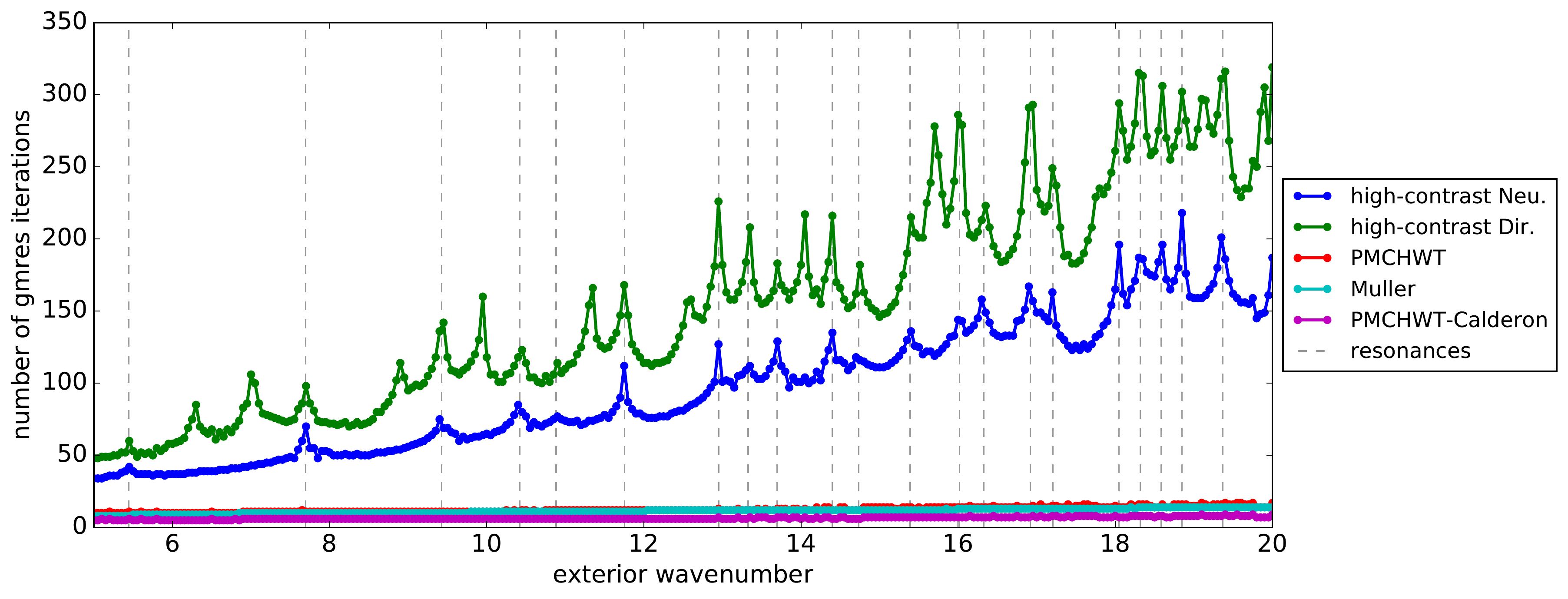}
		\caption{Exterior water and interior fat ($\kint/\kext=1.06$, $\rhoint/\rhoext=0.895$).}
	\end{subfigure}
	\begin{subfigure}[b]{0.99\columnwidth}
		\includegraphics[width=\textwidth]{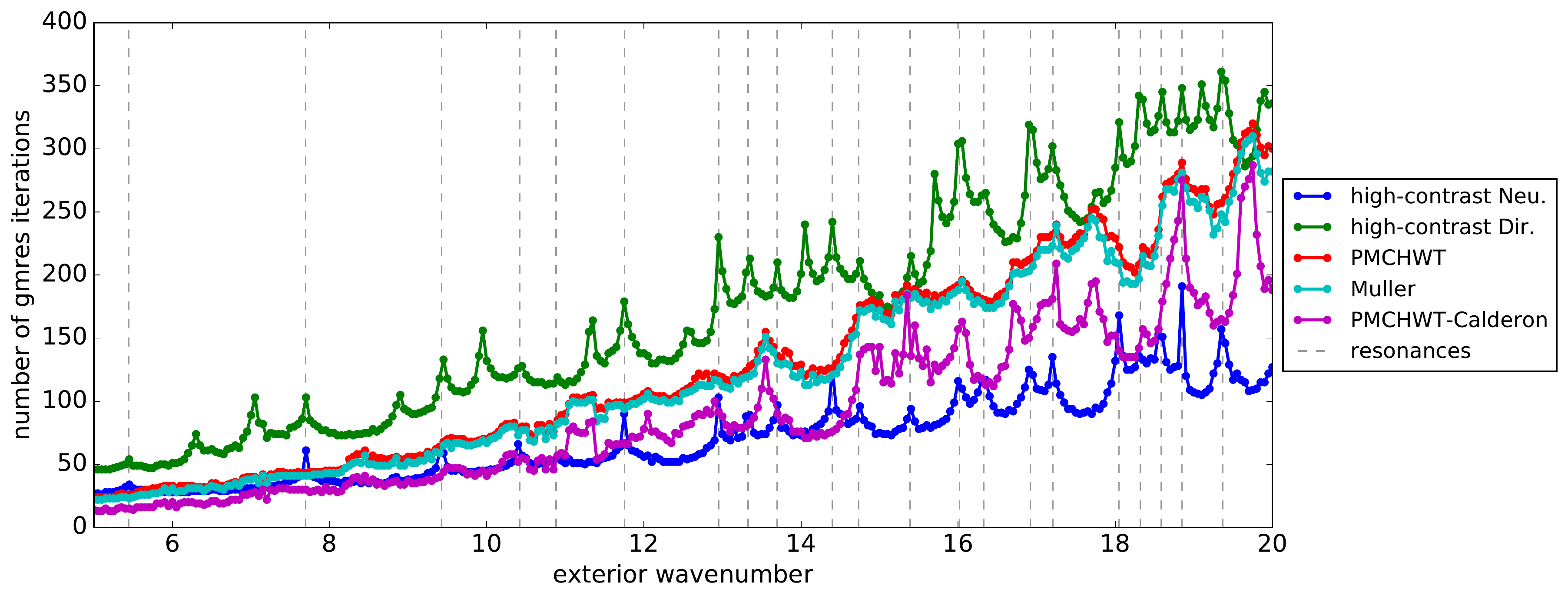}
		\caption{Exterior water and interior bone ($\kint/\kext=0.368$, $\rhoint/\rhoext=1.87$).}
	\end{subfigure}
	\begin{subfigure}[b]{0.99\columnwidth}
		\includegraphics[width=\textwidth]{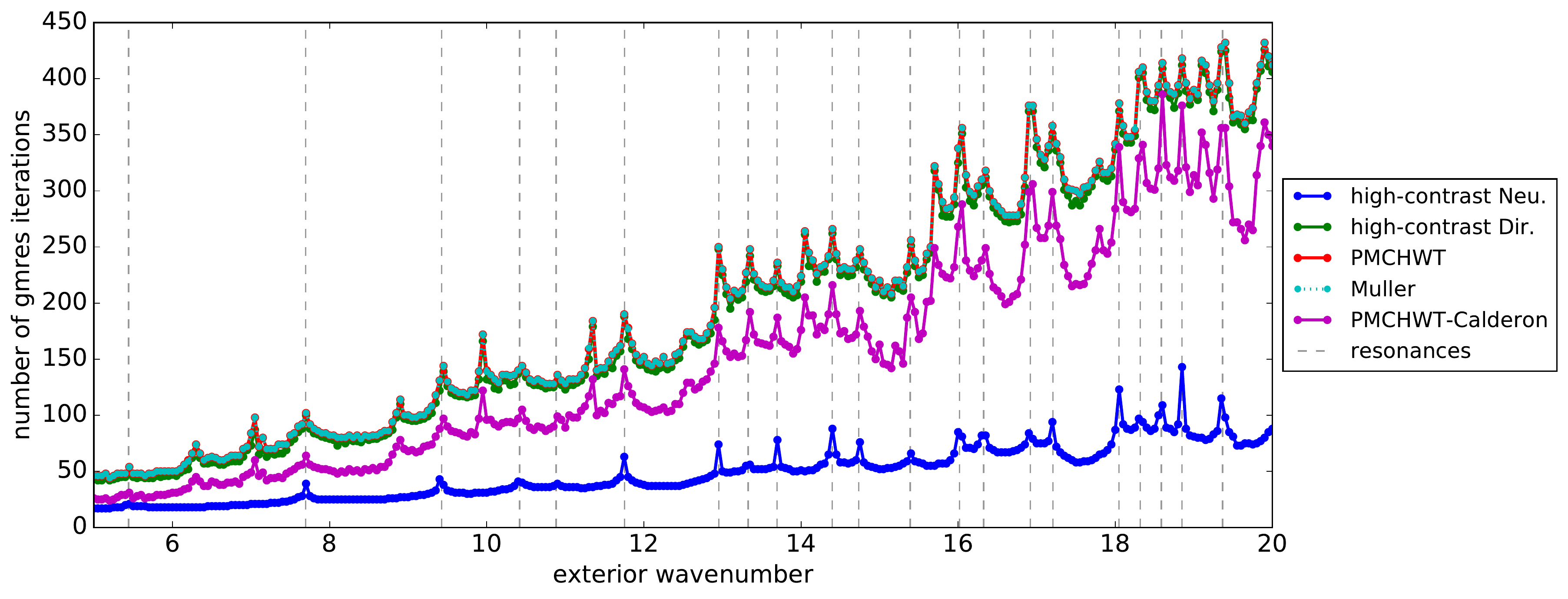}
		\caption{Exterior air and interior iron ($\kint/\kext=0.083$, $\rhoint/\rhoext=6306$).}
	\end{subfigure}
	\caption{The number of GMRES iterations with respect to the frequency of the wave field. The geometry is a unit cube with a mesh density of 6 triangles per wavelength.}
	\label{fig:freq:niter}
\end{figure}

The following benchmarks consider a unit cube, that is, the length of all edges is normalised to one. The incident wave field is a plane wave with direction vector $\begin{bmatrix} 1/\sqrt{3} & 1/\sqrt{3} & 1/\sqrt{3} \end{bmatrix}$. The triangular surface meshes are generated at each frequency, with at least 6 elements per wavelength. The smallest mesh has 202 vertices and the largest one 3068 vertices. Figures~\ref{fig:freq:cond} and~\ref{fig:freq:niter} present the condition number and the number of GMRES iterations, respectively. The vertical lines indicate the resonance frequencies of a rigid cube, given by
\begin{equation}
	k_n = \pi \sqrt{(n_x)^2 + (n_y)^2 + (n_z)^2} \qquad \text{for } n_x, n_y, n_z = 1,2,3,\dots.
\end{equation}
In the case of transmission problems, different resonances can occur as well~\cite{colton2010analytical, cossonniere2013surface}.

Comparing the condition number in Fig.~\ref{fig:freq:cond} with the GMRES convergence in Fig.~\ref{fig:freq:niter}, it is clear that the condition number of the system matrix is a good first estimate for the convergence behaviour of the GMRES algorithm. However, significant discrepancies between matrix conditioning and GMRES convergence can be observed as well. These differences are not surprising since the convergence of GMRES depends on the entire spectrum~\cite{antoine2021introduction} and might be quick for nearly singular systems~\cite{brown1997gmres}.

The different benchmarks at different material interfaces clearly demonstrate the influence of the material contrast on the efficiency of the boundary integral formulations. At low contrast in density and wavespeed, such as for the water-fat interface, the PMCHWT and Müller formulations are very efficient while the high-contrast formulations perform poorly. This is different when high material ratios are present in the configuration. For example, at the air-iron interface, the high-contrast Neumann formulation clearly outperforms all other formulations in terms of efficiency.

The frequency has a profound influence on the computational efficiency: the condition number increases and GMRES requires significant more iterations to converge. This deterioration in efficiency already starts at moderate frequencies and is on top of the increased number of degrees of freedom necessary to guarantee six elements per wavelength. This behaviour is well known and high-frequency simulations require specialised techniques from high-performance computing (cf.~\cite{betcke2017computationally}). The OSRC preconditioner can improve the convergence at high frequencies for first-kind formulations such as the PMCHWT~\cite{wout2021pmchwt} but cannot be applied directly to second-kind formulations such as the high-contrast and Müller formulations.

Another influence of the frequency on the conditioning of the system is the presence of spikes at specific frequencies. Boundary integral formulations are known to suffer from resonances in the acoustic transmission problem~\cite{hsiao2011system, hiptmair2021spurious}. In the case of rigid bodies, combined-field formulations such as the ones of Brakhage-Werner~\cite{brakhage1965dirichletsche} and Burton-Miller~\cite{burton1971application} can resolve spurious resonances~\cite{buffa2005regularized}. This approach has been extended to transmission problems (cf.~\cite{kleinman1988single, rapun2006indirect, claeys2015integral}) but a stability analysis at resonance frequencies is outside the scope of this study.

\FloatBarrier
\subsection{Multiple scattering}

All boundary integral formulations considered in this work can be extended to multiple scattering at disjoint simply-connected penetrable objects that are embedded in an exterior medium. Then, the system matrices have $2\ell \times 2\ell$ blocks of individual boundary integral operators, as in Eq.~\eqref{eq:system:multiplescattering} for the high-contrast formulations. Each of the objects can be composed of a different material and, therefore, material interfaces with a low contrast in material parameters as well as high-contrast interfaces can be present in the same configuration.

\begin{table}[!ht]
	\caption{Number of GMRES iterations for multiple scattering. All objects are spheres of radius 1 centered at location $(3j,0,0)$ for $j=0,1,2,\dots,\ell-1$. The incident plane wave field has direction $( 1/\sqrt{3} , 1/\sqrt{3} , 1/\sqrt{3} )$ and frequency $f=2353$ Hz. The entire mesh has 11\,270 nodes and at least 6 elements per wavelength. The exterior domain resembles water, and the interior materials either fat (\texttt{f}), bone (\texttt{b}), or iron (\texttt{i}). The wavenumbers are $k_\mathrm{water}=9.86$, $k_\mathrm{fat}=10.5$, $k_\mathrm{bone}=3.62$, and $k_\mathrm{iron}=3.61$.}
	\label{table:niter:multiplescattering}
	\centering
	\begin{tabular}{lrrrrr}
		\hline\hline
		interior & \multicolumn{1}{l}{high-contrast} & \multicolumn{1}{l}{high-contrast} & \multicolumn{1}{l}{PMCHWT} & \multicolumn{1}{l}{Müller} & \multicolumn{1}{l}{Calderón} \\
		materials & \multicolumn{1}{l}{Neumann} & \multicolumn{1}{l}{Dirichlet} & & & \multicolumn{1}{l}{PMCHWT} \\
		\hline
		\texttt{fffffff} & 95 & 442 & 36 & 36 & 19 \\
		\texttt{bbbbbbb} & 84 & 119 & 340 & 302 & 184 \\
		\texttt{iiiiiii} & 70 & 505 & 348 & 335 & 204 \\
		\texttt{fffbfff} & 105 & 648 & 98 & 96 & 53 \\
		\texttt{fffifff} & 102 & 666 & 110 & 103 & 67 \\
		\texttt{bbfffbb} & 119 & 845 & 258 & 232 & 142 \\
		\texttt{iifffii} & 115 & 985 & 292 & 282 & 183 \\
		\hline\hline
	\end{tabular}
\end{table}

Table~\ref{table:niter:multiplescattering} presents the number of GMRES iterations needed to solve the system matrix for a configuration of seven spheres. As expected, the standard formulations are very efficient at low-contrast problems involving water and fat only. When all interior materials resemble bone, the number of iterations increases for these formulations, which is even worse for the high-contrast water-iron interfaces. Differently, the high-contrast Neumann formulation is more efficient for high-contrast than for low-contrast transmission problems, which is consistent with the spectral analysis and the previous benchmarks.

The multiple scattering also allows for benchmarking the presence of interfaces with different material ratios. For example, the case `\texttt{fffbfff}' has six spheres made of fat and the middle one resembles bone. Even with only one high-contrast interface, the number of GMRES iterations already increases considerably for the standard formulations, in comparison with the case of fat materials only: almost three times the number of GMRES iterations are needed. The convergence of the high-contrast Neumann formulation deteriorates only slightly when different types of materials are present. This benchmark confirms that the high-contrast formulation can already have a superior efficiency when only few high-contrast interfaces are present in a multiple-scattering configuration. 

\subsection{Large-scale benchmark}

The final benchmark will be a large-scale simulation. As geometry, a \emph{mo'ai} statue~\cite{moai} will be used that is 2.42 meters tall, 1.05 meters wide and 83~cm thick. The acoustic parameters resemble basalt in the interior and air in the exterior medium. Hence, there is a large contrast in the wavespeed ($c_\mathrm{basalt}/c_\mathrm{air} = 9.85$) and density ($\rho_\mathrm{basalt}/\rho_\mathrm{air} = 2236.7$). The incident wave field is a plane wave with direction vector $\begin{bmatrix} 1/\sqrt{3} & 1/\sqrt{3} & -1/\sqrt{3} \end{bmatrix}$, unit amplitude, and a frequency of 3477~Hz ($\kext=64.3$ and $\kint=6.52$). The surface mesh has 35\,447 vertices and at least 7 elements per wavelength. The simulation was performed on a compute node with two 10-core Intel(R) Xeon(R) CPU E5-2630 v4 sockets, a clock speed of 2.4 GHz, hyperthreading activated (40 threads total), and a shared memory of 752~GB.

\begin{table}[!ht]
	\caption{The number of GMRES iterations and calculation (wall-clock) time for the large-scale benchmark. Time is indicated in the format h:mm:ss.}
	\label{table:largescale}
	\centering
	\begin{tabular}{lrrrr}
		\hline\hline
		formulation & \#iterations & \multicolumn{1}{l}{assembly} & solve & iteration \\
		\hline
		High-contrast Neumann   & 1852 & 36:06 &   21:05 & 0.68 s \\
		High-contrast Dirichlet & 4668 & 36:05 & 1:22:52 & 1.07 s \\
		PMCHWT                  & 4690 & 55:05 & 1:53:36 & 1.45 s \\
		Müller                  & 4690 & 55:05 & 2:01:43 & 1.56 s \\
		Calderón PMCHWT         & 4081 & 55:05 & 2:46:55 & 2.45 s \\
		\hline\hline
	\end{tabular}
\end{table}

\begin{figure}[!ht]
	\centering
	\begin{subfigure}[b]{0.55\columnwidth}
		\centering
		\includegraphics[width=\textwidth]{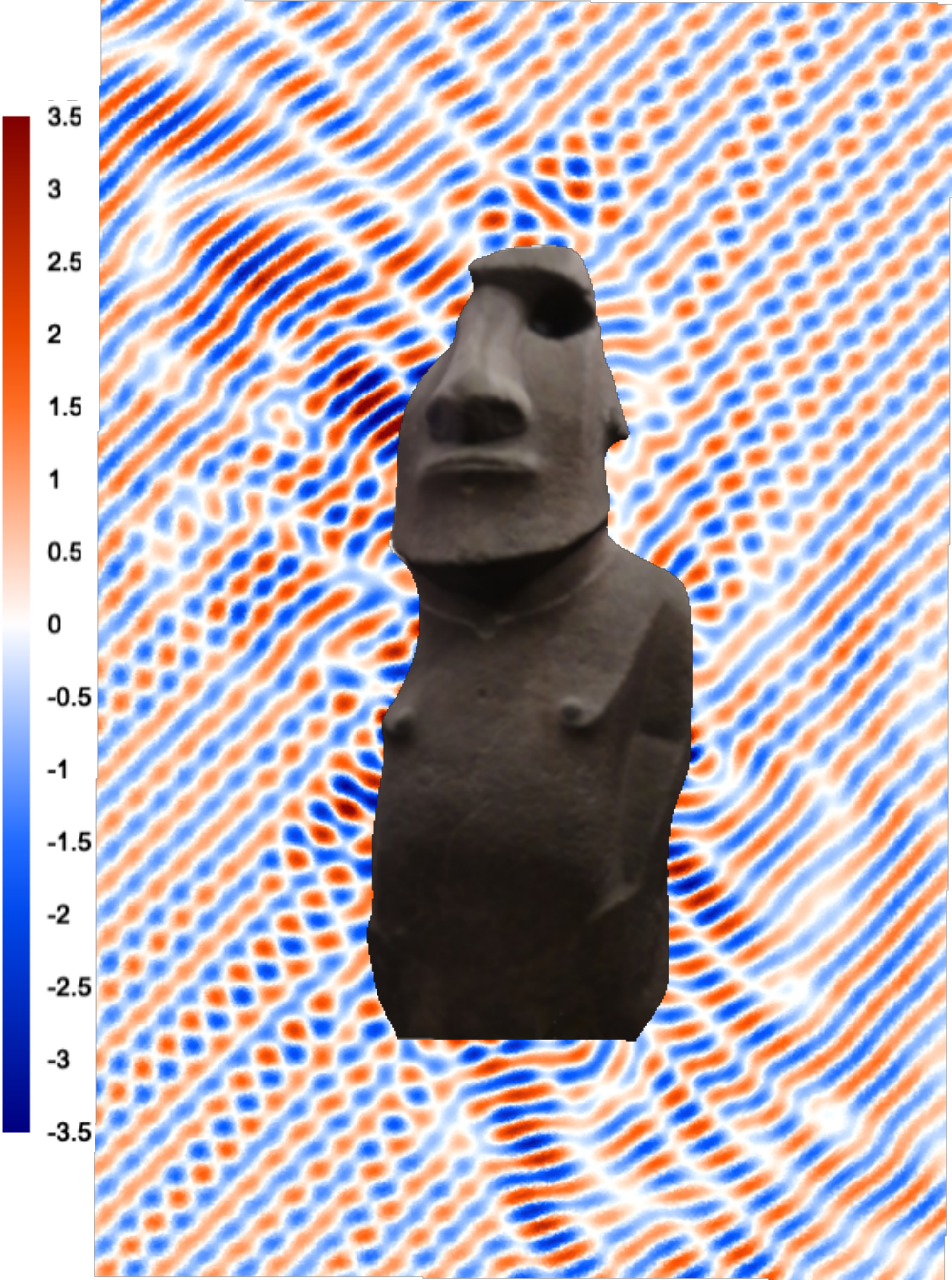}
		\caption{The real part of the acoustic pressure on a vertical plane in the exterior volume, calculated with the high-contrast Neumann formulation.}
	\end{subfigure}
	\hfill
	\begin{subfigure}[b]{0.43\columnwidth}
		\centering
		\includegraphics[width=.9\textwidth]{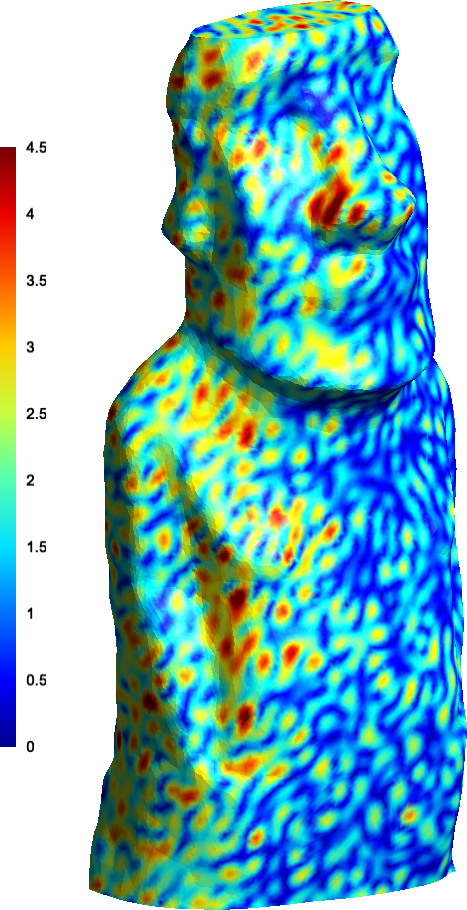}
		\caption{The magnitude of the acoustic pressure at the surface, calculated with the Calderón preconditioned PMCHWT formulation.}
	\end{subfigure}
	\caption{The acoustic field for the large-scale benchmark.}
	\label{fig:largescale}
\end{figure}

Table~\ref{table:largescale} presents the computational performance of the different formulations. The high-contrast Neumann formulation significantly outperforms all others on this high-contrast benchmark. The time to assemble the system matrix is only two thirds of the time for the standard formulations. Notice that it is not half of the assembly time since the adjoint double-layer operator does not require assembly: it is the transpose of the double-layer operator. Furthermore, the assembly time depends on the type of operator since the computational performance of the hierarchical matrix compression depends on the regularity of the kernel, with the single-layer operator the quickest and the hypersingular operator the most expensive operator to assemble in compressed format. Assembling the Calderón preconditioner does not incur any overhead since it equals the system matrix.

Concerning the time to solve the system, a large part the computation time of GMRES is spent on the matrix-vector multiplication of the compressed system matrices. As expected, the high-contrast formulations require only half the time per iteration compared to the standard formulations. Moreover, Calderón preconditioning doubles the calculation time of each GMRES iteration. While Calderón preconditioning improves GMRES convergence, this is not sufficient to improve the overall time to solve the system. The high-contrast Neumann formulation is by far the best conditioned system and requires, in total, only 57 minutes compared to almost 3 hours for the PMCHWT and Müller formulations. Figure~\ref{fig:largescale} depicts the acoustic field scattered from the penetrable structure.

\section{Conclusions}

This study analysed the influence of the mass density on the BEM's efficiency to solve acoustic transmission problems. High contrast in material parameters (density and wavespeed) between the homogeneous bounded object and the exterior medium causes ill-conditioning of the system matrix. The convergence of GMRES deteriorates quickly, and the BEM requires long computation times.

Novel high-contrast formulations were designed using a mix of a direct representation formula for the interior and an indirect representation for the exterior fields, with transmission modelled through Neumann-to-Dirichlet maps. Eight different boundary integral formulations can be devised with this approach, of which two are of the second kind and thus well-conditioned. The high-contrast Neumann and Dirichlet formulations' performance was compared against the PMCHWT, Müller and Calderón-preconditioned PMCHWT formulations.

A spectral analysis of the boundary integral formulations resulted in explicit expressions of the eigenvalue accumulation points in terms of the materials' density. The high-contrast Neumann and Dirichlet functions have two eigenvalue accumulation points each: $\tfrac12 \pm \tfrac\imath2 \sqrt{\rhoext/\rhoint}$ and $\tfrac12 \pm \tfrac\imath2 \sqrt{\rhoint/\rhoext}$, respectively. Hence, eigenvalues stay away from zero for any density ratio. For the high-contrast Neumann formulation, the accumulation points converge to a single point when the interior density is much higher than the exterior. For the Dirichlet version, this is the case when the density in the interior is much lower than in the exterior material. Furthermore, the Calderón preconditioned PMCHWT has a spectrum that accumulates around $\tfrac12 + \tfrac{\rhoext}{4\rhoint} + \tfrac{\rhoint}{4\rhoext}$ in the complex plane. These theoretical results were numerically validated on canonical test cases.

Extensive numerical benchmarks were performed for various materials, analysing the spectrum, condition number, and GMRES convergence. Small-scale simulations on a sphere show increasing ill-conditioning of the standard PMCHWT and Müller formulations at large density ratios while the high-contrast formulations and Calderón-preconditioned PMCHWT remain well-conditioned. For intermediate-scale simulations on a cube, the GMRES convergence slows down with increasing frequency. The traditional formulations are most efficient for low-contrast materials such as water and fat ($\rhoint/\rhoext = 0.895$). All formulations behave similarly for intermediate-contrast materials such as water and bone ($\rhoint/\rhoext = 1.87$). When considering high-contrast media such as air and iron ($\rhoint/\rhoext = 6306$), the high-contrast Neumann formulation significantly outperforms all other formulations. A benchmark simulation on seven spheres confirms that the high-contrast Neumann formulation can already outperform traditional boundary integral formulations when few high-contrast interfaces are present in a multiple-scattering configuration. An additional advantage of the high-contrast formulations is that they need less boundary integral operators, yielding faster matrix algebra and reducing memory consumption. Finally, a large-scale benchmark at a \emph{mo'ai} statue shows a reduction of hours in computation time when the novel high-contrast formulations simulate acoustic propagation in high-contrast media.

\section*{Acknowledgment}

This work was financially supported by CONICYT [FONDECYT 11160462], the Vicerrectoría de Investigación of the Pontificia Universidad Católica de Chile, and the EPSRC [EP/P012434/1].

\bibliographystyle{unsrt}
\bibliography{refs}

\begin{thebibliography}{10}

\bibitem{nedelec2001acoustic}
Jean-Claude N{\'e}d{\'e}lec.
\newblock {\em Acoustic and electromagnetic equations: integral representations
  for harmonic problems}, volume 144 of {\em Applied Mathematical Sciences}.
\newblock Springer, New York, 2001.

\bibitem{steinbach2008numerical}
Olaf Steinbach.
\newblock {\em Numerical approximation methods for elliptic boundary value
  problems: finite and boundary elements}.
\newblock Springer, New York, 2008.

\bibitem{sauter2010boundary}
Stefan~A Sauter and Christoph Schwab.
\newblock {\em Boundary Element Methods}, volume~39 of {\em Springer Series in
  Computational Mathematics}.
\newblock Springer, Berlin, 2010.

\bibitem{chew2001fast}
Weng~Cho Chew, Eric Michielssen, JM~Song, and Jian-Ming Jin.
\newblock {\em Fast and efficient algorithms in computational
  electromagnetics}.
\newblock Artech House, Inc., Norwood, MA, 2001.

\bibitem{betcke2017computationally}
Timo Betcke, Elwin van~'t Wout, and Pierre G{\'e}lat.
\newblock Computationally efficient boundary element methods for high-frequency
  {H}elmholtz problems in unbounded domains.
\newblock In Domenico Lahaye, Jok Tang, and Kees Vuik, editors, {\em Modern
  Solvers for {H}elmholtz Problems}, Geosystems Mathematics, pages 215--243.
  Birkh\"auser, Cham, 2017.

\bibitem{wout2021benchmarking}
Elwin van~'t Wout, Seyyed~R. Haqshenas, Pierre G{\'e}lat, Timo Betcke, and
  Nader Saffari.
\newblock Benchmarking preconditioned boundary integral formulations for
  acoustics.
\newblock {\em International Journal for Numerical Methods in Engineering},
  122(20):5873--5897, 2021.

\bibitem{poggio1973integral}
A.J. Poggio and E.K. Miller.
\newblock Integral equation solutions of three-dimensional scattering problems.
\newblock In R.~Mittra, editor, {\em Computer Techniques for Electromagnetics},
  International Series of Monographs in Electrical Engineering, chapter~4,
  pages 159--264. Pergamon, Oxford, UK, 1973.

\bibitem{chang1974surface}
Yu~Chang and Roger~F Harrington.
\newblock A surface formulation for characteristic modes of material bodies.
\newblock Technical report, Syracuse University, Syracuse, NY, 10 1974.
\newblock Technical Report TR-74-7.

\bibitem{wu1977scattering-bor}
Te-Kao Wu and Leonard~L Tsai.
\newblock Scattering from arbitrarily-shaped lossy dielectric bodies of
  revolution.
\newblock {\em Radio Science}, 12(5):709--718, 1977.

\bibitem{muller1957grundprobleme}
Claus M{\"u}ller.
\newblock {\em Grundprobleme der mathematischen {T}heorie elektromagnetischer
  {S}chwingungen}.
\newblock Springer, Berlin, 1957.

\bibitem{kress1978transmission}
R~Kress and GF~Roach.
\newblock Transmission problems for the {H}elmholtz equation.
\newblock {\em Journal of Mathematical Physics}, 19(6):1433--1437, 1978.

\bibitem{costabel1985direct}
Martin Costabel and Ernst Stephan.
\newblock A direct boundary integral equation method for transmission problems.
\newblock {\em Journal of Mathematical Analysis and Applications},
  106(2):367--413, 1985.

\bibitem{niino2012preconditioning}
Kazuki Niino and Naoshi Nishimura.
\newblock Preconditioning based on {C}alder{\'o}n's formulae for periodic fast
  multipole methods for {H}elmholtz' equation.
\newblock {\em Journal of Computational Physics}, 231(1):66--81, 2012.

\bibitem{haqshenas2021fast}
S.~R. Haqshenas, P.~G{\'e}lat, E.~van~'t Wout, T.~Betcke, and N.~Saffari.
\newblock A fast full-wave solver for calculating ultrasound propagation in the
  body.
\newblock {\em Ultrasonics}, 110:106240, 2021.

\bibitem{wout2021proximity}
Elwin van~'t Wout and Christopher Feuillade.
\newblock Proximity resonances of water-entrained air bubbles near acoustically
  reflecting boundaries.
\newblock {\em The Journal of the Acoustical Society of America},
  149(4):2477--–2491, 2021.

\bibitem{cummer2016controlling}
Steven~A Cummer, Johan Christensen, and Andrea Al{\`u}.
\newblock Controlling sound with acoustic metamaterials.
\newblock {\em Nature Reviews Materials}, 1(3):1--13, 2016.

\bibitem{gossye2018calderon}
Michiel Gossye, Martijn Huynen, Dries Vande~Ginste, Dani{\"e}l De~Zutter, and
  Hendrik Rogier.
\newblock A {C}alder{\'o}n preconditioner for high dielectric contrast media.
\newblock {\em IEEE Transactions on Antennas and Propagation}, 66(2):808--818,
  2018.

\bibitem{gossye2019electromagnetic}
Michiel Gossye, Dries Vande~Ginste, and Hendrik Rogier.
\newblock Electromagnetic modeling of high magnetic contrast media using
  {C}alder{\'o}n preconditioning.
\newblock {\em Computers \& Mathematics with Applications}, 77(6):1626--1638,
  2019.

\bibitem{antoine2008integral}
Xavier Antoine and Yassine Boubendir.
\newblock An integral preconditioner for solving the two-dimensional scattering
  transmission problem using integral equations.
\newblock {\em International Journal of Computer Mathematics},
  85(10):1473--1490, 2008.

\bibitem{yan2010comparative}
Su~Yan, Jian-Ming Jin, and Zaiping Nie.
\newblock A comparative study of {C}alder{\'o}n preconditioners for {PMCHWT}
  equations.
\newblock {\em IEEE Transactions on Antennas and Propagation},
  58(7):2375--2383, 2010.

\bibitem{cools2011calderon}
Kristof Cools, Francesco~P Andriulli, and Eric Michielssen.
\newblock A {C}alder{\'o}n multiplicative preconditioner for the {PMCHWT}
  integral equation.
\newblock {\em IEEE Transactions on Antennas and Propagation}, 59(12):4579,
  2011.

\bibitem{niino2012calderon}
K~Niino and N~Nishimura.
\newblock Calder{\'o}n preconditioning approaches for {PMCHWT} formulations for
  {M}axwell's equations.
\newblock {\em International Journal of Numerical Modelling: Electronic
  Networks, Devices and Fields}, 25(5-6):558--572, 2012.

\bibitem{betcke2020product}
Timo Betcke, Matthew~W Scroggs, and Wojciech {\'S}migaj.
\newblock Product algebras for {G}alerkin discretisations of boundary integral
  operators and their applications.
\newblock {\em ACM Transactions on Mathematical Software (TOMS)}, 46(1):1--22,
  2020.

\bibitem{marburg2003performance}
Steffen Marburg and Stefan Schneider.
\newblock Performance of iterative solvers for acoustic problems. {P}art {I}.
  {S}olvers and effect of diagonal preconditioning.
\newblock {\em Engineering Analysis with Boundary Elements}, 27(7):727--750,
  2003.

\bibitem{sakuma2008fast}
Tetsuya Sakuma, Stefan Schneider, and Yosuke Yasuda.
\newblock Fast solution methods.
\newblock In {\em Computational Acoustics of Noise Propagation in Fluids-Finite
  and Boundary Element Methods}, pages 333--366. Springer, Berlin, 2008.

\bibitem{saad1986gmres}
Youcef Saad and Martin~H Schultz.
\newblock {GMRES}: A generalized minimal residual algorithm for solving
  nonsymmetric linear systems.
\newblock {\em SIAM Journal on Scientific and Statistical Computing},
  7(3):856--869, 1986.

\bibitem{colton2013integral}
David Colton and Rainer Kress.
\newblock {\em Integral equation methods in scattering theory}.
\newblock SIAM, Philadelphia, PA, 2013.

\bibitem{costabel1988boundary}
Martin Costabel.
\newblock Boundary integral operators on {L}ipschitz domains: elementary
  results.
\newblock {\em SIAM Journal on Mathematical Analysis}, 19(3):613--626, 1988.

\bibitem{kleinman1988single}
RE~Kleinman and PA~Martin.
\newblock On single integral equations for the transmission problem of
  acoustics.
\newblock {\em SIAM Journal on Applied Mathematics}, 48(2):307--325, 1988.

\bibitem{claeys2013multi}
Xavier Claeys and Ralf Hiptmair.
\newblock Multi-trace boundary integral formulation for acoustic scattering by
  composite structures.
\newblock {\em Communications on Pure and Applied Mathematics},
  66(8):1163--1201, 2013.

\bibitem{boubendir2015integral}
Yassine Boubendir, Oscar Bruno, David Levadoux, and Catalin Turc.
\newblock Integral equations requiring small numbers of {K}rylov-subspace
  iterations for two-dimensional smooth penetrable scattering problems.
\newblock {\em Applied Numerical Mathematics}, 95:82--98, 2015.

\bibitem{lu2009phononic}
Ming-Hui Lu, Liang Feng, and Yan-Feng Chen.
\newblock Phononic crystals and acoustic metamaterials.
\newblock {\em Materials Today}, 12(12):34--42, 2009.

\bibitem{smigaj2015solving}
Wojciech {\'S}migaj, Timo Betcke, Simon Arridge, Joel Phillips, and Martin
  Schweiger.
\newblock Solving boundary integral problems with {BEM++}.
\newblock {\em ACM Transactions on Mathematical Software (TOMS)}, 41(2):6,
  2015.

\bibitem{scroggs2017software}
Matthew~W Scroggs, Timo Betcke, Erik Burman, Wojciech {\'S}migaj, and Elwin
  van~'t Wout.
\newblock Software frameworks for integral equations in electromagnetic
  scattering based on {C}alder{\'o}n identities.
\newblock {\em Computers \& Mathematics with Applications}, 74(11):2897--2914,
  2017.

\bibitem{geuzaine2009gmsh}
Christophe Geuzaine and Jean-Fran{\c{c}}ois Remacle.
\newblock Gmsh: A {3-D} finite element mesh generator with built-in pre- and
  post-processing facilities.
\newblock {\em International Journal for Numerical Methods in Engineering},
  79(11):1309--1331, 2009.

\bibitem{meshlab}
Paolo Cignoni, Marco Callieri, Massimiliano Corsini, Matteo Dellepiane, Fabio
  Ganovelli, and Guido Ranzuglia.
\newblock {MeshLab: an Open-Source Mesh Processing Tool}.
\newblock In Vittorio Scarano, Rosario~De Chiara, and Ugo Erra, editors, {\em
  Eurographics Italian Chapter Conference}. The Eurographics Association, 2008.

\bibitem{virtanen2020scipy}
Pauli Virtanen, Ralf Gommers, Travis~E. Oliphant, Matt Haberland, Tyler Reddy,
  David Cournapeau, Evgeni Burovski, Pearu Peterson, Warren Weckesser, Jonathan
  Bright, St{\'e}fan~J. {van der Walt}, Matthew Brett, Joshua Wilson, K.~Jarrod
  Millman, Nikolay Mayorov, Andrew R.~J. Nelson, Eric Jones, Robert Kern, Eric
  Larson, C~J Carey, {\.I}lhan Polat, Yu~Feng, Eric~W. Moore, Jake
  {VanderPlas}, Denis Laxalde, Josef Perktold, Robert Cimrman, Ian Henriksen,
  E.~A. Quintero, Charles~R. Harris, Anne~M. Archibald, Ant{\^o}nio~H. Ribeiro,
  Fabian Pedregosa, Paul {van Mulbregt}, and {SciPy 1.0 Contributors}.
\newblock {{SciPy} 1.0: Fundamental Algorithms for Scientific Computing in
  Python}.
\newblock {\em Nature Methods}, 17:261--272, 2020.

\bibitem{duck1990physical}
FA~Duck.
\newblock {\em Physical properties of tissue: a comprehensive reference book}.
\newblock Academic Press, London, UK, 1990.

\bibitem{itis2018}
{IT'IS Foundation}.
\newblock Tissue properties database, 2018.

\bibitem{gray1963american}
D.E. Gray.
\newblock {\em American Institute of Physics Handbook}.
\newblock McGraw-Hill, New York, NY, 1963.

\bibitem{wohletz1992volcanology}
Kenneth Wohletz and Grant Heiken.
\newblock {\em Volcanology and geothermal energy}, volume 432.
\newblock University of California Press, Berkeley, CA, 1992.

\bibitem{morse1986theoretical}
Philip~McCord Morse and K~Uno Ingard.
\newblock {\em Theoretical acoustics}.
\newblock Princeton University Press, Princeton, NJ, 1986.

\bibitem{colton2010analytical}
David Colton, Peter Monk, and Jiguang Sun.
\newblock Analytical and computational methods for transmission eigenvalues.
\newblock {\em Inverse Problems}, 26(4):045011, 2010.

\bibitem{cossonniere2013surface}
Anne Cossonni{\`e}re and Houssem Haddar.
\newblock Surface integral formulation of the interior transmission problem.
\newblock {\em Journal of Integral Equations and Applications}, 25(3):341--376,
  2013.

\bibitem{antoine2021introduction}
Xavier Antoine and Marion Darbas.
\newblock An introduction to operator preconditioning for the fast iterative
  integral equation solution of time-harmonic scattering problems.
\newblock {\em Multiscale Science and Engineering}, 3:1--35, 2021.

\bibitem{brown1997gmres}
Peter~N Brown and Homer~F Walker.
\newblock {GMRES} on (nearly) singular systems.
\newblock {\em SIAM Journal on Matrix Analysis and Applications}, 18(1):37--51,
  1997.

\bibitem{wout2021pmchwt}
E.~van~'t Wout, S.R. Haqshenas, P.~G{\'e}lat, T.~Betcke, and N.~Saffari.
\newblock Frequency-robust preconditioning of boundary integral equations for
  acoustic transmission.
\newblock Preprint available on Arxiv:2104.04609, 2021.

\bibitem{hsiao2011system}
George~C Hsiao and Liwei Xu.
\newblock A system of boundary integral equations for the transmission problem
  in acoustics.
\newblock {\em Applied Numerical Mathematics}, 61(9):1017--1029, 2011.

\bibitem{hiptmair2021spurious}
Ralf Hiptmair, Andrea Moiola, and Euan~A Spence.
\newblock Spurious quasi-resonances in boundary integral equations for the
  {H}elmholtz transmission problem.
\newblock {\em arXiv preprint:2109.08530}, 2021.

\bibitem{brakhage1965dirichletsche}
Helmut Brakhage and Peter Werner.
\newblock {\"U}ber das {D}irichletsche {A}u{\ss}enraumproblem f{\"u}r die
  {H}elmholtzsche {S}chwingungsgleichung.
\newblock {\em Archiv der Mathematik}, 16(1):325--329, 1965.

\bibitem{burton1971application}
AJ~Burton and GF~Miller.
\newblock The application of integral equation methods to the numerical
  solution of some exterior boundary-value problems.
\newblock {\em Proceedings of the Royal Society of London. A. Mathematical and
  Physical Sciences}, 323(1553):201--210, 1971.

\bibitem{buffa2005regularized}
Annalisa Buffa and Ralf Hiptmair.
\newblock Regularized combined field integral equations.
\newblock {\em Numerische Mathematik}, 100(1):1--19, 2005.

\bibitem{rapun2006indirect}
Mar{\'\i}a-Luisa Rap{\'u}n and Francisco-Javier Sayas.
\newblock Indirect methods with {B}rakhage-{W}erner potentials for {H}elmholtz
  transmission problems.
\newblock In {\em Numerical Mathematics and Advanced Applications}, pages
  1146--1154. Springer, 2006.

\bibitem{claeys2015integral}
Xavier Claeys and Ralf Hiptmair.
\newblock Integral equations for acoustic scattering by partially impenetrable
  composite objects.
\newblock {\em Integral Equations and Operator Theory}, 81(2):151--189, 2015.

\bibitem{moai}
The~British Museum.
\newblock Hoa {H}akananai'a, 2014.

\end{thebibliography}

\end{document}